\documentclass[11pt]{article}
\usepackage[margin=1in]{geometry}
\usepackage[T1]{fontenc}
\usepackage{lmodern}
\usepackage{amsmath,amssymb,amsthm,mathtools,bm}
\IfFileExists{bbm.sty}{%
	\usepackage{bbm}\newcommand{\ind}{\mathbbm 1}}{%
	\newcommand{\ind}{\boldsymbol 1}}
\usepackage[numbers,square,sort&compress]{natbib}
\usepackage[hidelinks]{hyperref}
\usepackage{enumitem,microtype}
\usepackage{booktabs}
\usepackage{graphicx}
\hypersetup{pdftitle={Peak intensities of general Gaussian fields on Riemannian surfaces: trace reduction, singular jet laws, and stability},pdfauthor={Dan Cheng}}
\numberwithin{equation}{section}
\setlist{itemsep=2pt,topsep=4pt,leftmargin=2em}
\newtheorem{theorem}{Theorem}[section]
\newtheorem{proposition}[theorem]{Proposition}
\newtheorem{lemma}[theorem]{Lemma}
\newtheorem{corollary}[theorem]{Corollary}
\theoremstyle{definition}

\newtheorem{assumption}[theorem]{Assumption}
\theoremstyle{remark}
\newtheorem{remark}[theorem]{Remark}

\newcommand{\R}{\mathbb R}
\newcommand{\E}{\mathbb E}
\newcommand{\PP}{\mathbb P}
\newcommand{\dd}{\,\mathrm d}
\newcommand{\Cov}{\operatorname{Cov}}
\newcommand{\Var}{\operatorname{Var}}
\newcommand{\tr}{\operatorname{tr}}
\newcommand{\diag}{\operatorname{diag}}

\newcommand{\Acal}{\mathcal A}

\newcommand{\Kcal}{\mathcal K}
\newcommand{\phin}[3]{\varphi(#1;#2,#3)}
\newcommand{\norm}[1]{\lVert#1\rVert}
\DeclareMathOperator{\vech}{vech}
\title{Peak intensities of general Gaussian fields on Riemannian
	surfaces: trace reduction, singular jet laws, and stability}
\author{Dan Cheng\thanks{School of Mathematical and Statistical Sciences,
		Arizona State University, Tempe, Arizona, USA.}
	\thanks{Corresponding author. Email: \texttt{chengdan@asu.edu}.}
	\and Zhou Yu\thanks{Department of Geography, University of Colorado Boulder,
		Boulder, Colorado, USA.}
	\and Guofeng Cao\footnotemark[3]}
\date{}
\begin{document}
	\maketitle
	
	\begin{abstract}
		We compute peak intensities, expected numbers of peaks, and local peak height distributions for smooth Gaussian fields
		on a Riemannian surface with general mean and covariance functions, and prove that local peak height laws are stable across the singular limits of the covariance structure. Splitting the covariant Hessian into trace and traceless parts reduces the Kac--Rice conditional determinant weight to a Gaussian trace average of a bivariate Gaussian norm moment, which is evaluated through noncentral chi square distribution functions. For a nondegenerate Hessian law, diagonalizing the determinant makes a spectral trace independent of the norm and reduces the weight to a single integral against one fixed distribution function. Assuming that height and gradient, and gradient and
		Hessian, are separately nondegenerate, we classify all singular laws of
		the jet formed by height, gradient and Hessian. For a singular jet law, given the height, the Hessian lies on an affine plane; the signature of its constraint operator makes the peak height support a half line or the whole real line, and its type (definite, rank-one or indefinite) determines the canonical weight.
		A singular jet law at a point means exactly that the field satisfies a linear second-order relation there almost surely, elliptic, parabolic or hyperbolic according to that type. Local peak height laws are stable under perturbations of the mean and covariance jets, even across degeneracy: they are Lipschitz continuous in total variation, with explicit rates for expected regional counts, and H\"older continuous of the optimal order $3/4$ in the uniform norm. Examples in the plane, on the sphere and on an ellipsoid illustrate the formulas in nondegenerate cases and realize all three constraint types in singular ones.
	\end{abstract}
	
	
	
	\section{Introduction}
	\label{sec:intro}
	Local maxima of smooth random fields are used to detect and test
	localized signals in images, spatial data, cosmological maps and rough
	surfaces
	\citep{TaylorWorsley2007,SchwartzmanGavrilovAdler2011,ChengSchwartzman2017,BardeenEtAl1986,BondEfstathiou1987,FuksCharnotskii2006}.
	Such inference requires the expected number of local maxima above a
	level and the distribution of the height of a peak at a given location.
	For a smooth Gaussian field $X$ on a Riemannian surface, the Kac--Rice
	formula expresses both through the conditional peak weight, also called
	the conditional determinant weight,
	\[
	J(t,x)=\E\big[|\det\nabla^2X(t)|\ind_{\{\nabla^2X(t)\prec0\}}\mid
	X(t)=x,\nabla X(t)=0\big]
	\]
	and its height average \citep{AdlerTaylor2007,AzaisWschebor2009,ChengSchwartzman2015}.
	This weight is explicit for isotropic fields in Euclidean space and on
	spheres \citep{ChengSchwartzman2015,ChengSchwartzman2018,AzaisDelmas2022,KurikiEtAl2026},
	and for stationary planar fields with trends under two
	uncorrelatedness assumptions \citep{Cheng2025}. For
	anisotropic stationary planar fields, Longuet-Higgins
	\citep{LonguetHiggins1957} obtained the total density of maxima, but
	height distributions only for narrow spectra. For nonstationary fields with a general mean, as produced for instance by
	smoothing data with a spatially varying covariance, $J(t,x)$ can be
	approximated by three-dimensional numerical integration over the
	Hessian law, but no reduction to low-dimensional integrals of standard
	distribution functions has been available.
	
	This paper computes $J(t,x)$, and hence the peak intensity, the expected
	number of peaks above a level and the local peak height density, for an
	arbitrary Gaussian law of the height, gradient and covariant Hessian at a
	point of a Riemannian surface. The starting point is classical: in
	dimension two the Hessian has eigenvalues $\tau/2\pm\norm\xi$, where
	$\tau$ is its trace and $\xi$ collects its two traceless coordinates; and
	closely related parametrizations by trace and ellipticity are standard in
	peak statistics
	\citep{BardeenEtAl1986,BondEfstathiou1987,PogosyanGayPichon2009}.
	Those works concern stationary isotropic fields, for which, given the
	height, rotation invariance makes the trace and the traceless part of the
	Hessian independent and the traceless part centered and circularly
	symmetric, so that in dimension two the conditional peak weight is a single Gaussian
	integral of the
	circular kernel~\eqref{eq:centeredcircularkernel}. For stationary planar fields with trends, Cheng~\citep{Cheng2025}
	derives computable formulas
	assuming uncorrelated gradient coordinates, rescaled to unit covariance,
	and a mixed Hessian entry uncorrelated with each diagonal entry. The present calculation allows arbitrary positive definite gradient
	covariance, general correlations among Hessian entries, and a general
	mean function, together with nonstationarity and nonconstant variance.
	Since $J(t,x)$ is defined
	intrinsically, every correct evaluation is independent of coordinates and
	frame; the difficulty lies in evaluating it for such Hessian laws, whose
	covariance may also be singular.
	
	Singular jet laws also make continuity a genuine question. When the
	height is an affine function of the gradient and Hessian and the
	constraint is semidefinite, the peak height support is a half line,
	whereas a nonsingular law gives a positive height density on the whole
	real line. Along nondegenerate models that approach such a singular one,
	densities that are positive everywhere can therefore converge only to a
	density that vanishes on a half line. At the same time the
	height-conditioned Hessian covariance loses rank, and the matrix square
	root, which underlies the natural coupling of Gaussian laws, is only
	$1/2$-H\"older continuous at singular matrices, so direct coupling bounds
	lose half the rate. Whether peak height laws, intensities and expected
	counts depend continuously on the mean and covariance across such
	limits, and at which rate, is therefore not evident. The question also matters in practice: these quantities are computed
	from estimated covariances, which need not share the singularity of the
	underlying model.
	
	The main contributions are as follows.
	\begin{enumerate}
		\item \textit{Low-dimensional representations of conditional peak weights for general Gaussian fields.}
		Theorem~\ref{thm:main} reduces $J(t,x)$ and its normalizer, for an
		arbitrary Gaussian Hessian law, to a Gaussian average over the trace
		of a lower partial moment of a two-dimensional Gaussian norm.
		Section~\ref{subsec:finiteCDF} and Appendix~\ref{app:normmoments}
		evaluate this kernel through ordinary or generalized noncentral
		chi square CDFs, including a closed form for the lower partial moment
		of centered norms with unequal variances
		(Proposition~\ref{prop:ordinarycentered}); at most two scalar
		integrations remain, and often one or none. For nondegenerate
		Hessian laws, Theorem~\ref{thm:spectral} diagonalizes the determinant,
		which makes a spectral trace independent of a norm with a fixed
		law and reduces the expected peak weight to a single integral against one
		distribution function. For centered fields,
		Corollary~\ref{cor:centerednormalizer} thereby expresses the normalizer, and hence
		the peak intensity and the expected total number of peaks, through the
		spectrum of the determinant quadratic form: a single integral of
		ordinary noncentral chi square CDFs, or a closed form when the two norm
		eigenvalues coincide.
		
		\item \textit{Classification of singular jet laws and peak height
			supports.} The law of the jet formed by $X(t)$, $\nabla X(t)$ and
		$\nabla^2X(t)$ is singular exactly when the height is an affine function of the gradient and
		Hessian; conditioning on the height then makes the Hessian covariance
		singular even though the Hessian given the gradient is nondegenerate.
		Theorem~\ref{thm:singularrank} shows that every such conditional law
		has rank two and is supported on an affine plane
		$\{\tr(QH)=x_*-x\}$. Theorem~\ref{thm:singularsupport} shows that the
		signature of the constraint operator $Q$ determines the peak height
		support: a half line when $Q$ is semidefinite and the whole real line
		when $Q$ is indefinite. In particular, a bounded peak height support
		cannot occur as long as the Hessian law given the gradient is
		nondegenerate. Theorem~\ref{thm:canonicalweights} gives canonical
		conditional peak weights for definite, rank-one and indefinite constraints.
		Such relations were known for Laplace eigenfunctions, whose trace
		is fixed given the height \citep[Proposition~3.10]{ChengSchwartzman2018},
		and as linear relations in the setting of \citep[Section~2.2]{Cheng2025};
		the classification here requires neither isotropy nor the two
		uncorrelatedness assumptions of \citep{Cheng2025}.
		Theorem~\ref{thm:singularpde} identifies a singular jet at a point
			with an almost sure linear second-order relation
		$Z+\tr(Q\nabla^2Z)+\langle\beta,\nabla Z\rangle=0$ there, which
			becomes an equation for the sample paths on open sets where singularity
			persists; its elliptic,
		parabolic or hyperbolic type is determined by the signature of $Q$, and
		peak heights are thus confined to a half line exactly where the
		equation is elliptic or parabolic. For stationary planar fields,
		singularity means that the spectral measure lives on an ellipse, two parallel lines or a hyperbola (Corollary~\ref{cor:spectralconic}).
		
		\item \textit{Stability of peak height laws across covariance degeneracy.}
		Conditioning on the Hessian instead of the height represents the local
		peak height as an affine function of a Hessian-tilted Gaussian vector
		plus independent Gaussian noise whose variance vanishes exactly in the
		singular case (Proposition~\ref{prop:tiltedrep}).
		Theorem~\ref{thm:singularcontinuity} uses this to show that local peak
		height laws are Lipschitz continuous in total variation with respect to
		the Gaussian regression parameters, uniformly on compact sets and across
		singular limits, and H\"older continuous of order $3/4$ in the
		uniform norm.
		Proposition~\ref{prop:contact} shows that the height density vanishes at
		the endpoint of a half-line support to order $4$ under a definite and to
		order $3/2$ under a rank-one constraint, and
		Proposition~\ref{prop:sharprate} uses the latter to show that the
		exponent $3/4$ is optimal.
		Corollary~\ref{cor:fieldparametercontinuity} gives the corresponding
		rates in terms of mean and covariance jets and for expected regional
		counts. 
	\end{enumerate}
	
	Explicit examples illustrate each part. A finite Fourier family realizes
	the three constraint types
	(Section~\ref{subsubsec:singularfamily}), and a nonstationary quadratic
	field passes through all of them as the location varies (Section~\ref{subsec:quadratic}). For a nonisotropic
	Gaussian quadrupole on the sphere, peak heights coincide with the
	largest eigenvalue of a random matrix, which yields an exact global check
	of the height-resolved local formulas
	(Section~\ref{subsec:sphericalchecks}). A quadratic form on an
	ellipsoid gives the same exact check on a surface of nonconstant
	curvature, where the second fundamental form enters the covariant
	Hessian (Section~\ref{subsec:ellipsoid}). Known isotropic laws and
	Longuet-Higgins' anisotropic planar counts are recovered as special cases
	and are stated only briefly. Second-order properties of critical points,
	such as correlations and fluctuations of counts
	\citep{BeliaevCammarotaWigman2019,BeliaevCammarotaWigman2020,Muirhead2020,CammarotaMarinucciWigman2016,CammarotaMarinucci2020},
	and critical points of high-dimensional landscapes
	\citep{Fyodorov2004,BrayDean2007} are not considered here.
	
	Sections~\ref{sec:framework} and~\ref{sec:regression} recall the
	Kac--Rice formulas and the Gaussian regression. Section~\ref{sec:trace}
	contains the trace reductions, Section~\ref{sec:singular} develops the singular
	classification, and Section~\ref{sec:continuity} establishes the
	stability of peak height laws across covariance degeneracy. Sections~\ref{sec:stationary}
	and~\ref{sec:sphere} treat the plane and curved surfaces (the sphere and an ellipsoid), and
	Section~\ref{sec:discussion} discusses higher dimensions.
	Appendix~\ref{app:normmoments} collects the norm moments.
	
	\section{Peak intensities}
	\label{sec:framework}
	
	\subsection{Model and assumptions}
	Let $(\mathcal M,\mathsf g)$ be a smooth, second countable Riemannian
	surface. Consider the Gaussian field
	\[
	X(t)=\mu(t)+Z(t),\qquad t\in\mathcal M,
	\]
	where $\mu$ is deterministic and $Z$ is centered, with covariance
	$C(t,t')=\E[Z(t)Z(t')]$. All gradients are defined with respect to $\mathsf g$, and all Hessians
	are covariant Hessians
	associated with the Levi--Civita connection.
	At each $t\in\mathcal M$, the covariant Hessian $(\nabla^2X)_t:T_t\mathcal M\times T_t\mathcal M\to\mathbb R$ is a symmetric bilinear form. Write $\dd A_{\mathsf g}$
	for Riemannian area and $d_{\mathsf g}$ for distance.
	
	On a coordinate neighborhood choose a smooth orthonormal tangent frame
	$(e_1,e_2)$. The coordinate vectors used in the Gaussian calculation are
	\[
	\nabla X(t)=(e_1X(t),e_2X(t))^\top,\qquad
	X_{jk}(t)=(\nabla^2X)_t(e_j,e_k),\qquad
	W=(X_{11}(t),X_{22}(t),X_{12}(t))^\top.
	\]
	Here $e_jX(t)$ denotes the directional derivative of $X$ along $e_j(t)$,
	and $X_{jk}(t)$ are the Hessian entries in the chosen frame. We also write $\nabla^2X(t)$ for the symmetric matrix $(X_{jk}(t))$ in
	this frame. Its trace is $\Delta_{\mathcal M}X(t)$, where the Laplace
	operator is the trace of the covariant Hessian. Its determinant is the
	determinant of the associated self-adjoint tangent endomorphism.
	For any symmetric
	$2\times2$ matrix $H$, set $\vech(H)=(H_{11},H_{22},H_{12})^\top$.
	The diagonal entries come first throughout the paper.
	We call the random vector $(X(t),\nabla X(t)^\top,W^\top)^\top\in\R^6$
	the jet of $X$ at $t$ and its Gaussian law the jet law; its mean and
	covariance, the mean and covariance jets, are given in
	Section~\ref{sec:regression}. A change of orthonormal frame acts on the
	jet by an invertible linear map. The jet law is called singular if its
	covariance matrix is singular.
	
	\begin{assumption}[Regularity and nondegeneracy]
		\label{ass:regularity}
		The mean $\mu$ belongs to $C^2(\mathcal M)$, and $Z$ has almost surely
		$C^2$ sample paths. In each smooth local orthonormal frame and on each
		compact subset $K$ of its domain, there are constants $M_K,\alpha_K>0$
		and $0< \delta_K<1$ such that, for $t,t'\in K$ and
		$0<d_{\mathsf g}(t,t')< \delta_K$,
		\[
		\max_{1\le i\le j\le2}
		\Var\{X_{ij}(t)-X_{ij}(t')\}
		\le M_K|\log d_{\mathsf g}(t,t')|^{-(1+\alpha_K)}.
		\]
		Moreover, for every $t\in\mathcal M$,
		\begin{equation}\label{eq:non-degenerate}
			\Cov\bigl((\nabla X(t)^\top,W^\top)^\top\bigr)\succ0,
			\qquad
			\Cov\bigl((X(t),\nabla X(t)^\top)^\top\bigr)\succ0.
		\end{equation}
	\end{assumption}
	
	The conditions are local and invariant under a change of
	orthonormal frame. Given the first condition in
	\eqref{eq:non-degenerate}, the second is equivalent to
	$\Var(X(t)\mid\nabla X(t))>0$; it is used only for statements at a
	fixed height. The increment bound and the first condition in
	\eqref{eq:non-degenerate} give the Morse property locally by the argument
	of \citep[Corollaries~11.3.2 and~11.3.5]{AdlerTaylor2007}: in a chart,
	ordinary and covariant jets are related by smooth invertible triangular
	maps that preserve the logarithmic bound, and a deterministic $C^2$ mean
	does not affect the density and modulus estimates used there. Thus every
	critical point has a nonsingular covariant Hessian almost surely. Almost
	sure $C^2$ paths and Fernique's theorem also make $C$ of class $C^{2,2}$,
	so the covariance derivatives below exist and are continuous. The
	assumptions concern height and gradient, and gradient and Hessian,
	separately; the jet law itself may be singular, and
	Theorem~\ref{thm:singularrank} shows that this happens exactly when the
	height is an affine function of the gradient and Hessian.
	
	Let $D\subset\mathcal M$ be a relatively compact Borel set and
	$B\subset\R$ be Borel. Write $N(D,B)$ for the number of local maxima
	of $X$ in $D$ with heights in $B$, $N(D)=N(D,\R)$, and
	$N(D,u)=N(D,(u,\infty))$. These are critical points of the field on the
	surface, including when $D$ has a boundary; they do not include maxima
	of the restriction to that boundary. We write $N_{\mathrm{min}}$ and
	$N_{\mathrm{sdl}}$ for the corresponding counts of minima and saddles. At a nondegenerate critical
	point these types correspond to a negative definite, positive definite,
	or indefinite covariant Hessian, respectively.
	
	\subsection{Kac--Rice formula for counts and heights of local maxima}\label{sec:Kac-Rice}
	Define the conditional peak weight (or conditional determinant weight), and its
	normalizer, by
	\begin{equation*}		\begin{split}
			J(t,x)&=\E\big[|\det\nabla^2X(t)|
			\ind_{\{\nabla^2X(t)\prec0\}}\mid X(t)=x,\nabla X(t)=0\big],\\
			J(t)&=\E\big[|\det\nabla^2X(t)|
			\ind_{\{\nabla^2X(t)\prec0\}}\mid\nabla X(t)=0\big].
		\end{split}
	\end{equation*}
	Here $\nabla X(t)=(e_1X(t),e_2X(t))^\top$ and
	$\nabla^2X(t)=(X_{ij}(t))_{1\le i,j\le2}$ denote, respectively, the vector of
	the Riemannian gradient and the matrix of the covariant Hessian
	in the chosen local orthonormal tangent frame. The condition
	$\nabla X(t)=0$, the determinant $\det\nabla^2X(t)$, and the
	negative definiteness of $\nabla^2X(t)$ are independent of the chosen frame.
	Conditional expectations given $\{X(t)=x,\nabla X(t)=0\}$ are the
	Gaussian regression versions of Section~\ref{sec:regression}, which are
	continuous in $(t,x)$; this is the version required by the Kac--Rice
	formula, and it remains meaningful when the joint law of height,
	gradient and Hessian is singular.
	Under Assumption~\ref{ass:regularity}, applying the Kac--Rice formula
	\citep[Theorems~6.2 and~6.4]{AzaisWschebor2009} in local charts gives
	its intrinsic form \citep{AdlerTaylor2007,ChengSchwartzman2015}:
	\begin{equation}\label{eq:KR}
		\E[N(D,B)]=\int_D\int_B\rho(t,x)\dd x\,\dd A_{\mathsf g}(t), \qquad \E[N(D)]=\int_D\ell(t)\,\dd A_{\mathsf g}(t),
	\end{equation}
	where
	\begin{equation}\label{eq:rho}
		\begin{split}
			\rho(t,x)&=p_{\nabla X(t)}(0)
			p_{X(t)\mid\nabla X(t)=0}(x)J(t,x),\\
			\ell(t)&=\int_\R\rho(t,x)\dd x
			=p_{\nabla X(t)}(0)J(t).
		\end{split}
	\end{equation}
	Here $p_{\nabla X(t)}(0)$ is the density of the two orthonormal
	gradient coordinates. Thus $\rho(t,x)$ is the peak intensity per unit
	Riemannian area and unit height, whereas $\ell(t)$ is the intensity
	per unit Riemannian area. These functions do not depend on the frame. In particular,
	\eqref{eq:KR} gives the expected number of peaks above level $u_0$:
	\begin{equation*}
		\E[N(D,u_0)]=\int_D\int_{u_0}^\infty\rho(t,x)\dd x\,\dd A_{\mathsf g}(t).
	\end{equation*}
	Normalizing by the spatial intensity gives the peak height
	density at $t$ \citep{ChengSchwartzman2015}:
	\begin{equation}\label{eq:peak-height}
		h_t(x) = \frac{\rho(t,x)}{\ell(t)} = \frac{J(t,x)}{J(t)}\cdot p_{X(t)\mid\nabla X(t)=0}(x), \qquad x\in \R.
	\end{equation}
	
	In a chart with metric matrix $\mathsf g_{\mathrm c}$, the coordinate
	gradient $\partial_qX$ and the covariant Hessian
	$(H_{\mathrm c})_{ij}=\partial_{ij}X-\sum_k\Gamma^k_{ij}\partial_kX$ satisfy
	$\nabla X=E_{\mathrm c}^\top\partial_qX$ and $H=E_{\mathrm c}^\top H_{\mathrm c}E_{\mathrm c}$, where the columns of $E_{\mathrm c}$ are the
	coordinate components of an orthonormal frame. Since
	$\det H=\det H_{\mathrm c}/\det\mathsf g_{\mathrm c}$,
	$p_{\nabla X}(0)=(\det\mathsf g_{\mathrm c})^{1/2}p_{\partial_qX}(0)$ and
	$\dd A_{\mathsf g}=(\det\mathsf g_{\mathrm c})^{1/2}\dd q$, the metric factors
	cancel in the count integrand, recovering the coordinate Kac--Rice
	formula; covariances of $H_{\mathrm c}$ require the connection terms.
	
	\section{Covariance derivatives and Gaussian regression}
	\label{sec:regression}
	
	\subsection{Intrinsic covariance derivatives}
	Fix $t$ and a local orthonormal tangent frame, and suppress $t$ when
	convenient. Write $v=\Var(X(t))$, $r=\Cov(\nabla X(t),X(t))$,
	$G=\Cov(\nabla X(t))$, $\iota=\Cov(W,X(t))$,
	$A=\Cov(W,\nabla X(t))$, and $V=\Cov(W)$. Then
	\begin{equation}\label{eq:Sigma}
		\Sigma(t)=\Cov\begin{pmatrix}X(t)\\\nabla X(t)\\W\end{pmatrix}
		=\begin{pmatrix}v&r^\top&\iota^\top\\r&G&A^\top\\\iota&A&V\end{pmatrix}.
	\end{equation}
	The corresponding means are $\mu=\mu(t)$, $\nabla\mu=\nabla\mu(t)$ and $m_W=\vech(\nabla^2\mu(t))$.
	We call $(\mu(t),\nabla\mu(t),m_W)$ and $\Sigma(t)$ the mean and
	covariance jets of $X$ at $t$, and write $\mathrm{jet}(t)$ for the pair.
	
	The blocks in \eqref{eq:Sigma} are obtained by applying the required
	covariant derivatives separately in the two arguments of $C(t,t')$ and
	then restricting to the diagonal. For example,
	\[
	G_{ij}=(\nabla_t\nabla_{t'} C)(e_i,e_j)|_{t'=t},\qquad
	\iota_{ij}=(\nabla_t^2 C)(e_i,e_j)|_{t'=t},
	\]
	and $V$ uses two covariant derivatives in each argument. One must retain
	the connection terms when differentiating in a general chart.
	Differentiating the diagonal variance function gives
	\begin{equation*}
		r_j=\tfrac12 e_jv,\qquad
		\iota_{jk}=\tfrac12(\nabla^2v)(e_j,e_k)-G_{jk}.
	\end{equation*}
	Here $\iota_{jk}=\Cov(X_{jk}(t),X(t))$ denotes the corresponding entry of
	$\iota$. The first identity follows by symmetry of $C$; the covariant
	product rule gives the second. Thus constant variance implies $r=0$,
	but need not imply $A=0$. Curvature enters only through these blocks: once they are expressed
	in an orthonormal basis of $T_t\mathcal M$, the conditioning below and the
	support classification of Section~\ref{sec:singular} involve only
	finite-dimensional matrices.
	
	\subsection{Two stages of conditioning}
	We write $\mathsf N_k(m,\Sigma)$ for the Gaussian distribution on $\R^k$
	with mean $m\in\R^k$ and covariance
	$\Sigma\in\R^{k\times k}$, allowing $\Sigma$ to be singular; and omit $k$ when $k=1$.
	First condition on $\nabla X(t)=0$ and define
	\begin{equation*}		\begin{split}
			\bar\mu&=\mu-r^\top G^{-1}\nabla\mu,\qquad
			s^2=v-r^\top G^{-1}r,
			\\
			\bar m&=m_W-AG^{-1}\nabla\mu, \qquad
			\widetilde\iota=\iota-AG^{-1}r,\qquad
			\widetilde V=V-AG^{-1}A^\top.
		\end{split}
	\end{equation*}
	Gaussian conditioning gives
	\begin{equation}\label{eq:Wgradient}
		\begin{pmatrix}X(t)\\W\end{pmatrix}
		\,\Big|\,\nabla X(t)=0
		\sim
		\mathsf N_4\!\left(
		\begin{pmatrix}\bar\mu\\\bar m\end{pmatrix},
		\begin{pmatrix}
			s^2&\widetilde\iota^\top\\
			\widetilde\iota&\widetilde V
		\end{pmatrix}
		\right).
	\end{equation}
	In particular, $W\mid\nabla X(t)=0\sim\mathsf N_3(\bar m,\widetilde V)$. To condition further on height, subtract the part of $W$ explained by
	$X(t)$. Since $s^2>0$, define the residual $W^\perp
	=W-\bar m-\widetilde\iota\bigl(X(t)-\bar\mu\bigr)/s^2$.
	Under the law conditional on $\nabla X(t)=0$, this residual is centered
	and independent of $X(t)$, with covariance $\widetilde V-\widetilde\iota\widetilde\iota^\top/s^2$. The resulting regression decomposition is $W=\bar m+\widetilde\iota\bigl(X(t)-\bar\mu\bigr)/s^2
	+W^\perp$.
	Consequently, conditioning further on $X(t)=x$ by Gaussian regression gives
	\begin{equation}\label{eq:WgivenHeightGradient}
		\begin{split}
			W\mid\{X(t)=x,\nabla X(t)=0\}
			\sim\mathsf N_3(m(x),\Omega),
		\end{split}
	\end{equation}
	where
	\[
	m(x)=\bar m+\frac{\widetilde\iota}{s^2}(x-\bar\mu),\qquad
	\Omega=\widetilde V-
	\frac{\widetilde\iota\widetilde\iota^\top}{s^2}.
	\]
	The conditional mean depends affinely on the observed height,
	whereas the conditional covariance is independent of $x$.
	These formulas remain valid when $\Omega$ is singular.
	
	Assumption~\ref{ass:regularity} ensures that $G\succ0$, $s^2>0$,
	and $\widetilde V\succ0$, but allows $\Omega\succeq0$ to be singular.
	The regression formulas require only the inverses of $G$ and $s^2$,
	so they apply in all these cases. Theorem~\ref{thm:singularrank}
	classifies the possible singularities.
	
	Let $\phi$ and $\Phi$ denote the standard normal density and distribution
	function. Writing $\varphi(x;a,b)$ for the
	$\mathsf N(a,b)$ density, the remaining factors
	in \eqref{eq:rho} are
	\begin{equation*}
		p_{X(t)\mid \nabla X(t)=0}(x)=\phin{x}{\bar\mu}{s^2},\qquad
		p_{\nabla X(t)}(0)=\frac{\exp\{-(\nabla\mu)^\top G^{-1}\nabla\mu/2\}}{2\pi\sqrt{\det G}}.
	\end{equation*}

	Table~\ref{tab:notation} collects the principal notation.
	
	\begin{table}[htbp]
		\centering
		\small
		\setlength{\tabcolsep}{4pt}
		\renewcommand{\arraystretch}{1.15}
		\caption{Principal notation.}
		\label{tab:notation}
		\begin{tabular}{@{}p{0.30\textwidth}p{0.52\textwidth}l@{}}
			\toprule
			Symbol & Meaning & Location \\
			\midrule
			$D$, $N(D,B)$, $\rho$, $\ell$, $h_t$ & domain, peak counts, intensities, local height density & \S\ref{sec:framework} \\
			$J(t,x)$, $J(t)$ & conditional peak weight and its normalizer & \S\ref{sec:Kac-Rice} \\
			$v,r,G,\iota,A,V$ & covariance blocks of the jet & \eqref{eq:Sigma} \\
			$\bar\mu,s^2,\bar m,\widetilde\iota,\widetilde V$ & regression parameters given $\nabla X(t)=0$ & \S\ref{sec:regression} \\
			$m(x),\Omega$ & Hessian law given $X(t)=x$ and $\nabla X(t)=0$ & \eqref{eq:WgivenHeightGradient} \\
			$\tau,\xi$, $L$ & trace and traceless coordinates, $(\tau,\xi^\top)^\top=LW$ & \eqref{eq:L} \\
			$m_\tau,m_\xi,\sigma_\tau^2,\sigma_{\xi\tau},S_\xi$ & mean and covariance of the trace coordinates & \eqref{eq:m-tau} \\
			$\eta(z),\Lambda$ & law of $\xi$ given $\tau=z$ & \eqref{eq:traceconditional} \\
			$\Acal,\Acal_{\mathrm{sdl}}$; $\mathcal H$ & norm kernels for peaks and saddles; scalar kernel & \eqref{eq:kernelA}, \eqref{eq:kernelH} \\
			$\mathfrak D$, $\Kcal$, $\mathfrak D_{\mathrm{sdl}}$, $\Kcal_{\mathrm{sdl}}$ & peak weight, expected peak weight; saddle analogues & \S\ref{sec:trace} \\
			$\Upsilon,\mathsf R,\vartheta,\lambda_1,\lambda_2$ & determinant quadratic form and its spectrum & \eqref{eq:determinantspectrum} \\
			$\tau_*,\xi_*,\eta_*,\Lambda_*,\breve m,\varsigma_0$ & spectral trace, its norm coordinates & Thm.~\ref{thm:spectral} \\
			$b_*,\varepsilon_*$ & regression of $X(t)$ on $W$, residual variance & Thm.~\ref{thm:singularrank} \\
			$Q,x_*$ & constraint operator and endpoint & \S\ref{subsec:singularsupport} \\
			$\mathcal Q_t$ & second-order operator of a singular jet & \eqref{eq:pointwisepde} \\
			$P,M,\varpi_P,\delta,\varsigma$ & canonical coordinates & \eqref{eq:sigma} \\
			$\theta,\Theta$ & regression parameter vectors and their set  & \S\ref{subsec:singularconvergence} \\
			$\nu_\theta,W^\circ,\mathcal L_\theta$ & tilted law, tilted variable, height map & Prop. \ref{prop:tiltedrep} \\
			\bottomrule
		\end{tabular}
	\end{table}
	
	\section{Trace reductions and peak intensity formulas}
	\label{sec:trace}
	We now evaluate the expected peak weight for an arbitrary Gaussian Hessian
	law. Appendix~\ref{app:normmoments} collects the CDF and moment
	calculations.
	
	\subsection{Trace coordinates and the peak kernel}
	Decompose the Hessian into its trace and traceless parts by setting
	\begin{equation*}		\tau=X_{11}(t)+X_{22}(t),\qquad
		U=\frac{X_{11}(t)-X_{22}(t)}2,\qquad
		X_{12}=X_{12}(t),\qquad
		\xi=(U,X_{12})^\top.
	\end{equation*}
	The resulting decomposition is
	\begin{equation}\label{eq:Hessian-decomp}
		\nabla^2X(t)
		=
		\underbrace{\frac{\tau}{2}I_2}_{\text{trace part}}
		+
		\underbrace{
			\begin{pmatrix}
				U&X_{12}\\
				X_{12}&-U
		\end{pmatrix}}_{\text{traceless part}}.
	\end{equation}
	Here $\tau=\Delta_{\mathcal M}X(t)$. Thus $\tau$ specifies the trace part, and $\xi$ contains the two
	coordinates of the traceless part. In the Hessian entry order
	$W=(X_{11}(t),X_{22}(t),X_{12}(t))^\top$, the transformation is
	\begin{equation}\label{eq:L}
		\begin{pmatrix}\tau\\ \xi\end{pmatrix}
		=LW,\qquad
		L=\begin{pmatrix}
			1&1&0\\
			1/2&-1/2&0\\
			0&0&1
		\end{pmatrix}.
	\end{equation}
	The matrix $L$ is invertible, so the transformation preserves the full
	Hessian information.
	
	\begin{lemma}[Negative definite Hessians]
		\label{lem:cones}
		The Hessian eigenvalues are
		$\tau/2\pm\norm\xi$, so
		\begin{equation}
			\det\nabla^2X(t)=\frac{\tau^2}{4}-\norm\xi^2,\qquad
			\nabla^2X(t)\prec0
			\quad\Longleftrightarrow\quad\tau<0,\ \norm\xi^2<\tau^2/4.
			\label{eq:dettrace}
		\end{equation}
	\end{lemma}
	\begin{proof}
		By \eqref{eq:Hessian-decomp} and $\norm\xi^2=U^2+X_{12}^2$, the characteristic polynomial of the Hessian is
		\begin{equation*}
			\det(y I_2-\nabla^2X(t))
			= (y-\tau/2-U)(y-\tau/2+U)-X_{12}^2
			= (y-\tau/2)^2-\norm\xi^2.
		\end{equation*}
		Thus the two eigenvalues are
		$y^{\pm}=\tau/2\pm \norm\xi$.
		Their product gives the determinant in \eqref{eq:dettrace}.
		Since the Hessian is symmetric, it is negative definite if and only if
		$y^+<0$, that is, $\norm\xi<-\tau/2$, which holds if
		and only if $\tau<0$ and $\norm\xi^2<\tau^2/4$.
	\end{proof}
	
	The peak constraint and determinant depend on the three
	Hessian entries only through two scalar quantities: the trace
	$\tau$ and the squared norm $\norm\xi^2$ of the traceless
	coordinates. For each fixed $\tau<0$, the admissible values of
	$\xi$ form a disk of radius $-\tau/2$, within which the determinant
	weight $\tau^2/4-\norm\xi^2$ decreases from $\tau^2/4$ at the
	center to zero at the boundary. Conditioning on
	$\tau$ therefore isolates the remaining calculation in the distribution
	of the squared Gaussian norm, including its nonzero mean and general
	covariance.
	
	For $\xi\sim\mathsf N_2(\eta,\Lambda)$ with $\Lambda\succeq0$, define
	the Gaussian norm kernel
	\begin{equation}\label{eq:kernelA}
		\Acal(a;\eta,\Lambda)=\E[(a-\norm\xi^2)_+],\qquad a\ge0.
	\end{equation}
	Here $(y)_+=\max\{y,0\}$. The kernel is a lower partial moment of
	$\norm\xi^2$, whose distribution is a possibly degenerate generalized
	chi square law. If $F_{\eta,\Lambda}$ denotes the CDF of $\norm\xi^2$, then
	\begin{equation}
		\Acal(a;\eta,\Lambda)=\int_0^aF_{\eta,\Lambda}(y)\dd y.
		\label{eq:Alayer}
	\end{equation}
	Indeed, with $Y=\norm\xi^2$, the identity
	$(a-Y)_+=\int_0^a\ind_{\{Y\le y\}}\dd y$ and Fubini's theorem
	give \eqref{eq:Alayer}.
	
	Let $O(2)$ denote the group of real $2\times2$ orthogonal matrices
	satisfying $O^\top O=I_2$. For a real matrix $M=(M_{ij})\in\mathbb R^{m\times n}$, its
	Frobenius norm is $\|M\|_{\mathrm F}
	=\big(\sum_{i=1}^m\sum_{j=1}^n M_{ij}^2\big)^{1/2}
	=\{\operatorname{tr}(M^\top M)\}^{1/2}$.
	
	\begin{remark}[Frame invariance of the norm kernel]
		\label{rem:frameinvariance}
		Under an orthonormal frame change, the Hessian matrix
		$H=\nabla^2X(t)$ becomes $H'=O^\top HO$ with $O\in O(2)$; the trace is
		unchanged and $\xi'=R_O\xi$, where $R_O$ is a
		rotation through $-2\psi$ if $O$ is a rotation through $\psi$, and
		$R_O=\diag(1,-1)$ if $O=\diag(1,-1)$. Indeed, $\tr(O^\top HO)=\tr H$, and identifying the traceless part of $H$
		with $U+\mathrm iX_{12}$, conjugation by a rotation through $\psi$
		multiplies it by $e^{-2\mathrm i\psi}$, while conjugation by
		$\diag(1,-1)$ replaces it by its complex conjugate. Consequently, if $\xi\sim\mathsf N_2(\eta,\Lambda)$ in one frame, then
		$\xi'=R_O\xi\sim\mathsf N_2(R_O\eta,R_O\Lambda R_O^\top)$ in the other, and
		$\norm{\xi'}=\norm\xi$, so that
		$\Acal(a;R_O\eta,R_O\Lambda R_O^\top)=\Acal(a;\eta,\Lambda)$.
		Because the trace is the same in both frames, the same identity holds for
		the conditional law of $\xi$ given $\tau=z$; hence each integrand in
		\eqref{eq:substitution} is frame independent. The frame independence of
		$J(t,x)$ itself does not rely on this: $J(t,x)$ is the conditional
		expectation of $|\det H|\ind_{\{H\prec0\}}$, which is unchanged under
		$H\mapsto O^\top HO$, given the event $\{X(t)=x,\nabla X(t)=0\}$ which does
		not depend on the frame.
	\end{remark}
	
	\subsection{Reduction theorem and peak formulas}
	\label{subsec:reduction}
	
	For $w=(w_1,w_2,w_3)^\top\in\R^3$, let $H(w)$ be the symmetric
	$2\times2$ matrix with $\vech H(w)=w$; for $m=(m_1,m_2,m_3)^\top$ and
	$S\succeq0$, define the peak weight and the expected peak weight by
	\[
	\mathfrak D(w)=|\det H(w)|\ind_{\{H(w)\prec0\}},\qquad
	\Kcal(m,S)=\E[\mathfrak D(W)],\qquad W\sim\mathsf N_3(m,S).
	\]
	Thus $\Kcal(m,S)$ is the expected peak weight of the Hessian $H(W)$
	when its entries have the law $\mathsf N_3(m,S)$. We give the trace
	parameters for a general pair $(m,S)$, using $(m(x),\Omega)$ for the
	conditional peak weight $J(t,x)$ and $(\bar m,\widetilde V)$ for its normalizer $J(t)$.
	
	The trace coordinates in \eqref{eq:L} have mean and covariance
	\begin{equation}\label{eq:traceblocks}
		Lm=\binom{m_\tau}{m_\xi},\qquad
		LSL^\top=
		\begin{pmatrix}\sigma_\tau^2&\sigma_{\xi\tau}^\top\\\sigma_{\xi\tau}&S_\xi\end{pmatrix},
	\end{equation}
	where
	\begin{equation}\label{eq:m-tau}
		\begin{split}
			m_\tau&=m_1+m_2,\quad
			m_\xi=((m_1-m_2)/2,m_3)^\top,\quad
			\sigma_\tau^2=S_{11}+S_{22}+2S_{12},\\
			\sigma_{\xi\tau}&=((S_{11}-S_{22})/2,S_{13}+S_{23})^\top,\quad
			S_\xi=\begin{pmatrix}
				(S_{11}+S_{22}-2S_{12})/4&(S_{13}-S_{23})/2\\
				(S_{13}-S_{23})/2&S_{33}
			\end{pmatrix}.
		\end{split}
	\end{equation}
	When $\sigma_\tau^2>0$, Gaussian conditioning gives
	\begin{equation}\label{eq:traceconditional}
		\xi\mid\tau=z\sim\mathsf N_2(\eta(z),\Lambda),\qquad
		\eta(z)=m_\xi+\frac{\sigma_{\xi\tau}}{\sigma_\tau^2}(z-m_\tau),\qquad
		\Lambda=S_\xi-\frac{\sigma_{\xi\tau}\sigma_{\xi\tau}^\top}{\sigma_\tau^2},
	\end{equation}
	where the conditional covariance matrix $\Lambda\succeq0$ may be singular (Proposition~\ref{prop:singularcdfnorm}).
	
	\begin{theorem}[Trace reduction of the conditional peak weight]
		\label{thm:main}
		Fix $t\in\mathcal M$ and choose any orthonormal tangent frame.
		Under Assumption~\ref{ass:regularity}, evaluate
		\eqref{eq:traceblocks}--\eqref{eq:traceconditional} at
		$(m,S)=(m(x),\Omega)$. With the kernel $\Acal$ defined in
		\eqref{eq:kernelA}, the conditional peak weight is
		\begin{equation}\label{eq:substitution}
			J(t,x)=\Kcal(m(x),\Omega)=
			\begin{cases}
				\displaystyle\int_{-\infty}^0
				\phin{z}{m_\tau}{\sigma_\tau^2}
				\Acal(z^2/4;\eta(z),\Lambda)\dd z,
				&\sigma_\tau^2>0,\\[6pt]
				\displaystyle\ind_{\{m_\tau<0\}}
				\Acal(m_\tau^2/4;m_\xi,S_\xi),
				&\sigma_\tau^2=0.
			\end{cases}
		\end{equation}
		The normalizer follows either from the same reduction
		with $(m,S)=(\bar m,\widetilde V)$ or by averaging $J(t,x)$ over height:
		\begin{equation}\label{eq:normalizerweight}
			J(t)=\Kcal(\bar m,\widetilde V) \qquad {\rm or } \qquad	J(t)=\int_\R\phin{x}{\bar\mu}{s^2}J(t,x)\dd x.
		\end{equation}
		Both $J(t,x)$ and $J(t)$ are finite, and $J(t)>0$. Combined with
		\eqref{eq:KR}--\eqref{eq:peak-height}, they give the peak intensity,
		the expected number of peaks above any level, and the local peak
		height density of a Gaussian field with arbitrary mean and
		covariance.
	\end{theorem}
	
	\begin{proof}
		The conditioning laws \eqref{eq:Wgradient} and
		\eqref{eq:WgivenHeightGradient} identify
		$J(t,x)=\Kcal(m(x),\Omega)$ and $J(t)=\Kcal(\bar m,\widetilde V)$.
		Let $W\sim\mathsf N_3(m(x),\Omega)$ and $(\tau,\xi^\top)^\top=LW$.
		Lemma~\ref{lem:cones} gives the pointwise identity
		\[
		\mathfrak D(W)
		=\ind_{\{\tau<0\}}(\tau^2/4-\norm\xi^2)_+.
		\]
		If $\sigma_\tau^2>0$, conditioning on the trace $\tau=z$ and using
		\eqref{eq:traceconditional} and \eqref{eq:kernelA} gives
		\begin{align*}
			\Kcal(m(x),\Omega)
			&=\int_{-\infty}^0\phin{z}{m_\tau}{\sigma_\tau^2}
			\E\big[(z^2/4-\norm\xi^2)_+\mid\tau=z\big]\dd z\\
			&=\int_{-\infty}^0\phin{z}{m_\tau}{\sigma_\tau^2}
			\Acal(z^2/4;\eta(z),\Lambda)\dd z,
		\end{align*}
		proving
		the first branch of \eqref{eq:substitution}. If
		$\sigma_\tau^2=0$, positive semidefiniteness forces $\sigma_{\xi\tau}=0$.
		Then $\tau=m_\tau$ almost surely and
		$\xi\sim\mathsf N_2(m_\xi,S_\xi)$, giving the second branch.
		
		The bound $0\le J(t,x)\le (m_\tau^2+\sigma_\tau^2)/4$, in which
		$m_\tau=m_\tau(x)$ is affine in $x$, is
		integrable against $\varphi(x;\bar\mu,s^2)$. Conditional total expectation
		therefore gives \eqref{eq:normalizerweight} and finiteness
		of $J(t)$. The Hessian law given zero gradient has covariance
		$\widetilde V\succ0$, so its density is positive on the open
		negative definite cone and $J(t)>0$.
	\end{proof}
	The reduction holds for every $S\succeq0$, including the singular
	height-conditioned laws of Section~\ref{sec:singular}; for $S\succ0$,
	Theorem~\ref{thm:spectral} below gives a form in which the spectral
	trace is independent of the norm.
	
	\subsection{Norm kernels and their evaluation}
	\label{subsec:finiteCDF}
	\label{subsec:quadrature}
	The reduction uses the two-dimensional kernel $\Acal$ and its scalar counterpart
	\begin{equation}\label{eq:kernelH}
		\mathcal H(b;\mu_0,\lambda)
		=\E\big[(b-(\mu_0+\sqrt\lambda\,\mathcal N)^2)_+\big],
		\qquad \mathcal N\sim\mathsf N(0,1),
	\end{equation}
	where $b,\mu_0\in\R$ and $\lambda\ge0$ is a variance.
	The scalar kernel is even in $\mu_0$, vanishes for $b\le0$, and
	satisfies $\mathcal H(b;\mu_0,0)=(b-\mu_0^2)_+$.
	For a centered circular Gaussian norm, the exponential law gives $F_{0,\lambda I_2}(y)=1-e^{-y/(2\lambda)}$, $y\ge 0$, and hence
	\begin{equation}\label{eq:centeredcircularkernel}
		\begin{split}
			\Acal(a;0,\lambda I_2)=a-2\lambda+2\lambda e^{-a/(2\lambda)},
			\qquad a\ge0,\quad\lambda>0.
		\end{split}
	\end{equation}
	
	For $\Lambda\succ0$, choose an orthogonal matrix $O$ such that
	$\Lambda=O\diag(\lambda_1,\lambda_2)O^\top$ and write
	$\widehat\eta=O^\top\eta$. The coordinates of $O^\top\xi$ are independent,
	so conditioning on the first gives
	\begin{equation}\label{eq:normconvolution}
		\Acal(a;\eta,\Lambda)
		=\int_{-\sqrt a}^{\sqrt a}
		\phin{w}{\widehat\eta_1}{\lambda_1}
		\mathcal H(a-w^2;\widehat\eta_2,\lambda_2)\dd w,
		\qquad a\ge0,
	\end{equation}
	and the roles of the two coordinates may be interchanged.
	
	Appendix~\ref{app:normmoments} gives closed forms for $\Acal$ in three
	cases. For equal norm eigenvalues (Proposition~\ref{prop:ordinaryequal})
	or a centered norm with unequal eigenvalues
	(Proposition~\ref{prop:ordinarycentered}), $\Acal$ is a finite combination
	of ordinary noncentral chi square CDFs, and a rank-one norm covariance
	reduces $\Acal$ to $\mathcal H$ (Proposition~\ref{prop:singularcdfnorm}).
	For a noncentered norm with unequal eigenvalues,
	\eqref{eq:normconvolution} leaves one scalar integral of ordinary
	noncentral chi square CDFs; alternatively, $\Acal$ is a finite combination
	of generalized noncentral chi square CDFs by \eqref{eq:normgeneralizedcdf}, or a convergent series of central
	chi square CDFs by Ruben's expansion \citep{Ruben1962}.
	
	Consequently, $J(t,x)$ in \eqref{eq:substitution} requires at most two
	scalar integrations, over the trace and over one norm coordinate; the
	finite CDF formulas remove the second, and a deterministic trace removes
	the first.
	
	\subsection{Spectral trace reduction and centered normalizers}
	\label{subsec:centeredcounts}
	Theorem~\ref{thm:main} conditions on the trace, so the law of the
	traceless coordinates moves with the trace value. For a nondegenerate
	Hessian law, diagonalizing the determinant instead produces a spectral
	trace $\tau_*$ that is independent of a norm with a fixed law. The
	spectral trace is a coordinate adapted to the determinant quadratic form;
	in general it differs from the geometric trace $\tr H$. Write
	$\det H(w)=w^\top\Upsilon w$, where
	\begin{equation}\label{eq:determinantspectrum}
		\Upsilon=\begin{pmatrix}
			0&1/2&0\\1/2&0&0\\0&0&-1
		\end{pmatrix}.
	\end{equation}
	
	\begin{theorem}[Spectral trace reduction]
		\label{thm:spectral}
		Let $m\in\R^3$ and $S\succ0$. The matrix
		$\mathsf R=S^{1/2}\Upsilon S^{1/2}$ has eigenvalues $\vartheta>0$,
		$-\lambda_1$ and $-\lambda_2$ with $\lambda_1\ge\lambda_2>0$. Let
		$o_0,o_1,o_2$ be corresponding orthonormal eigenvectors, set
		$\breve m_j=o_j^\top S^{-1/2}m$, and let $\varsigma_0\in\{-1,1\}$ be the
		sign of $\tr H(S^{1/2}o_0)$. Let
		$\tau_*\sim\mathsf N(2\varsigma_0\sqrt\vartheta\,\breve m_0,4\vartheta)$ and
		$\xi_*\sim\mathsf N_2(\eta_*,\Lambda_*)$ be independent, where
		$\eta_*=(\sqrt{\lambda_1}\,\breve m_1,\sqrt{\lambda_2}\,\breve m_2)^\top$ and
		$\Lambda_*=\diag(\lambda_1,\lambda_2)$. Then
		\begin{equation}\label{eq:spectral}
			\begin{aligned}
				\Kcal(m,S)&=\E\bigl[(\tau_*^2/4-\norm{\xi_*}^2)_+\ind_{\{\tau_*<0\}}\bigr]
				=\int_{-\infty}^0\phin{z}{2\varsigma_0\sqrt\vartheta\,\breve m_0}{4\vartheta}\,
				\Acal(z^2/4;\eta_*,\Lambda_*)\dd z\\
				&=\int_0^\infty F_{\eta_*,\Lambda_*}(y)\,
				\Phi\Bigl(-\frac{\sqrt y+\varsigma_0\sqrt\vartheta\,\breve m_0}{\sqrt\vartheta}\Bigr)\dd y .
			\end{aligned}
		\end{equation}
	\end{theorem}
	\begin{proof}
		The eigenvalue signs follow from Sylvester's law of inertia, since
		$\Upsilon $ has one positive and two negative eigenvalues. Write
		$W=m+S^{1/2}\mathcal V$ with $\mathcal V\sim\mathsf N_3(0,I_3)$, and let
		$\mathcal N_j=o_j^\top\mathcal V$, $j=0,1,2$, which are independent
		$\mathsf N(0,1)$. Since $S^{-1/2}m=\sum_j\breve m_jo_j$, we have
		$S^{-1/2}W=\sum_j(\breve m_j+\mathcal N_j)o_j$, and
		$\det H(W)=W^\top\Upsilon  W=(S^{-1/2}W)^\top\mathsf R(S^{-1/2}W)$ gives
		\[
		\det H(W)=\vartheta(\breve m_0+\mathcal N_0)^2-\lambda_1(\breve m_1+\mathcal N_1)^2
		-\lambda_2(\breve m_2+\mathcal N_2)^2=\tau_*^2/4-\norm{\xi_*}^2,
		\]
		where $\tau_*=2\varsigma_0\sqrt\vartheta(\breve m_0+\mathcal N_0)$ and
		$\xi_*=(\sqrt{\lambda_1}(\breve m_1+\mathcal N_1),\sqrt{\lambda_2}(\breve m_2+\mathcal N_2))^\top$.
		These are independent and have the laws stated in the theorem. It
		remains to express the event $\{H(W)\prec0\}$ through $\tau_*$.
		
		The sign $\varsigma_0$ is well defined, because
		$\det H(S^{1/2}o_0)=o_0^\top\mathsf Ro_0=\vartheta>0$ makes
		$H(S^{1/2}o_0)$ definite, with nonzero trace. On
		$\{\det H(W)>0\}$ the matrix $H(W)$ is definite, so $H(W)\prec0$ exactly
		when $\tr H(W)<0$. With
		$\varrho=\{(\lambda_1(\breve m_1+\mathcal N_1)^2+\lambda_2(\breve m_2+\mathcal N_2)^2)/\vartheta\}^{1/2}$,
		the event $\{\det H(W)>0\}$ is the union of $\{\breve m_0+\mathcal N_0>\varrho\}$ and
		$\{\breve m_0+\mathcal N_0<-\varrho\}$. Each part is connected, being the
		region above or below the graph of a continuous function of
		$(\mathcal N_1,\mathcal N_2)$, and $\tr H(W)$ does not vanish on it, so the
		trace has constant sign on each part. By linearity,
		$\tr H(W)=\sum_j(\breve m_j+\mathcal N_j)\tr H(S^{1/2}o_j)$, and letting
		$\mathcal N_0\to\pm\infty$ with $\mathcal N_1,\mathcal N_2$ fixed shows that
		this sign is $\varsigma_0$ on the first part and $-\varsigma_0$ on the
		second. Hence $H(W)\prec0$ exactly when $\det H(W)>0$ and $\tau_*<0$, and
		then $|\det H(W)|=\tau_*^2/4-\norm{\xi_*}^2$, which gives the first
		expression in \eqref{eq:spectral}. Conditioning on $\tau_*$ and using
		\eqref{eq:kernelA} gives the second. For the third, let
		$Y=\norm{\xi_*}^2$; since
		$(\tau_*^2/4-Y)_+\ind_{\{\tau_*<0\}}
		=\int_0^\infty\ind_{\{Y\le y\}}\ind_{\{\tau_*<-2\sqrt y\}}\dd y$, Fubini's
		theorem and the independence of $\tau_*$ and $Y$ give
		$\int_0^\infty F_{\eta_*,\Lambda_*}(y)\,\PP(\tau_*<-2\sqrt y)\dd y$, with
		$\PP(\tau_*<-2\sqrt y)=\Phi(-(\sqrt y+\varsigma_0\sqrt\vartheta\,\breve m_0)/\sqrt\vartheta)$.
	\end{proof}
	
	The expressions do not depend on the signs of the eigenvectors: replacing
	$o_0$ by $-o_0$ changes the signs of both $\breve m_0$ and $\varsigma_0$,
	and replacing $o_1$ or $o_2$ by its negative changes only the sign of a
	coordinate of $\eta_*$, which leaves the law of $\norm{\xi_*}$ unchanged,
	as does any choice of eigenvectors when $\lambda_1=\lambda_2$. Unlike in
	Theorem~\ref{thm:main}, the norm law in \eqref{eq:spectral} does not
	depend on the trace value and has diagonal covariance, so the last
	expression is a single integral against one fixed distribution function.
	That function is an ordinary noncentral chi square CDF when
	$\lambda_1=\lambda_2$ (Proposition~\ref{prop:ordinaryequal}), a finite
	combination of such CDFs when $\eta_*=0$
	(Proposition~\ref{prop:ordinarycentered}), and in general the generalized
	chi square CDF $F^{\mathrm{gen}}_{1,1}$ of
	Remark~\ref{remark:ordinarynoncentered}, which conditioning on one
	coordinate writes as one integral of ordinary noncentral chi square CDFs.
	Theorem~\ref{thm:spectral} applies to the normalizer
	$J(t)=\Kcal(\bar m,\widetilde V)$, since $\widetilde V\succ0$, and to
	$J(t,x)=\Kcal(m(x),\Omega)$ when the jet law is nonsingular; then
	$\mathsf R$ and its eigenvectors do not depend on $x$, and
	$\breve m=(\breve m_0,\breve m_1,\breve m_2)^\top$ is affine in $x$. As
	$\Omega$ approaches a singular matrix, an eigenvalue of $\mathsf R$ tends
	to zero and the corresponding coordinate of $\breve m$ may grow, so the
	representation becomes ill-conditioned; there Theorems~\ref{thm:main}
	and~\ref{thm:canonicalweights} are the appropriate tools. For centered
	laws, $\breve m=0$, and the theorem reduces to the following normalizer.
	By Proposition~\ref{prop:ordinarycentered}, it gives the peak intensity
	$\ell(t)=p_{\nabla X(t)}(0)J(t)$ of a centered field as a single scalar
	integral of ordinary noncentral chi square CDFs, or in closed form when the
	two eigenvalues coincide, and \eqref{eq:KR} turns it into the expected
	total number of peaks.
	\begin{corollary}[Normalizer for a centered Hessian law]
		\label{cor:centerednormalizer}
		Under Assumption~\ref{ass:regularity}, suppose $\bar m(t)=0$; this holds, for instance, if $X$ is centered
		or if $X$ is a stationary planar field plus an affine trend.
		Let $\vartheta(t)>0$, $-\lambda_1(t)$ and $-\lambda_2(t)$, with
		$\lambda_1(t)\ge\lambda_2(t)>0$, be the eigenvalues of
		$\mathsf R(t)=\widetilde V(t)^{1/2}\Upsilon \widetilde V(t)^{1/2}$.
		Suppressing $t$ and setting
		$\Lambda_*=\diag(\lambda_1,\lambda_2)$, the normalizer is
		\begin{equation}\label{eq:centerednormalizer}
			J(t)=\int_{-\infty}^0\phin{z}{0}{4\vartheta}
			\Acal(z^2/4;0,\Lambda_*)\dd z.
		\end{equation}
		If $\lambda_1>\lambda_2$, the kernel is the finite ordinary CDF
		expression \eqref{eq:centeredordinarykernel}.
		If $\lambda_1=\lambda_2=\lambda$, then
		\begin{equation}\label{eq:centerednormalizercircular}
			J(t)=\frac{\vartheta}{2}-\lambda
			+\lambda\sqrt{\frac{\lambda}{\lambda+\vartheta}}.
		\end{equation}
	\end{corollary}
	\begin{proof}
		Theorem~\ref{thm:spectral} with $m=\bar m=0$ and $S=\widetilde V$
		gives $\breve m=0$, so $\eta_*=0$ and $\tau_*\sim\mathsf N(0,4\vartheta)$, and
		its second expression is \eqref{eq:centerednormalizer}.
		When $\lambda_1=\lambda_2=\lambda$, plugging
		\eqref{eq:centeredcircularkernel} gives
		\eqref{eq:centerednormalizercircular}.
	\end{proof}
	
	\subsection{Local minima and saddle points}
	\label{subsec:saddles}
	A local minimum of $X$ is a local maximum of $-X$, which has mean
	$-\mu$ and the same covariance. Hence $J^X_{\mathrm{min}}(t,x)=J^{-X}(t,-x)$,
	$h^X_{\mathrm{min}}(t,x)=h^{-X}(t,-x)$ and
	$N^X_{\mathrm{min}}(D,B)=N^{-X}(D,-B)$; reversing the height alone does not
	give the minimum law of a noncentered field.
	
	For $w\in\R^3$, with $H(w)$ as in Section~\ref{subsec:reduction}, and
	for $m\in\R^3$ and $S\succeq0$, define the saddle weight and the expected saddle weight by
	\[
	\mathfrak D_{\mathrm{sdl}}(w)=|\det H(w)|\ind_{\{\det H(w)<0\}},\qquad
	\Kcal_{\mathrm{sdl}}(m,S)=\E[\mathfrak D_{\mathrm{sdl}}(W)],\qquad W\sim\mathsf N_3(m,S).
	\]
	A saddle on a surface has one positive and one negative Hessian eigenvalue,
	equivalently a negative determinant. Thus \eqref{eq:dettrace} gives
	the pointwise identity
	\begin{equation*}
		\mathfrak D_{\mathrm{sdl}}(W)=\big(\norm\xi^2-\tau^2/4\big)_+,\qquad(\tau,\xi^\top)^\top=LW.
	\end{equation*}
	For $\xi\sim\mathsf N_2(\eta,\Lambda)$ with $\Lambda\succeq0$ and
	$Y=\norm\xi^2$, the identity $(Y-a)_+=Y-a+(a-Y)_+$ gives the
	complementary kernel
	$\Acal_{\mathrm{sdl}}(a;\eta,\Lambda)=\E[(Y-a)_+]
	=\tr(\Lambda)+\norm\eta^2-a+\Acal(a;\eta,\Lambda)$, $a\ge0$.
	Conditioning on $\tau=z$ gives, for $\sigma_\tau^2>0$,
	\begin{equation}
		\Kcal_{\mathrm{sdl}}(m,S)=\int_\R\phin{z}{m_\tau}{\sigma_\tau^2}
		\Acal_{\mathrm{sdl}}(z^2/4;\eta(z),\Lambda)\dd z,
		\label{eq:Ksaddle}
	\end{equation}
	with integration over all real traces, and
	$\Kcal_{\mathrm{sdl}}(m,S)=\Acal_{\mathrm{sdl}}(m_\tau^2/4;m_\xi,S_\xi)$
	when $\sigma_\tau^2=0$. The conditional saddle weight and its normalizer are
	$J_{\mathrm{sdl}}(t,x)=\Kcal_{\mathrm{sdl}}(m(x),\Omega)$ and
	$J_{\mathrm{sdl}}(t)=\Kcal_{\mathrm{sdl}}(\bar m,\widetilde V)>0$, and
	\eqref{eq:rho}--\eqref{eq:peak-height} give the saddle intensities and
	height density. For $S\succ0$, the proof of Theorem~\ref{thm:spectral}
	gives likewise $\Kcal_{\mathrm{sdl}}(m,S)=\E[(\norm{\xi_*}^2-\tau_*^2/4)_+]$,
	with no sign condition on $\tau_*$.

	\section{Singular jet laws and support geometry}
	\label{sec:singular}
	At a fixed location, a singular jet law confines the
	height-conditioned Hessian law to an affine plane, and the intersection of this plane with the negative
	definite cone determines the possible peak heights. Throughout this
	section $s^2>0$ and $\widetilde V\succ0$ (Assumption~\ref{ass:regularity}).
	Consequently a singular conditional Hessian covariance has rank two, and
	singularity of the law is compatible with the Morse property: the
	Hessian at a critical point is still nonsingular almost surely. Under
	$H\mapsto O^\top HO$ the constraint operator $Q$ below transforms as
	$Q\mapsto O^\top QO$, so its signature and the endpoint $x_*$ do not
	depend on the frame. Statements are local at $t$ unless a field is
	specified.
	
	\subsection{Admissible singular laws and affine support}
	\label{subsec:singular}
	Recall from \eqref{eq:WgivenHeightGradient} that
	$W\mid\{X(t)=x,\nabla X(t)=0\}\sim\mathsf N_3(m(x),\Omega)$, with
	$s^2>0$ and $\widetilde V\succ0$ under Assumption~\ref{ass:regularity}. The next theorem
	shows that the rank of $\Omega$ can only be two or three.
	
	\begin{theorem}[Classification of singular jet laws]
		\label{thm:singularrank}
		Under Assumption~\ref{ass:regularity}, let
		\begin{equation*}			b_* = \widetilde V^{-1}\widetilde\iota,\qquad
			\varepsilon_* = \Var(X(t)\mid \nabla X(t), W) = s^2-
			\widetilde\iota^\top\widetilde V^{-1}\widetilde\iota.
		\end{equation*}
		Then $\varepsilon_*\ge0$. If $\varepsilon_*>0$, then
		$\Omega\succ0$. If $\varepsilon_*=0$, then
		$\operatorname{rank}(\Omega)=2$, $b_*\ne0$, and
		$\ker(\Omega)=\operatorname{span}\{b_*\}$. In the singular case $\varepsilon_*=0$, under the law conditional on $\nabla X(t)=0$,
		\begin{equation}
			X(t)-\bar\mu=b_*^\top(W-\bar m)\qquad\text{almost surely}.
			\label{eq:singularheightrelation}
		\end{equation}
		Consequently, for every $x\in\R$, the conditional law in
		\eqref{eq:WgivenHeightGradient} has support equal to the affine plane
		\begin{equation}\label{eq:singularheightplane}
			\left\{w\in\R^3:
			b_*^\top(w-\bar m)=x-\bar\mu\right\}.
		\end{equation}
		Moreover, $\varepsilon_*=\det\Sigma(t)/\det\Cov\bigl((\nabla X(t)^\top,W^\top)^\top\bigr)$,
		so the jet law is singular if and only if $\varepsilon_*=0$.
	\end{theorem}
	\begin{proof}
		Conditional on $\nabla X(t)=0$, the covariance matrix of
		$(X(t),W^\top)^\top$ is
		\[
		\begin{pmatrix}
			s^2&\widetilde\iota^\top\\
			\widetilde\iota&\widetilde V
		\end{pmatrix}\succeq0.
		\]
		Since $\widetilde V\succ0$, its Schur complement
		$\varepsilon_*$ is nonnegative. To classify $\Omega$, set
		$u_*=\widetilde V^{-1/2}\widetilde\iota/s$. By
		\eqref{eq:WgivenHeightGradient},
		\[
		\Omega=\widetilde V^{1/2}
		(I_3-u_*u_*^\top)\widetilde V^{1/2},
		\qquad
		1-\norm{u_*}^2=\frac{\varepsilon_*}{s^2}.
		\]
		The eigenvalues of $I_3-u_*u_*^\top$ are
		$1,1,\varepsilon_*/s^2$. Since congruence by an invertible
		matrix preserves rank and positive definiteness, $\Omega\succ0$
		when $\varepsilon_*>0$, whereas $\operatorname{rank}(\Omega)=2$
		when $\varepsilon_*=0$.
		
		In the latter case,
		$\widetilde\iota^\top\widetilde V^{-1}\widetilde\iota=s^2>0$,
		so $\widetilde\iota\ne0$ and hence $b_*\ne0$. Moreover, $\Omega b_*
		=\widetilde\iota\left(
		1-\widetilde\iota^\top\widetilde V^{-1}\widetilde\iota/s^2
		\right)=0$.
		Thus $b_*$ spans the one-dimensional nullspace of $\Omega$. Gaussian regression gives
		\begin{equation*}
			\E[X(t)\mid W,\nabla X(t)=0]
			=\bar\mu+b_*^\top(W-\bar m),\qquad
			\Var(X(t)\mid W,\nabla X(t)=0)=\varepsilon_*.
		\end{equation*}
		Thus $\varepsilon_*=0$ implies
		\eqref{eq:singularheightrelation}. Fixing the height at $X(t)=x$
		then forces the Hessian entries to satisfy
		$b_*^\top(W-\bar m)=x-\bar\mu$. Moreover, since $\operatorname{rank}(\Omega)=2$, the conditional Gaussian law has an affine support of dimension two. Symmetry of $\Omega$ gives $b_*^\top\Omega=0$, and the conditional mean satisfies $b_*^\top(m(x)-\bar m)=x-\bar\mu$. Consequently, this support is contained in the plane in \eqref{eq:singularheightplane}. Both are affine planes of dimension two, so they coincide. Finally, $\varepsilon_*=\Var(X(t)\mid\nabla X(t),W)$ is the Schur complement of $\Cov\bigl((\nabla X(t)^\top,W^\top)^\top\bigr)$ in $\Sigma(t)$, which gives the determinant formula.
	\end{proof}
	
	The condition $\widetilde V\succ0$ is essential to this completeness
	statement. It excludes lower-dimensional Hessian laws, such as
	separable sums $X(t_1,t_2)=Y_1(t_1)+Y_2(t_2)$ with identically zero
	mixed derivative. For a stationary planar field, $r=A=0$, so the
	singularity criterion simplifies to
	$\varepsilon_*=v-\iota^\top V^{-1}\iota=0$.
	
	\subsection{Peak height support geometry}
	\label{subsec:singularsupport}
	For a nonsingular conditional law, $J(t,x)>0$ at every height.
	For a singular law, the affine support must meet the negative definite
	cone. The next theorem turns this intersection into an exact support
	classification.
	
	Write
	$b_*=(b_{*,1},b_{*,2},b_{*,3})^\top$ and
	$\bar m=(\bar m_1,\bar m_2,\bar m_3)^\top$, and define
	\[
	\bar M=\begin{pmatrix}
		\bar m_1 & \bar m_3\\ \bar m_3 & \bar m_2
	\end{pmatrix},\qquad Q=-\begin{pmatrix}
		b_{*,1}&b_{*,3}/2\\ b_{*,3}/2&b_{*,2}
	\end{pmatrix},\qquad
	x_* =\bar\mu+\tr(Q\bar M)=\bar\mu-b_*^\top\bar m.
	\]
	Here $\bar M$ is the conditional mean of $\nabla^2X(t)$ given
	$\nabla X(t)=0$, and $Q$ is the symmetric matrix with
	$b_*^\top\vech H=-\tr(QH)$ for every symmetric $2\times2$ matrix $H$.
	With $H=\nabla^2X(t)$ and $W=\vech H$, the relation
	\eqref{eq:singularheightrelation} therefore reads
	$\tr(QH)=x_*-X(t)$ almost surely under the law conditional on
	$\nabla X(t)=0$, and the plane \eqref{eq:singularheightplane} is
	$\{H:\tr(QH)=x_*-x\}$.
	\begin{theorem}[Height supports under a linear Hessian constraint]
		\label{thm:singularsupport}
		Fix $t$ and suppose $\varepsilon_*=0$.  Then $Q\ne0$, and the conditional Hessian has full 
		support on $\{H:\tr(QH)=x_*-x\}$. Moreover,
		\begin{equation}\label{eq:generalheightpositivity}
			\{x:h_t(x)>0\}=
			\begin{cases}
				(x_*,\infty),&Q\succeq0,\\
				(-\infty,x_*),&Q\preceq0,\\
				\R,&Q\text{ is indefinite}.
			\end{cases}
		\end{equation}
		The height support is the closure of the corresponding set,
		with $h_t(x_*)=0$ in either semidefinite case.
		For minima, the two half lines are reversed; the saddle height
		density is strictly positive on $\R$.
	\end{theorem}
	\begin{proof}
		Since $b_*\ne0$ by Theorem~\ref{thm:singularrank}, $Q\ne0$. The support plane follows from
		\eqref{eq:singularheightplane}.
		The peak weight $\mathfrak D$ is positive exactly on the open negative
		definite cone. Thus $J(t,x)>0$ exactly when this plane contains a
		negative definite matrix. If $Q\succeq0$ is nonzero, then
		$\tr(QH)<0$ for every $H\prec0$. Conversely, every $q<0$ is attained by
		$H=(q/\tr Q)I_2\prec0$. Hence a peak Hessian is possible exactly
		when $x_*-x<0$, which proves the first case of
		\eqref{eq:generalheightpositivity}. Replacing $Q$ by $-Q$ proves
		the second. If $Q$ is indefinite,
		diagonalize it with eigenvalues $q_1>0>q_2$. The values
		$-q_1a_1-q_2a_2$, for $a_1,a_2>0$, cover $\R$, proving the third case.
		Positive definite Hessians give the assertions for minima.

		For saddles, work in the eigenbasis of $Q$ and choose
		$a_1,a_2$ with $q_1a_1+q_2a_2=x_*-x$.
		A matrix with these diagonal entries remains in the plane
		for every off diagonal entry $a_3$, and has negative
		determinant when $|a_3|$ is sufficiently large.
		Thus every plane meets the open saddle cone.
		The Gaussian height factor and normalizers are positive,
		giving the density statements.
	\end{proof}
	
	If the field $X$ is centered, then $x_*=0$. The classification
	admits only half-line and whole-line supports, so bounded peak height support
	cannot occur when $\widetilde V\succ0$ (see
	Remark~\ref{rem:degeneratehessian}). Both definite and rank-one
	semidefinite constraints give a half line, but their conditional peak weights
	differ, and the height density vanishes at $x_*$ to order $4$ in the
	first case and to order $3/2$ in the second
	(Proposition~\ref{prop:contact}). An indefinite constraint permits every height despite the
	singular Hessian covariance. Section~\ref{subsubsec:singularfamily}
	realizes all three constraint types in stationary planar fields; the
	quadratic field in Section~\ref{subsec:quadratic} realizes them as
	the spatial location varies.
	
	\begin{remark}[Degenerate Hessian laws]
		\label{rem:degeneratehessian}
		Assumption~\ref{ass:regularity} requires $\widetilde V\succ0$; the
		singular laws treated here are singular jet laws, not singular Hessian
		laws given the gradient. If $\widetilde V$ is singular,
		Theorem~\ref{thm:main} and the kernel formulas still apply, since they
		hold for every covariance $S\succeq0$, and conditioning on the height
		lowers the rank by at most one, so that
		$\operatorname{rank}\Omega\in\{\operatorname{rank}\widetilde V-1,
		\operatorname{rank}\widetilde V\}$. By convexity of the negative definite
		cone, the positivity set
		$\{x:(m(x)+\operatorname{range}\Omega)\cap\{H\prec0\}\ne\emptyset\}$ is
		still an open interval, but it may be empty or bounded. It is empty
		for a Gaussian harmonic field, where $\det H=-\norm\xi^2\le0$, and
		bounded in the following example. For
		$X(t)=\frac c4\norm t^2+\mathcal N_1\bigl(1+\frac{t_1^2-t_2^2}2\bigr)
		+\mathcal N_2t_1+\mathcal N_3t_2+\mathcal N_4t_1t_2$ with $c<0$ and
		independent standard normal $\mathcal N_j$, the trace is fixed at $c$;
		at $t=0$ the Hessian given $X(0)=x$ and $\nabla X(0)=0$ has rows
		$(c/2+x,\mathcal N_4)$ and $(\mathcal N_4,c/2-x)$, so the peak height
		support is the bounded interval $[-|c|/2,|c|/2]$. In general, the Morse
		property must then be verified separately, and continuity may fail
		when the normalizer vanishes in the limit; we do not pursue these
		cases.
	\end{remark}
	
	\paragraph{Laplace eigenfunction constraints and curvature.}
	Suppose $Z$ is centered and satisfies
	$\Delta_{\mathcal M}Z=-\omega_\Delta Z$ almost surely for a fixed $\omega_\Delta>0$,
	and let $X=\mu+Z$. The identity
	\begin{equation}\label{eq:generaleigenfunctionconstraint}
		\tr\nabla^2X=-\omega_\Delta X+(\Delta_{\mathcal M}+\omega_\Delta)\mu
		=-\omega_\Delta(X-x_*),\qquad
		x_* =\mu+\frac{\Delta_{\mathcal M}\mu}{\omega_\Delta},
	\end{equation}
	holds before conditioning. When Assumption~\ref{ass:regularity} holds,
	it gives the definite constraint $Q=I_2/\omega_\Delta$ in
	Theorem~\ref{thm:singularsupport}, with peak height support
	$[x_*,\infty)$. For a field of degree-$j$ spherical harmonics on the
	sphere of radius $R$, $\omega_\Delta=j(j+1)/R^2$, so curvature sets the scale of
	$Q=R^2I_2/\{j(j+1)\}$ but not its signature; Section~\ref{subsec:sphericaltrend}
	treats the quadrupole with a trend.
	
	\subsection{Singular jets and second-order equations}
	\label{subsec:singularpde}
	The constraint operator has an intrinsic meaning: a singular jet law
	means that, at the given point, the field satisfies a linear
	second-order relation whose principal part is $Q$. For a symmetric
	endomorphism $Q_t$ of $T_t\mathcal M$ and a vector $\beta_t\in T_t\mathcal M$,
	write
	\begin{equation}\label{eq:pointwisepde}
		\mathcal Q_tf=f(t)+\tr\bigl(Q_t\nabla^2f(t)\bigr)
		+\langle\beta_t,\nabla f(t)\rangle,\qquad f\in C^2(\mathcal M),
	\end{equation}
	where the trace is computed in any orthonormal frame. We call
	$\mathcal Q_t$ elliptic, parabolic or hyperbolic according as $Q_t$ is
	definite, has rank one, or is indefinite. This is the usual
	classification of second-order operators in two variables by their
	principal symbol $k\mapsto\langle Q_tk,k\rangle$.
	
	\begin{theorem}[Singular jets and second-order equations]
		\label{thm:singularpde}
		Suppose Assumption~\ref{ass:regularity} holds, and let
		$U\subset\mathcal M$ be open.
		\begin{enumerate}[label=(\roman*),leftmargin=*]
			\item The jet law at $t$ is singular if and only if $Z$ satisfies at $t$ the second-order relation
			$\mathcal Q_tZ=0$ almost surely, with $\mathcal Q_t$ as in
			\eqref{eq:pointwisepde}, for some symmetric $Q_t$ and some $\beta_t$. Then $Q_t\ne0$, the pair $(Q_t,\beta_t)$ is unique, $Q_t$ is
			the constraint operator of Theorem~\ref{thm:singularsupport}, and
			$x_*(t)=\mathcal Q_t\mu$.
			\item If the jet law is singular at every point of the open set $U$, then
			$t\mapsto(Q_t,\beta_t)$ is continuous on $U$, and almost surely the
			sample path of $Z$ solves the second-order equation
			$\mathcal Q_tZ=0$ for all $t\in U$. 
			\item At a point where the jet law is singular, the positivity set of
			$h_t$ is a half line if $\mathcal Q_t$ is elliptic or parabolic, and
			$\R$ if $\mathcal Q_t$ is hyperbolic.
		\end{enumerate}
	\end{theorem}
	\begin{proof}
		(i) Since $\Cov\bigl((\nabla X(t)^\top,W^\top)^\top\bigr)\succ0$, every nonzero vector in
		the kernel of $\Sigma(t)$ has a nonzero first coordinate, and the kernel
		has dimension at most one. A kernel vector $(1,\beta^\top,c^\top)^\top$
		is the same as an almost sure relation
		$Z(t)+\langle\beta,\nabla Z(t)\rangle+c^\top\vech\nabla^2Z(t)=0$, and
		$c^\top\vech H=\tr(Q_tH)$ for the symmetric matrix $Q_t$ with diagonal
		$(c_1,c_2)$ and off-diagonal entry $c_3/2$. This proves the equivalence
		and the uniqueness. If $Q_t=0$, the relation contradicts
		$\Cov\bigl((X(t),\nabla X(t)^\top)^\top\bigr)\succ0$. With $X=\mu+Z$, the relation reads
		$X(t)=\mathcal Q_t\mu-\tr(Q_t\nabla^2X(t))-\langle\beta_t,\nabla X(t)\rangle$,
		so that $X(t)=\mathcal Q_t\mu-\tr(Q_tH)$ almost surely under the law
		conditional on $\nabla X(t)=0$. Under this law, Theorems~\ref{thm:singularrank} and~\ref{thm:singularsupport}
		also give $X(t)=x_*(t)-\tr(QH)$ almost surely, with $Q$ the constraint
		operator. Since $W$ has the nondegenerate covariance $\widetilde V$, the
		coefficients of an affine function of $W$ are uniquely determined, and
		comparing the two representations gives $Q_t=Q$ and $x_*(t)=\mathcal Q_t\mu$.
		
		(ii) In a local orthonormal frame, the covariance jet $\Sigma(t)$ is
		continuous and, by (i), has rank five at every point of $U$. Its kernel is therefore a line that
		depends continuously on $t$, and so does the kernel vector normalized as
		in (i); this gives the continuity of $(Q_t,\beta_t)$ in a local
		orthonormal frame. For each $t$, $\mathcal Q_tZ=0$ almost surely; hence the relation holds simultaneously on a countable dense
		subset of $U$, and then on all of $U$, because $t\mapsto\mathcal Q_tZ$ is
		continuous for $C^2$ sample paths. 
		
		(iii) In dimension two, a nonzero symmetric $Q_t$ is semidefinite
		exactly when it is definite or has rank one, so (iii) restates
		Theorem~\ref{thm:singularsupport}.
	\end{proof}
	
	Part~(i) is a relation at a single point; part~(ii) turns it into an
	equation for the sample paths only on an open set where singularity
	persists, and the example after Corollary~\ref{cor:spectralconic},
	singular only along a line, shows that the first does not imply the
	second. Thus peak heights are confined to a half line where the field satisfies an elliptic or parabolic relation, and range over the whole real line where
	it satisfies a hyperbolic one; the endpoint $x_*=\mathcal Q_t\mu$ is the
	operator applied to the mean. The Laplace eigenfunction relation
	\eqref{eq:generaleigenfunctionconstraint} is the elliptic case
	$Q_t=\omega_\Delta^{-1}I_2$ and $\beta_t=0$, for which
	$\mathcal Q_tf=f(t)+\omega_\Delta^{-1}\Delta_{\mathcal M}f(t)$. The
	quadratic field of Section~\ref{subsec:quadratic} solves an equation of
	mixed type, and the ellipsoid field of
	Proposition~\ref{prop:ellipsoid}(ii) solves an elliptic equation whose
	principal part is proportional to the inverse second fundamental form.
	
	For stationary fields, the equation has constant coefficients and is determined by the support of the spectral measure.
	\begin{corollary}[Stationary fields and spectral conics]
		\label{cor:spectralconic}
		Let $Z$ be a centered stationary field on $\R^2$ that satisfies
		Assumption~\ref{ass:regularity}, and let $F$ be its spectral measure.
		The jet law is
		singular if and only if $F$ is concentrated on the conic
		$\mathcal C_Q=\{k\in\R^2:k^\top Qk=1\}$ for some symmetric $Q\ne0$. In
		that case, almost surely $Z+\tr(Q\nabla^2Z)=0$ on $\R^2$, $Q$ is the
		constraint operator at every point, and $\mathcal C_Q$ is an ellipse, a
		pair of parallel lines or a hyperbola according as $Q\succ0$,
		$Q\succeq0$ has rank one, or $Q$ is indefinite, the equation being
		correspondingly elliptic, parabolic or hyperbolic. Since $\mathcal C_Q$ is empty when
		$Q\preceq0$, negative semidefinite constraints do not occur.
	\end{corollary}
	\begin{proof}
		Stationarity gives $r=A=0$, so $\nabla Z(t)$ is uncorrelated with
		$(Z(t),\nabla^2Z(t))$, and the variance of $\mathcal Q_tZ$ equals
		$\Var\{Z(t)+\tr(Q_t\nabla^2Z(t))\}+\beta_t^\top G\beta_t$. Hence
		$\beta_t=0$ in Theorem~\ref{thm:singularpde}(i), and by stationarity and
		uniqueness $Q_t=Q$ does not depend on $t$. Write
		$Z(t)=\int_{\R^2}e^{\mathrm i\langle k,t\rangle}\,\mathcal Z(\dd k)$ with an
		orthogonal random measure $\mathcal Z$ of control measure $F$. Since $F$
		has finite fourth moments,
		$\tr(Q\nabla^2Z(t))=-\int_{\R^2}k^\top Qk\,e^{\mathrm i\langle k,t\rangle}\,
		\mathcal Z(\dd k)$, and therefore
		\[
		\Var\{Z(t)+\tr(Q\nabla^2Z(t))\}=\int_{\R^2}(1-k^\top Qk)^2\,F(\dd k),
		\]
		which vanishes exactly when $F(\R^2\setminus\mathcal C_Q)=0$.
		Theorem~\ref{thm:singularpde} gives the remaining statements, and
		diagonalizing $Q$ identifies the conic.
	\end{proof}
	After a linear change of coordinates $y$, the elliptic, hyperbolic and
	parabolic equations become $\Delta_yZ+Z=0$, the Klein--Gordon equation
	$\partial_{y_1}^2Z-\partial_{y_2}^2Z+Z=0$, and $\partial_{y_1}^2Z+Z=0$,
	whose solutions are exactly sinusoidal in $y_1$; the first case contains
	the random plane wave \citep{Berry1977}.
	
	\paragraph{Singular sets with empty interior.} The jet law may be
	singular on a set with empty interior, where only the pointwise relation
	of Theorem~\ref{thm:singularpde}(i) holds. If $Z$ is singular everywhere
	with operators $\mathcal Q_t$ and $Y$ is an independent centered field with
	$C^2$ paths, then $Z+Y$ is singular at $t$ exactly when $\mathcal Q_tY=0$
	almost surely, since $\ker(\Sigma_Z+\Sigma_Y)=\ker\Sigma_Z\cap\ker\Sigma_Y$.
	For $Z$ in Section~\ref{subsubsec:singularfamily} and
	$X(t)=Z(t)+t_1\mathcal N_0$, with $\mathcal N_0\sim\mathsf N(0,1)$ independent
	of $Z$, $\mathcal Q_t(t_1\mathcal N_0)=t_1\mathcal N_0$, so this happens
	exactly on the line $\{t_1=0\}$, and Assumption~\ref{ass:regularity}
	holds. Thus $h_t$ has the support of Theorem~\ref{thm:singularsupport} on
	the line and is positive on $\R$ off it,
	Theorem~\ref{thm:singularcontinuity} makes $t\mapsto h_t$ continuous
	across the line, and $X$ solves no second-order equation on any open set.

	\subsection{Canonical conditional peak weights}
	\label{subsec:canonicalweights}
	We now choose coordinates adapted to the constraint to evaluate conditional peak weights.
	For a Hessian $H$ conditional on height $x$ and zero gradient,
	with $\varepsilon_*=0$, choose an invertible matrix $P$ satisfying
	\begin{equation}\label{eq:sigma}
		(P^{-1})^\top QP^{-1}=\varsigma I_2,\quad
		\varsigma\diag(1,0),\quad\text{or}\quad\diag(1,-1),
	\end{equation}
	according as $Q$ is definite, has rank one, or is indefinite; here
	$\varsigma\in\{-1,1\}$. Such a matrix is obtained by diagonalizing $Q$
	and rescaling the coordinates corresponding to its nonzero eigenvalues.
	Set $M=PHP^\top$ and $\varpi_P=(\det P)^2$. By \eqref{eq:sigma}, the conditional Hessian support in Theorem \ref{thm:singularsupport}, $\tr(QH)=\tr[(P^{-1})^\top QP^{-1}M]=x_*-x$ becomes, respectively,
	\[
	\tr M=-\delta,\qquad M_{11}=-\delta,\qquad
	M_{11}-M_{22}=-\delta,
	\]
	where $\delta=\delta(x)=\varsigma(x-x_*)$, the signed distance from the endpoint, with $\varsigma=1$ in the indefinite case. Congruence preserves negative definiteness and gives
	$\det H=\det M/\varpi_P$.
	Conditional on height and zero gradient, $M$ has full Gaussian
	support on its two-dimensional constraint plane. Set
	$\tau_M=\tr M$ and $B_0=M_{12}$. The pair $((M_{11}-M_{22})/2,B_0)$ in the definite case~(i), and the pair
	$(\tau_M,B_0)$ in the rank-one and indefinite cases~(ii) and~(iii) of
	Theorem~\ref{thm:canonicalweights} below, are affine coordinates on the
	respective constraint planes. Their covariance matrices are therefore
	positive definite, which gives $\sigma_{\tau_M}^2>0$, $\lambda_{B_0}>0$
	and $\Lambda_M\succ0$ below.
	
	All moments in the following displays are taken conditional on
	$X(t)=x$ and $\nabla X(t)=0$. The means are affine functions of $x$;
	the covariance parameters depend on $t$ and $P$, but not on $x$.
	In cases (ii) and (iii), write
	\[
	\mu_{\tau_M}(x)=\E[\tau_M\mid X(t)=x,\nabla X(t)=0], \qquad \mu_{B_0}(x)=\E[B_0\mid X(t)=x,\nabla X(t)=0].
	\]
	Gaussian conditioning gives $\tau_M\sim\mathsf N(\mu_{\tau_M}(x),\sigma_{\tau_M}^2)$ with $\sigma_{\tau_M}^2>0$, and $B_0\mid\tau_M=z\sim\mathsf N\bigl(\beta_{B_0}(z;x),\lambda_{B_0}\bigr)$,
	where 
	\[
	\begin{split}
		\beta_{B_0}(z;x)&=\mu_{B_0}(x)
		+\Cov(B_0,\tau_M)
		\bigl(z-\mu_{\tau_M}(x)\bigr)/\sigma_{\tau_M}^2,\\
		\lambda_{B_0}&=\Var(B_0)
		-\Cov(B_0,\tau_M)^2/\sigma_{\tau_M}^2>0.
	\end{split}
	\]
	These parameters follow from \eqref{eq:WgivenHeightGradient} through
	$H\mapsto PHP^\top$, and only the means depend on $x$.
	
	\begin{theorem}[Conditional peak weights in canonical coordinates]
		\label{thm:canonicalweights}
		Fix $t$ and suppose $\varepsilon_*=0$. In the coordinates above, with $\delta=\delta(x)=\varsigma(x-x_*)$,
		one has $J(t,x)=\widehat J(x)/\varpi_P$, where:
		
		(i) If $Q$ is definite, so that $\tr M=-\delta$, let
		$\eta_M(x)$ and $\Lambda_M\succ0$ denote the mean and covariance of
		$((M_{11}-M_{22})/2,M_{12})^\top$, respectively. Then 
		\begin{equation}\label{eq:Q-definite}
			\widehat J(x)=\ind_{\{\delta>0\}}
			\Acal(\delta^2/4;\eta_M(x),\Lambda_M).
		\end{equation}
		
		(ii) If $Q$ has rank one, so that $M_{11}=-\delta$, then
		\begin{equation}\label{eq:rankonesimplifiedweight}
			\widehat J(x)=\ind_{\{\delta>0\}}
			\int_{-\infty}^{-\delta}
			\phin{z}{\mu_{\tau_M}(x)}{\sigma_{\tau_M}^2}
			\mathcal H\bigl(-\delta(\delta+z);
			\beta_{B_0}(z; x),\lambda_{B_0}\bigr)\dd z.
		\end{equation}
		
		(iii) If $Q$ is indefinite, so that $M_{11}-M_{22}=-\delta$, then
		\begin{equation}\label{eq:indefinitesimplifiedweight}
			\widehat J(x)=\int_{-\infty}^{-|\delta|}
			\phin{z}{\mu_{\tau_M}(x)}{\sigma_{\tau_M}^2}
			\mathcal H\!\left(\frac{z^2-\delta^2}{4};
			\beta_{B_0}(z; x),\lambda_{B_0}\right)\dd z.
		\end{equation}
	\end{theorem}
	\begin{proof}
		Congruence preserves negative definiteness and gives
		$\det M=\varpi_P\det H$, so it suffices to calculate the conditional peak
		weight for $M$. In case (i), its deterministic trace $-\delta$
		and Lemma~\ref{lem:cones} give the displayed norm kernel in \eqref{eq:Q-definite} directly.
		In (ii), $M\prec0$ requires $\delta>0$. Given $\tau_M=z$,
		$M_{22}=z+\delta$ and $\det M=-\delta(\delta+z)-B_0^2$.
		Thus, given also $\tau_M=z$, the conditional expectation of the peak weight of $M$ is
		$\mathcal H(-\delta(\delta+z);
		\beta_{B_0}(z; x),\lambda_{B_0})$. For $\delta>0$ it vanishes
		unless $z<-\delta$, proving \eqref{eq:rankonesimplifiedweight}.
		
		In (iii), $\det M=(\tau_M^2-\delta^2)/4-B_0^2$.
		For a negative trace, the peak weight is the positive part of
		this expression and vanishes unless $\tau_M<-|\delta|$.
		Conditioning on $\tau_M=z$ gives
		\eqref{eq:indefinitesimplifiedweight}.
	\end{proof}
	
	\begin{table}[!t]
		\centering
		\small
		\setlength{\tabcolsep}{4pt}
		\renewcommand{\arraystretch}{1.25}
		\caption{Support geometry and canonical constraints for $M=PHP^\top$.
		}
		\label{tab:singularsummary}
		\begin{tabular}{@{}lll@{}}
			\toprule
			$Q$ & $\{x:h_t(x)>0\}$ & Constraint on $M$ \\
			\midrule
			$Q\succ0$ & $(x_*,\infty)$ & $\tr M=-(x-x_*)$ \\
			$Q\prec0$ & $(-\infty,x_*)$ & $\tr M=x-x_*$ \\
			$Q\succeq0$, rank one & $(x_*,\infty)$ & $M_{11}=-(x-x_*)$ \\
			$Q\preceq0$, rank one & $(-\infty,x_*)$ & $M_{11}=x-x_*$ \\
			$Q$ indefinite & $\R$ & $M_{11}-M_{22}=x_*-x$ \\
			\bottomrule
		\end{tabular}
	\end{table}
	
\begin{figure}[!t]
		\centering
		\includegraphics[width=\textwidth]{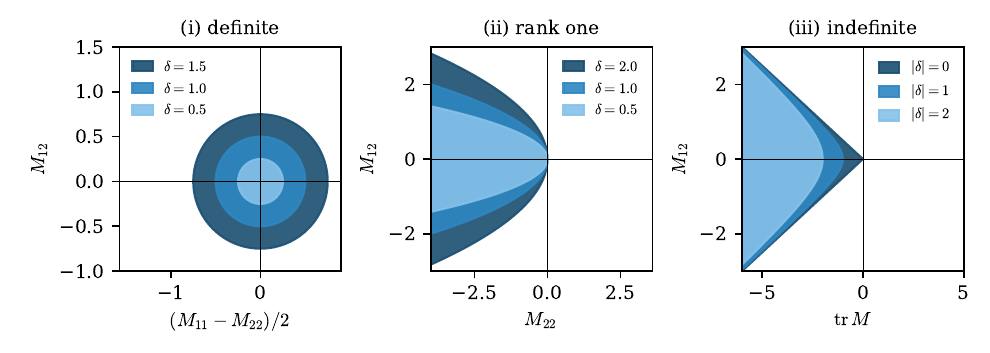}
		\caption{Sections of the open negative definite cone by the
			constraint planes of Theorem~\ref{thm:canonicalweights}, in
			canonical coordinates. (i) Definite, $\tr M=-\delta$: a disk of
			radius $\delta/2$ in the traceless coordinates, empty for
			$\delta\le0$. (ii) Rank one, $M_{11}=-\delta$: the region
			$\{M_{22}<0,\ M_{12}^2<-\delta M_{22}\}$, empty for $\delta\le0$.
			(iii) Indefinite, $M_{11}-M_{22}=-\delta$: the region
			$\{\tr M<-|\delta|,\ M_{12}^2<((\tr M)^2-\delta^2)/4\}$, nonempty
			for every $\delta$.}
		\label{fig:slices}
	\end{figure}

	The formulas follow the support geometry summarized in
	Table~\ref{tab:singularsummary}. A definite constraint becomes a fixed
	trace and requires one norm kernel evaluation. Rank one and indefinite
	constraints leave a random trace and require one scalar Gaussian
	integral of $\mathcal H$. Circular symmetry further simplifies the definite case. Figure~\ref{fig:slices} shows the sections of the negative definite
	cone by the three canonical constraint planes.
	
	\begin{corollary}[Closed formulas under circular symmetry]
		\label{cor:circularsingular}
		Suppose $\varepsilon_*=0$ and $Q$ is definite. Choose
		$\varsigma\in\{-1,1\}$ so that $Q_+=\varsigma Q\succ0$, and use
		the coordinates of Theorem~\ref{thm:canonicalweights}(i).
		Assume $\eta_M(x)=0$ and $\Lambda_M=\lambda I_2$, $\lambda>0$,
		for every height $x$. With $\delta=\varsigma(x-x_*)$, we have
		$\varpi_P=\det Q_+$ and
		\begin{align}\label{eq:circularsingularpeak}
			J(t,x)=\frac{\ind_{\{\delta>0\}}}{\det Q_+}
			\left[\frac{\delta^2}{4}-2\lambda+
			2\lambda e^{-\delta^2/(8\lambda)}\right],\qquad
			J_{\mathrm{sdl}}(t,x)=\frac{2\lambda}{\det Q_+}
			e^{-\delta^2/(8\lambda)}.
		\end{align}
	\end{corollary}
	\begin{proof}
		For a centered circular Gaussian law, \eqref{eq:centeredcircularkernel} gives $\Acal(a;0,\lambda I_2)
		=a-2\lambda+2\lambda e^{-a/(2\lambda)}$.
		Its complementary saddle moment is
		$2\lambda-a+\Acal(a;0,\lambda I_2)
		=2\lambda e^{-a/(2\lambda)}$.
		Applying Theorem~\ref{thm:canonicalweights}(i) at
		$a=\delta^2/4$ proves the first part in \eqref{eq:circularsingularpeak}; the
		complementary trace formula gives
		the second one.
	\end{proof}

	\subsection{Fields realizing the three constraint types}
	\label{subsubsec:singularfamily}
	The following finite Fourier family realizes all three canonical
	constraints in Theorem~\ref{thm:canonicalweights}. For $\gamma>0$,
	let $k_1=(1,0)^\top$, $k_2=(\gamma,1)^\top$, and
	$k_3=(\gamma,-1)^\top$, and define
	\[
	Z(t)=\frac1{\sqrt3}\sum_{j=1}^3
	\bigl[\mathcal{N}^{(\mathrm c)}_j\cos(k_j^\top t)+\mathcal{N}^{(\mathrm s)}_j\sin(k_j^\top t)\bigr],
	\qquad t\in\R^2,
	\]
	where the six coefficients are independent standard normal variables.
	This centered stationary field has unit variance and covariance
	\[
	C_0(h)=\frac13\bigl[\cos h_1+
	\cos(\gamma h_1+h_2)+\cos(\gamma h_1-h_2)\bigr].
	\]
	Its paths are analytic and its Hessian increment variances are
	$O(\norm{t-t'}^2)$ as $t'\to t$. Direct differentiation gives
	\[
	G=\diag\left(\frac{1+2\gamma^2}{3},\frac23\right),\quad \iota=\begin{pmatrix}
		-(1+2\gamma^2)/3\\
		-2/3\\
		0
	\end{pmatrix}, \quad
	V=\begin{pmatrix}
		(1+2\gamma^4)/3&2\gamma^2/3&0\\
		2\gamma^2/3&2/3&0\\
		0&0&2\gamma^2/3
	\end{pmatrix}.
	\]
	Since $V\succ0$, with $\det V=4\gamma^2/27$, and $r=A=0$, both
	conditions in \eqref{eq:non-degenerate} hold.
	Each frequency satisfies
	$k_{j,1}^2+(1-\gamma^2)k_{j,2}^2=1$, yielding the exact relation $Z_{11}+(1-\gamma^2)Z_{22}=-Z$.
	In the terms of Corollary~\ref{cor:spectralconic}, the spectral
	measure lives on the conic $k_1^2+(1-\gamma^2)k_2^2=1$, which is an
	ellipse for $\gamma<1$, the pair of lines $k_1=\pm1$ for $\gamma=1$ and a
	hyperbola for $\gamma>1$, and $Z$ solves
	$\partial_1^2Z+(1-\gamma^2)\partial_2^2Z+Z=0$, an equation of elliptic,
	parabolic or hyperbolic type, respectively.
	Consequently,
	\[
	b_*=V^{-1}\iota=(-1,\gamma^2-1,0)^\top,\qquad
	\varepsilon_*=0,\qquad Q=\diag(1,1-\gamma^2),\qquad
	\sigma_\tau^2=\frac{2\gamma^4}{9}>0.
	\]
	Thus $\Omega$ has rank two for every $\gamma>0$, and, with $X=Z$,
	\[
	\tau\mid\{X(t)=x,\nabla X(t)=0\}
	\sim\mathsf N\left(-\frac{3+2\gamma^2}{3}x,\frac{2\gamma^4}{9}\right).
	\]
	The peak intensity in every case is $\rho_\gamma(x)=
	3\phi(x)J_\gamma(x)/\{2\pi\sqrt{2(1+2\gamma^2)}\}$.
	Under the conditional Hessian law, the constraint becomes
	$H_{11}+(1-\gamma^2)H_{22}=-x$. We now specialize $\gamma$ to
	obtain definite, rank-one, and indefinite examples.
	
	\paragraph{Definite constraint: $\gamma=1/2$.}
	Here $Q=\diag(1,3/4)$ and the original conditional trace variance is
	$1/72$. In the canonical coordinates
	$M=Q^{1/2}HQ^{1/2}$, the trace is fixed at $-x$, and the
	traceless coordinates have mean zero and covariance $I_2/8$.
	Since $\varpi_P=\det Q=3/4$, Corollary~\ref{cor:circularsingular} with
	$\lambda=1/8$ and $\delta=x$ gives
	\begin{align*}
		J_{1/2}(x)
		=\frac43\ind_{\{x>0\}}\Bigl(\frac{x^2}{4}-\frac14+\frac14e^{-x^2}\Bigr)
		=\frac13\bigl(x^2-1+e^{-x^2}\bigr)\ind_{\{x>0\}}.
	\end{align*}
	Consequently, $\rho_{1/2}(x)=\ind_{\{x>0\}}\phi(x)
	(x^2-1+e^{-x^2})/(2\pi\sqrt3)$. The normalizer is $J=1/(6\sqrt3)$ and the total peak intensity is
	$\ell=1/(12\pi)$. Comparing with \eqref{eq:isotropicboundaryheight}, this is the
	random plane wave peak height law, and by
	Corollary~\ref{cor:circularsingular} the saddle height law is
	$\mathsf N(0,1/3)$.
	
	\paragraph{Constraint of rank one: $\gamma=1$.}
	The constraint is $H_{11}=-x$, and the trace law is
	$\mathsf N(-5x/3,2/9)$. Conditional on $\tau=z$, the traceless
	coordinates are $(-x-z/2,B_0)$, where $B_0\sim\mathsf N(0,2/3)$.
	Theorem~\ref{thm:canonicalweights}(ii), with $P=I_2$, gives
	\[
	J_1(x)=\ind_{\{x>0\}}\int_{-\infty}^{-x}
	\phin{z}{-5x/3}{2/9}
	\mathcal H\bigl(-x(x+z);0,2/3\bigr)\dd z,
	\]
	where \eqref{eq:Hcentralcdf} evaluates $\mathcal H$ through central chi
	square CDFs with one and three degrees of freedom.
	
\begin{figure}[htbp]
		\centering
		\includegraphics[width=\textwidth]{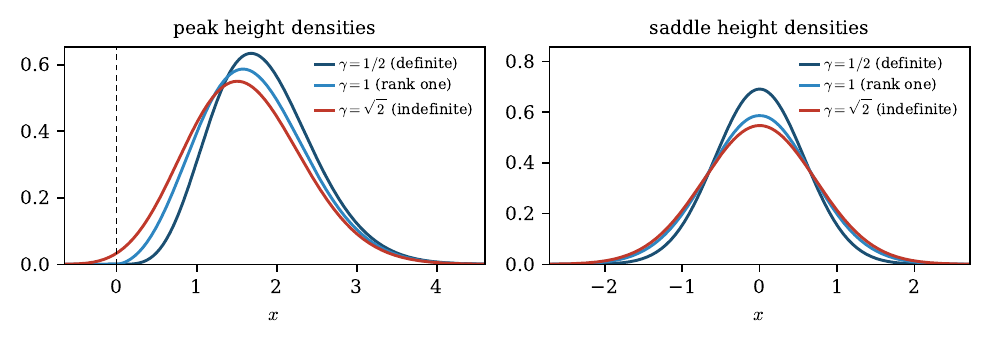}
		\caption{Peak height densities (left) and saddle height densities
			(right) of the Fourier family of
			Section~\ref{subsubsec:singularfamily} for the definite
			($\gamma=1/2$), rank-one ($\gamma=1$) and indefinite
			($\gamma=\sqrt2$) constraints, whose spectral measures lie on an
			ellipse, a pair of parallel lines and a hyperbola. The first two
			peak densities
			vanish on $(-\infty,0]$, whose endpoint is marked by the dashed
			vertical line; the indefinite one is positive on $\R$.}
		\label{fig:fourier}
	\end{figure}

	\paragraph{Indefinite constraint: $\gamma=\sqrt2$.}
	Now $H_{11}-H_{22}=-x$. The trace
	$\tau\sim\mathsf N(-7x/3,8/9)$ is independent of
	$B_0\sim\mathsf N(0,4/3)$, so Theorem~\ref{thm:canonicalweights}(iii)
	gives
	\[
	J_{\sqrt2}(x)=\int_{-\infty}^{-|x|}
	\phin{z}{-7x/3}{8/9}
	\mathcal H\!\left(\frac{z^2-x^2}{4};0,\frac43\right)\dd z.
	\]
	This conditional peak weight is positive at every height and uses the same two
	central CDFs in one trace integral. For both of the latter cases,
	the normalizer $J$ can be obtained through
	Corollary~\ref{cor:centerednormalizer},
	and the local height density is $\phi(x)J_\gamma(x)/J$.
	Numerically, $J\approx0.22683$ and $\ell\approx0.044215$ for
	$\gamma=1$, and $J\approx0.43144$ and $\ell\approx0.065142$ for
	$\gamma=\sqrt2$.
	Adding a deterministic $C^2$ trend preserves the covariance
	classification; the semidefinite cases have cutoff $x_*(t)$,
	while the indefinite case retains all real heights.
	Figure~\ref{fig:fourier} shows the three peak height densities and
	the corresponding saddle height densities.
	The saddle densities are $\phi(x)J_{\mathrm{sdl}}(x)/J_{\mathrm{sdl}}$,
	where $J_{\mathrm{sdl}}(x)=\Kcal_{\mathrm{sdl}}(m(x),\Omega)$ is computed from
	\eqref{eq:Ksaddle} with the rank-one $\Lambda$ and
	Proposition~\ref{prop:singularcdfnorm}, and $J_{\mathrm{sdl}}=2J$ by
	Proposition~\ref{prop:stationarycounts}; for $\gamma=1/2$ this is the
	$\mathsf N(0,1/3)$ law above.

	\section{Stability of peak height laws across covariance degeneracy}
	\label{sec:continuity}
	This section establishes the stability of local peak height laws with
	respect to the
	Gaussian regression parameters, uniformly on compact sets and through
	the covariance degeneracy of Section~\ref{sec:singular}.
	Section~\ref{subsec:singularconvergence} proves Lipschitz continuity in
	total variation and H\"older continuity of order $3/4$ in the uniform
	norm, and translates both into rates for mean and covariance jets and
	regional counts. Section~\ref{subsec:endpoint} describes how the height
	density vanishes at the endpoint of a half-line support and shows that
	the exponent $3/4$ is optimal.
	
	\subsection{Stability bounds and singular limits}
	\label{subsec:singularconvergence}
	Singular conditional laws may confine peak heights to a half line,
	whereas nonsingular laws give positive height densities on $\R$. We show that the height laws nevertheless depend Lipschitz continuously
	on the Gaussian regression parameters in total variation, and
	H\"older continuously of order $3/4$ in the uniform norm, even across
	singular limits. All jets at a fixed $t$ are expressed in one orthonormal frame;
	the conclusions do not depend on this choice. Throughout, for a function $g$ on $\R^k$ we write
	$\norm g_{L^1}=\int_{\R^k}|g|$, or $\norm g_{L^1(\R)}$ when $k=1$, and
	$\norm g_\infty=\sup|g|$. Total variation between densities $h,h'$ on
	$\R$ is measured by $\norm{h-h'}_{L^1(\R)}$, which is twice the total
	variation distance $\sup_B|\int_B(h-h')|$ between the corresponding laws. We begin with the
	expected peak weight itself.
	
	\begin{lemma}[Coupling bound for the expected peak weight]
		\label{lem:kernelcoupling}
		The map $(m,S)\mapsto\Kcal(m,S)$ is continuous on
		$\R^3\times\{S:S\succeq0\}$. More precisely, for
		$m,m'\in\R^3$ and $S,S'\succeq0$, define
		\[
		\mathcal R_2(m,S)=\sqrt{\norm m^2+\tr S},\qquad
		\mathcal D=\sqrt{\norm{m-m'}^2+
			\norm{S^{1/2}-{S'}^{1/2}}_{\mathrm F}^2},
		\]
		where $S^{1/2}$ is the positive semidefinite square root and
		$\norm{\cdot}_{\mathrm F}$ is the Frobenius norm. Then
		\begin{equation}\label{eq:kernelcouplingbound}
			|\Kcal(m,S)-\Kcal(m',S')|
			\le2\{\mathcal R_2(m,S)+\mathcal R_2(m',S')\}\mathcal D.
		\end{equation}
	\end{lemma}
	
	\begin{proof}
		For $w=(w_1,w_2,w_3)^\top$, recall from Section~\ref{subsec:reduction} that
		\begin{equation}\label{eq:F-max}
			H(w)=\begin{pmatrix}w_1&w_3\\w_3&w_2\end{pmatrix},\qquad
			\mathfrak D(w)=|\det H(w)|\ind_{\{H(w)\prec0\}}=(\lambda_{\min}(w))_-(\lambda_{\max}(w))_-,
		\end{equation}
		where $\lambda_{\min}(w)\le\lambda_{\max}(w)$ are the ordered eigenvalues of $H(w)$
		and $a_-=\max\{-a,0\}$. 
		This identity includes the boundary of the negative definite cone,
		where the peak weight is zero. The variational characterization of the
		eigenvalues gives, for $i\in\{\min,\max\}$,
		\[
		|\lambda_i(w)-\lambda_i(w')|
		\le\norm{H(w)-H(w')}_{\mathrm F}
		\le\sqrt2\norm{w-w'},\qquad
		|\lambda_i(w)|\le\sqrt2\norm w.
		\]
		Since $|a_--b_-|\le|a-b|$, taking the difference of the two products
		yields $|\mathfrak D(w)-\mathfrak D(w')|
		\le2(\norm w+\norm{w'})\norm{w-w'}$.
		Also,
		\[
		0\le \mathfrak D(w)
		\le|w_1w_2-w_3^2|
		\le\frac{w_1^2+w_2^2}{2}+w_3^2
		\le\norm w^2.
		\]
		
		To compare two Gaussian laws, use a common Gaussian vector
		$\mathcal N\sim\mathsf N_3(0,I_3)$ and set $\mathcal W=m+S^{1/2}\mathcal N$ and $\mathcal W'=m'+{S'}^{1/2}\mathcal N$.
		This construction remains valid for singular covariances,
		and $\E\norm {\mathcal W}^2=\mathcal R_2(m,S)^2$, $\E\norm{\mathcal W'}^2=\mathcal R_2(m',S')^2$ and $\E\norm{\mathcal W-\mathcal W'}^2=\mathcal D^2$.
		The difference bound, followed by Cauchy--Schwarz and the
		triangle inequality in $L^2$, now gives
		\begin{align*}
			|\Kcal(m,S)-\Kcal(m',S')|&=\left|\E[\mathfrak D(\mathcal W)-\mathfrak D(\mathcal W')]\right|
			\le2\E\big[(\norm {\mathcal W}+\norm{\mathcal W'})\norm{\mathcal W-\mathcal W'}\big]\\
			&\le2\{\mathcal R_2(m,S)+\mathcal R_2(m',S')\}\mathcal D.
		\end{align*}
		This proves \eqref{eq:kernelcouplingbound}. The positive semidefinite square
		root is continuous in the matrix entries, so this bound proves continuity
		of $\Kcal$, including at singular covariance matrices.
	\end{proof}
	
	At singular matrices the square root is only $1/2$-H\"older continuous,
	so \eqref{eq:kernelcouplingbound} is of order $\norm{S-S'}^{1/2}$, and this
	exponent cannot be improved for $\Kcal$ over unrestricted covariances: if
	$H(m)=\diag(-1,0)$, then
	$\Kcal(m,0)=0$, whereas $\Kcal(m,\epsilon I_3)=\sqrt{\epsilon/(2\pi)}+O(\epsilon)$.
	Taking the infimum over all couplings of the two laws shows that $\mathcal D$
	may be replaced by the Wasserstein distance between $\mathsf N_3(m,S)$ and
	$\mathsf N_3(m',S')$. The example uses the covariance $S=0$,
	whereas the conditional Hessian covariances $\Omega$ in
	Theorem~\ref{thm:singularcontinuity} have rank at least two by
	Theorem~\ref{thm:singularrank}. Better rates for height laws must use this
	structure, which the following representation exposes; with it,
	Theorem~\ref{thm:singularcontinuity} obtains a linear rate in total
	variation and the uniform exponent $3/4$, which
	Proposition~\ref{prop:sharprate} shows to be optimal.
	
	For regression parameters
	$\theta=(\bar\mu,s^2,\bar m,\widetilde\iota,\widetilde V)$ as in
	Section~\ref{sec:regression}, write $J_\theta(x)=\Kcal(m(x),\Omega)$,
	$J_\theta=\Kcal(\bar m,\widetilde V)$ and $h_\theta$ for the height density
	\eqref{eq:peak-height}. Let $\Theta$ be the set of parameters with $s^2>0$,
	$\widetilde V\succ0$ and
	$\varepsilon_*=s^2-\widetilde\iota^\top\widetilde V^{-1}\widetilde\iota\ge0$,
	equipped with the Euclidean norm of the entries, and let
	$\mathfrak D$ be the peak weight of Section~\ref{subsec:reduction}, written through eigenvalues in \eqref{eq:F-max}.
	
	\begin{proposition}[Hessian-tilted representation of peak heights]
		\label{prop:tiltedrep}
		Let $\theta\in\Theta$, $b_*=\widetilde V^{-1}\widetilde\iota$ and
		$\mathcal L_\theta(w)=\bar\mu+b_*^\top(w-\bar m)$, and let $\nu_\theta$ be the
		probability measure on $\R^3$ with density
		\begin{equation}\label{eq:tiltedlaw}
			f_\theta(w)=\mathfrak D(w)\,\varphi_{\widetilde V}(w-\bar m)/J_\theta,
		\end{equation}
		where $\varphi_{\widetilde V}$ is the $\mathsf N_3(0,\widetilde V)$ density.
		Then $h_\theta$ is the density of
		\begin{equation}\label{eq:tiltedrep}
			Y_\theta=\mathcal L_\theta(W^\circ)+\sqrt{\varepsilon_*}\,\mathcal N,\qquad
			W^\circ\sim\nu_\theta,\quad\mathcal N\sim\mathsf N(0,1)\ \text{independent}.
		\end{equation}
		In particular, $h_\theta$ is real-analytic if $\varepsilon_*>0$. If
		$\varepsilon_*=0$, the support of $h_\theta$ is the closure of
		$\{x_*-\tr(QH):H\preceq0\}$, with $Q$ and $x_*$ as in
		Theorem~\ref{thm:singularsupport}.
	\end{proposition}
	\begin{proof}
		Let $(X_0,W^\top)^\top$ be Gaussian with mean $(\bar\mu,\bar m^\top)^\top$
		and covariance with blocks $s^2$, $\widetilde\iota$ and $\widetilde V$, which
		is positive semidefinite because $\varepsilon_*\ge0$; by
		\eqref{eq:Wgradient}, this is the law of $(X(t),W^\top)^\top$ given
		$\nabla X(t)=0$ when $\theta$ comes from a field. Then
		$W\sim\mathsf N_3(\bar m,\widetilde V)$ and, by the regression in the
		proof of Theorem~\ref{thm:singularrank},
		$X_0\mid W\sim\mathsf N(\mathcal L_\theta(W),\varepsilon_*)$; that is,
		$X_0\overset{\mathrm d}{=}\mathcal L_\theta(W)+\sqrt{\varepsilon_*}\,\mathcal N$
		with $\mathcal N\sim\mathsf N(0,1)$ independent of $W$. Fix a bounded
		measurable $g$ and set $I_g=\E\bigl[\mathfrak D(W)\,g(X_0)\bigr]$,
		which is finite because $g$ is bounded and $\mathfrak D(w)\le\norm w^2$. We
		compute $I_g$ in two ways, interchanging integrals by Fubini's theorem. Conditioning first on $W$ and using \eqref{eq:tiltedlaw}, we obtain
		\begin{align*}
			I_g&=\int_{\R^3}\mathfrak D(w)\,\varphi_{\widetilde V}(w-\bar m)\,
			\E\bigl[g\bigl(\mathcal L_\theta(w)+\sqrt{\varepsilon_*}\,\mathcal N\bigr)\bigr]\dd w\\
			&=J_\theta\int_{\R^3}f_\theta(w)\,
			\E\bigl[g\bigl(\mathcal L_\theta(w)+\sqrt{\varepsilon_*}\,\mathcal N\bigr)\bigr]\dd w
			=J_\theta\,\E[g(Y_\theta)],
		\end{align*}
		where the inner expectations are over $\mathcal N$ only and $Y_\theta$ is
		as in \eqref{eq:tiltedrep}. Conditioning first on $X_0$ instead, whose
		density is $\phin{x}{\bar\mu}{s^2}$, we have
		$\E[\mathfrak D(W)\mid X_0=x]=\Kcal(m(x),\Omega)=J_\theta(x)$
		by the regression leading to \eqref{eq:WgivenHeightGradient}, and hence, by \eqref{eq:peak-height},
		\[
		I_g=\int_\R g(x)\,\phin{x}{\bar\mu}{s^2}\,J_\theta(x)\dd x
		=J_\theta\int_\R g(x)\,h_\theta(x)\dd x .
		\]
		Since $J_\theta>0$, the two expressions give
		$\E[g(Y_\theta)]=\int_\R g\,h_\theta$ for every bounded measurable $g$, so
		$Y_\theta$ has density $h_\theta$.
		
		If $\varepsilon_*>0$, $h_\theta$ is the convolution of a probability
		measure with a nondegenerate Gaussian density. If $\varepsilon_*=0$,
		$h_\theta$ is the density of the image of $\nu_\theta$ under
		$\mathcal L_\theta$; the support of $\nu_\theta$ is the closed negative
		semidefinite cone, and $\mathcal L_\theta(\vech H)=x_*-\tr(QH)$ since $\tr(QH)=-b_*^\top\vech H$ and
		$x_*=\bar\mu-b_*^\top\bar m$.
	\end{proof}
	Thus a singular law is precisely one in which the Gaussian smoothing in
	\eqref{eq:tiltedrep} disappears, and the support statement of
	Theorem~\ref{thm:singularsupport} expresses whether $\tr(QH)$ has constant
	sign on $\{H\preceq0\}$. The tilted law $\nu_\theta$ depends smoothly on
	$(\bar m,\widetilde V)$ because its Gaussian factor stays nondegenerate.
	
	\begin{theorem}[Stability of peak height laws]
		\label{thm:singularcontinuity}
		Let $\Theta_0\subset\Theta$ be compact. There is a constant $L_0$ such that,
		for all $\theta,\theta'\in\Theta_0$,
		\begin{equation}\label{eq:lipschitzTV}
			\norm{h_\theta-h_{\theta'}}_{L^1(\R)}+|J_\theta-J_{\theta'}|
			\le L_0\norm{\theta-\theta'},\qquad
			\norm{h_\theta-h_{\theta'}}_\infty\le L_0\norm{\theta-\theta'}^{3/4}.
		\end{equation}
		In particular, if $\theta_n\to\theta$ in $\Theta$, then $h_{\theta_n}\to h_\theta$
		in total variation at rate $O(\norm{\theta_n-\theta})$ and uniformly at rate
		$O(\norm{\theta_n-\theta}^{3/4})$, also when $\Omega_n\succ0$ and $\Omega$
		is singular. The exponents~$1$ and~$3/4$ in \eqref{eq:lipschitzTV} are both optimal,
		by Remark~\ref{rem:optimalityscope} and Proposition~\ref{prop:sharprate}, respectively.
	\end{theorem}
	
	\begin{proof}
		Constants depend only on $\Theta_0$. By compactness and
		Lemma~\ref{lem:kernelcoupling}, there is $c_0>0$ such that $s^2\ge c_0$,
		$\widetilde V\succeq c_0I_3$, $J_\theta\ge c_0$ and $\norm\theta\le c_0^{-1}$
		on $\Theta_0$. The maps $\theta\mapsto\varepsilon_*$, $b_*$, $m(x)$,
		$\Omega$ are Lipschitz on $\Theta_0$, with $\norm{m_\theta(x)-m_{\theta'}(x)}
		\le c(1+|x|)\norm{\theta-\theta'}$.
		
		\emph{Step 1: regularity of the tilted density.} By
		Lemma~\ref{lem:cones}, $\mathfrak D=p\,\ind_{\mathcal S_-}$, where
		$p(w)=w_1w_2-w_3^2=\det H(w)$ and $\mathcal S_-=\{w:H(w)\prec0\}$ is, in the
		coordinates $Lw=(\tau,\xi)$, the right circular cone
		$\{\tau<0,\norm\xi<-\tau/2\}$. Since $p$ vanishes on $\partial\mathcal S_-$,
		$\mathfrak D$ is locally Lipschitz with $\nabla \mathfrak D=\ind_{\mathcal S_-}\nabla p$, and for
		every unit vector $e$, in the sense of distributions,
		\[
		\partial_e^2\mathfrak D=\ind_{\mathcal S_-}\partial_e^2p
		-(\partial_ep)(e^\top n)\,\sigma_\partial,
		\]
		with $n$ the outward unit normal and $\sigma_\partial$ the surface measure on
		$\partial\mathcal S_-$. Since $0\le \mathfrak D(w)\le\norm w^2$,
		$\norm{\nabla \mathfrak D(w)}\le2\norm w$ and $\norm{\bar m}\le c_0^{-1}$, the
		tilted density satisfies
		$f_\theta+\norm{\nabla f_\theta}\le c(1+\norm w)^3e^{-c'\norm w^2}$.
		With $\sigma_\partial(\{\norm w\le R\})\le cR^2$, this gives, uniformly over
		$\theta,\theta'\in\Theta_0$ and unit vectors $e$,
		\begin{equation}\label{eq:tiltedregularity}
			\int_{\R^3}(1+\norm w)\bigl(f_\theta+\norm{\nabla f_\theta}\bigr)\dd w\le c,
			\qquad |\partial_e^2f_\theta|(\R^3)\le c,\qquad
			\int_{\R^3}|f_\theta-f_{\theta'}|\dd w\le c\norm{\theta-\theta'}.
		\end{equation}
		The first bound follows from the Gaussian envelope above; the
		affine transport estimate in Step~4 uses this weighted gradient bound.
		The last bound holds because $f_\theta$ depends on $\theta$ only through
		$(\bar m,\widetilde V)$ and $J_\theta$, smoothly and with Gaussian bounds;
		the same argument gives $|J_\theta-J_{\theta'}|\le c\norm{\theta-\theta'}$.
		
		\emph{Step 2: reduction.} Write $\varepsilon,b$ and $\varepsilon',b'$ for
		$\varepsilon_*,b_*$ at $\theta,\theta'$. Pushforwards and independent noise do
		not increase total variation, so by \eqref{eq:tiltedrep} and
		\eqref{eq:tiltedregularity} it suffices to compare the laws of
		$\mathcal L_\theta(W^\circ)+\sqrt\varepsilon\mathcal N$ and
		$\mathcal L_{\theta'}(W^\circ)+\sqrt{\varepsilon'}\mathcal N$, where $W^\circ$ has density
		$f=f_{\theta'}$. If $\varepsilon_*\le3c_0/4$, then
		$\widetilde\iota^\top\widetilde V^{-1}\widetilde\iota\ge c_0/4$, so
		$\norm{b_*}\ge b_0:=c_0/2$. As $\varepsilon_*$ is Lipschitz, there is
		$\bar c>0$ such that $\norm{\theta-\theta'}<\bar c$ implies
		$\varepsilon,\varepsilon'\ge c_0/4$ or $\norm b,\norm{b'}\ge b_0$; larger
		distances are trivial.
		
		\emph{Step 3: smoothed regime.} If $\varepsilon,\varepsilon'\ge c_0/4$, the
		two laws are mixtures over $W^\circ$ of $\mathsf N(\mathcal L_\theta(W^\circ),\varepsilon)$
		and $\mathsf N(\mathcal L_{\theta'}(W^\circ),\varepsilon')$, whose total variation
		distance is at most $c(|\mathcal L_\theta(W^\circ)-\mathcal L_{\theta'}(W^\circ)|+|\varepsilon-\varepsilon'|)
		\le c(1+\norm{W^\circ})\norm{\theta-\theta'}$; \eqref{eq:tiltedregularity} bounds
		its expectation.
		
		\emph{Step 4: near-singular regime.} Let $\norm b,\norm{b'}\ge b_0$, and
		set $\check b=b/\norm b^2$, $\check b'=b'/\norm{b'}^2$, so that
		$\mathcal L_\theta(w+\alpha\check b)=\mathcal L_\theta(w)+\alpha$ for $\alpha\in\R$, and similarly for $\theta'$.
		Then $\mathcal L_\theta(W^\circ)+\sqrt\varepsilon\mathcal N=\mathcal L_\theta(\mathcal U)$
		with $\mathcal U=W^\circ+\sqrt\varepsilon\mathcal N\check b$ of density
		$f*\mathcal G^{\check b}_\varepsilon$, where $\mathcal G^{\check b}_r$ is the law of
		$\sqrt r\mathcal N\check b$; likewise with $\mathcal U'=W^\circ+\sqrt{\varepsilon'}\mathcal N\check b'$.
		Let $\Pi=I_3+b'(b-b')^\top/\norm{b'}^2$, so that $\Pi^\top b'=b$, and choose $w_0$
		parallel to $b'$ with $\mathcal L_{\theta'}(\Pi w+w_0)=\mathcal L_\theta(w)$;
		then $\norm{\Pi-I_3}+\norm{w_0}\le c\norm{\theta-\theta'}$, and it suffices to
		bound
		\[
		\norm{(f*\mathcal G^{\check b}_\varepsilon)_{\Pi,w_0}-f*\mathcal G^{\check b}_\varepsilon}_{L^1}
		+\norm{f*\mathcal G^{\check b}_\varepsilon-f*\mathcal G^{\check b'}_\varepsilon}_{L^1}
		+\norm{f*\mathcal G^{\check b'}_\varepsilon-f*\mathcal G^{\check b'}_{\varepsilon'}}_{L^1},
		\]
		where $g_{\Pi,w_0}$ is the density of $\Pi\mathcal U+w_0$ when $\mathcal U$ has
		density $g$. Since $1+\norm{w+v}\le(1+\norm w)(1+\norm v)$,
		$\norm{\check b}\le1/b_0$ and $\varepsilon\le c_0^{-1}$, the density
		$g=f*\mathcal G^{\check b}_\varepsilon$ inherits the weighted bounds in
		\eqref{eq:tiltedregularity}, the one on $\norm{\nabla g}$ because
		$\nabla g=(\nabla f)*\mathcal G^{\check b}_\varepsilon$. Shrinking $\bar c$ so that
		$\norm{\Pi-I_3}\le1/2$, the densities $g_\alpha$ of $A_\alpha(\mathcal U)$, with
		$A_\alpha(w)=w+\alpha\{(\Pi-I_3)w+w_0\}$, satisfy the same bound and solve
		$\partial_\alpha g_\alpha+\nabla\cdot(g_\alpha v_\alpha)=0$, where
		$\norm{v_\alpha(y)}\le c\{\norm{\Pi-I_3}(1+\norm y)+\norm{w_0}\}$ and
		$|\nabla\cdot v_\alpha|\le c\norm{\Pi-I_3}$; integrating over $\alpha\in[0,1]$ bounds
		the first term by $c(\norm{\Pi-I_3}+\norm{w_0})$. The second term is at most
		$\sqrt\varepsilon\,\E|\mathcal N|\,\norm{\check b-\check b'}\int\norm{\nabla f}
		\le c\norm{\theta-\theta'}$. For the third, $r\mapsto f*\mathcal G^{\check b'}_r$
		solves $\partial_rg=\frac12\partial_{\check b'}^2g$, so with $e'=b'/\norm{b'}$ it is at
		most $\frac12|\varepsilon-\varepsilon'|\,\norm{\check b'}^2|\partial_{e'}^2f|(\R^3)
		\le c\norm{\theta-\theta'}$. The residual variance thus enters linearly,
		also when $\varepsilon'=0$, which gives the linear rate.
		
		\emph{Step 5: uniform bound.} Since
		$h_\theta(x)=\phin{x}{\bar\mu}{s^2}J_\theta(x)/J_\theta$, with $J_\theta\ge c_0$,
		$J_\theta(x)\le c(1+x^2)$, and the Gaussian factor and $J_\theta$ Lipschitz in
		$\theta$, it suffices to show that, for $\norm{\theta-\theta'}<\bar c$,
		\begin{equation}\label{eq:Jholder}
			|J_\theta(x)-J_{\theta'}(x)|\le c(1+|x|)^{5/2}\norm{\theta-\theta'}^{3/4};
		\end{equation}
		larger distances are trivial because $h_\theta\le c$. As in Step~2, either
		$\varepsilon,\varepsilon'\ge c_0/4$ or $\varepsilon,\varepsilon'\le3c_0/4$.
		
		In the first case, with $u_*=\widetilde V^{-1/2}\widetilde\iota/s$ as in the
		proof of Theorem~\ref{thm:singularrank}, the eigenvalues $1,1,\varepsilon_*/s^2$
		of $I_3-u_*u_*^\top$ give $\Omega,\Omega'\succeq c_0^3I_3/4$. On
		$\{S\succeq cI_3\}$ the square root is Lipschitz,
		$\norm{S^{1/2}-S'^{1/2}}_{\mathrm F}\le\norm{S-S'}_{\mathrm F}/(2\sqrt c)$, as
		one sees by writing $S^{1/2}-S'^{1/2}$ in eigenbases of $S$ and $S'$; hence
		\eqref{eq:kernelcouplingbound} gives \eqref{eq:Jholder} with exponent $1$.
		
		In the second case, $\norm{u_*}^2=\widetilde\iota^\top\widetilde V^{-1}\widetilde\iota/s^2\ge c_0^2/4$
		at $\theta$ and $\theta'$. Let $\hat u=u_*/\norm{u_*}$,
		$e_*=\widetilde V^{1/2}\hat u$ and $\varkappa=(\varepsilon_*/s^2)^{1/2}\in[0,1]$.
		Since $I_3-u_*u_*^\top=(I_3-\hat u\hat u^\top)+\varkappa^2\hat u\hat u^\top$,
		the proof of Theorem~\ref{thm:singularrank} gives
		$\Omega=\widetilde V^{1/2}(I_3-\hat u\hat u^\top)\widetilde V^{1/2}+\varkappa^2e_*e_*^\top$.
		Let $\mathcal V\sim\mathsf N_3(0,I_3)$ and $\mathcal N\sim\mathsf N(0,1)$ be
		independent, and set
		\[
		\mathcal Y_\theta=m(x)+\widetilde V^{1/2}(I_3-\hat u\hat u^\top)\mathcal V,\qquad
		\Psi_\theta(r)=\E[\mathfrak D(\mathcal Y_\theta+re_*)],\quad r\in\R.
		\]
		Then $\mathcal Y_\theta+\varkappa\mathcal Ne_*\sim\mathsf N_3(m(x),\Omega)$, so that
		$J_\theta(x)=\E[\Psi_\theta(\varkappa\mathcal N)]$. The vector $\mathcal Y_\theta$
		lives on the plane $\mathcal P=m(x)+\widetilde V^{1/2}\hat u^\perp$, where its density
		is at most $(2\pi c_0)^{-1}\exp\{-c_0\norm{y-m(x)}^2/2\}$, because the nonzero
		eigenvalues of its covariance lie in $[c_0,c_0^{-1}]$. The vector $e_*$ is
		transversal to $\mathcal P$: the unit normal
		$n_{\mathcal P}=\widetilde V^{-1/2}\hat u/\norm{\widetilde V^{-1/2}\hat u}$ satisfies
		$n_{\mathcal P}^\top e_*\ge c_0^{1/2}$. Finally, $m(x)$,
		$\widetilde V^{1/2}(I_3-\hat u\hat u^\top)$, $e_*$ and $\varkappa^2$ are
		Lipschitz in $\theta$, the first with constant $c(1+|x|)$.
		
		We claim that, for every such $\theta$, $\Psi=\Psi_\theta$ is continuously
		differentiable and
		\begin{equation}\label{eq:Psiholder}
			|\Psi'(r)-\Psi'(r')|\le c(1+|x|+|r|+|r'|)^{5/2}|r-r'|^{1/2},\qquad r,r'\in\R.
		\end{equation}
		By Step~1, $\mathfrak D$ is locally Lipschitz with
		$\nabla \mathfrak D=\ind_{\mathcal S_-}\nabla p$ off $\partial\mathcal S_-$, and
		$W_r=\mathcal Y_\theta+re_*$ lies on $\partial\mathcal S_-$ with probability
		zero, so dominated convergence gives
		$\Psi'(r)=\E[\ind_{\mathcal S_-}(W_r)\,\partial_{e_*}p(W_r)]$, which is
		continuous in $r$. Since $\partial_{e_*}p$ is linear with
		$|\partial_{e_*}p(w)|\le2\norm{e_*}\norm w$,
		\[
		|\Psi'(r)-\Psi'(r')|\le2\norm{e_*}^2|r-r'|
		+2\norm{e_*}\,\E\bigl[\norm{W_{r'}}\,|\ind_{\mathcal S_-}(W_r)-\ind_{\mathcal S_-}(W_{r'})|\bigr].
		\]
		If the indicators differ, the segment from $W_{r'}$ to $W_r$ meets
		$\partial\mathcal S_-$, so $\mathcal Y_\theta$ lies in the projection onto
		$\mathcal P$ along $e_*$ of $\partial\mathcal S_-\cap\mathcal K$, where $\mathcal K$
		is the slab between $\mathcal P+r'e_*$ and $\mathcal P+re_*$, of width
		$\eta=n_{\mathcal P}^\top e_*|r-r'|$. In the coordinates $(\tau,\xi)=Lw$,
		$\partial\mathcal S_-=\{\varrho(-2,\cos\phi,\sin\phi):\varrho\ge0,\ \phi\in[0,2\pi)\}$
		is a round cone with area element $\sqrt5\,\varrho\dd\varrho\dd\phi$, and a slab
		$\{|n_L^\top(\tau,\xi)-d_L|\le\eta_L\}$ with unit normal $n_L$ meets it where
		$|\varrho g(\phi)-d_L|\le\eta_L$, with $g(\phi)=g_0+g_1\cos(\phi-\phi_0)$ and
		$\max\{|g_0|/2,|g_1|\}\ge2^{-1/2}$. Hence
		$|\{\phi:|g(\phi)|\le\lambda\}|\le c\lambda^{1/2}$ for all $\lambda>0$: if
		$|g_1|\ge1/4$, this follows from $g'^2+g''^2=g_1^2$ and the van der Corput
		sublevel estimate on the arcs where $|g'|$ or $|g''|$ is at least
		$|g_1|/\sqrt2$, and otherwise $|g|>1$. For fixed $\phi$, the admissible
		$\varrho\le R$ form an interval of length at most
		$\min\{R,2\eta_L/|g(\phi)|\}$, so the part of the cone in the slab with
		$\varrho\le R$ has area at most
		\[
		cR\int_0^{2\pi}\min\{R,2\eta_L/|g(\phi)|\}\dd\phi
		=cR\int_0^R\bigl|\{\phi:|g(\phi)|<2\eta_L/a\}\bigr|\dd a
		\le cR^{3/2}\eta_L^{1/2}.
		\]
		Projection along $e_*$ multiplies areas by at most $\norm{e_*}/n_{\mathcal P}^\top e_*$,
		and $\eta_L\le c\eta$. Splitting according to
		$\norm{\mathcal Y_\theta-m(x)}\in[2^{j-1}R_0,2^jR_0)$, with
		$R_0=1+\norm{m(x)}+|r|+|r'|$, and using the Gaussian bound on the density of
		$\mathcal Y_\theta$, we obtain
		$\E[\norm{W_{r'}}\,|\ind_{\mathcal S_-}(W_r)-\ind_{\mathcal S_-}(W_{r'})|]
		\le cR_0^{5/2}|r-r'|^{1/2}$, which proves \eqref{eq:Psiholder}.
		
		Now let $0\le\varkappa'<\varkappa\le1$ and $\bar\varkappa=(\varkappa^2-\varkappa'^2)^{1/2}$; the bound below is symmetric in $\varkappa$ and $\varkappa'$.
		For independent standard normal $\mathcal N',\mathcal N''$,
		$\varkappa\mathcal N\overset{\mathrm d}{=}\varkappa'\mathcal N'+\bar\varkappa\mathcal N''$ and
		$\E[\Psi'(\varkappa'\mathcal N')\mathcal N'']=0$, so \eqref{eq:Psiholder} gives
		\[
		|\E\Psi(\varkappa\mathcal N)-\E\Psi(\varkappa'\mathcal N)|
		=\bigl|\E\bigl[\Psi(\varkappa'\mathcal N'+\bar\varkappa\mathcal N'')-\Psi(\varkappa'\mathcal N')
		-\Psi'(\varkappa'\mathcal N')\,\bar\varkappa\mathcal N''\bigr]\bigr|
		\le c(1+|x|)^{5/2}|\varkappa^2-\varkappa'^2|^{3/4}.
		\]
		Coupling $\theta$ and $\theta'$ through the same $(\mathcal V,\mathcal N)$, the
		difference bound for $\mathfrak D$ in the proof of
		Lemma~\ref{lem:kernelcoupling} gives
		$|\E[\Psi_\theta(\varkappa\mathcal N)-\Psi_{\theta'}(\varkappa\mathcal N)]|
		\le c(1+|x|)^2\norm{\theta-\theta'}$. Since
		$|\varkappa^2-\varkappa'^2|\le c\norm{\theta-\theta'}$, the decomposition
		\[
		J_\theta(x)-J_{\theta'}(x)
		=\E[\Psi_\theta(\varkappa\mathcal N)-\Psi_{\theta'}(\varkappa\mathcal N)]
		+\E[\Psi_{\theta'}(\varkappa\mathcal N)-\Psi_{\theta'}(\varkappa'\mathcal N)],
		\]
		in which $\varkappa'$ is the value at $\theta'$, proves \eqref{eq:Jholder} and
		hence the uniform bound.
	\end{proof}
	The proof of Theorem~\ref{thm:singularcontinuity} uses $\Theta_0$ only through the constant $c_0$. On the set
	$\Theta(c_0)=\{\theta\in\Theta:s^2\ge c_0,\ \widetilde V\succeq c_0I_3,\
	J_\theta\ge c_0,\ \norm\theta\le c_0^{-1}\}$, which is compact by
	Lemma~\ref{lem:kernelcoupling}, the Lipschitz constants of
	$\theta\mapsto\varepsilon_*,b_*,m(x),\Omega$ are functions of $c_0$ as well.
	Hence $L_0$ can be taken as a nonincreasing function of $c_0$ alone, and
	\eqref{eq:lipschitzTV} holds uniformly on $\Theta(c_0)$, which contains any
	given compact $\Theta_0\subset\Theta$ once $c_0$ is small enough.

	\begin{corollary}[Rates from mean and covariance jets; regional counts]
		\label{cor:fieldparametercontinuity}
		Let $X_n$ and $X$ satisfy Assumption~\ref{ass:regularity}, with
		mean and covariance jets $\mathrm{jet}_n(t)$ and $\mathrm{jet}(t)$ as in
		Section~\ref{sec:regression}. Let $D\subset\mathcal M$ be relatively compact, and suppose
		that, in finitely many smooth local orthonormal frames covering $\bar D$,
		$d_n=\sup_{t\in\bar D}\norm{\mathrm{jet}_n(t)-\mathrm{jet}(t)}\to0$.
		Then there is $c$ such that, for all large $n$ and all $t\in\bar D$,
		\[
		\norm{h_{n,t}-h_t}_{L^1(\R)}+|\ell_n(t)-\ell(t)|
		+\norm{\rho_n(t,\cdot)-\rho(t,\cdot)}_{L^1(\R)}\le cd_n,\qquad
		\norm{h_{n,t}-h_t}_\infty\le cd_n^{3/4},
		\]
		and $|\E[N_n(D,B)]-\E[N(D,B)]|\le c\,|D|_{\mathsf g}\,d_n$ for every
		Borel set $B\subset\R$, where $|D|_{\mathsf g}$ is the Riemannian area of
		$D$. The pointwise bounds also hold at a single location $t$ when only
		$\mathrm{jet}_n(t)\to\mathrm{jet}(t)$.
	\end{corollary}
	\begin{proof}
		The jets of $X$ are continuous on the compact set $\bar D$ and satisfy
		\eqref{eq:non-degenerate}, so there is a compact set of jets containing
		$\mathrm{jet}(t)$ and, for large $n$, $\mathrm{jet}_n(t)$, $t\in\bar D$, on
		which $G$, $s^2$ and $\widetilde V$ are uniformly positive definite. There
		the regression parameters of Section~\ref{sec:regression} and
		$p_{\nabla X(t)}(0)$ are smooth, hence Lipschitz, functions of the jet,
		with values in a compact subset of $\Theta$.
		Theorem~\ref{thm:singularcontinuity} gives the bounds for $h$, and
		$\ell=p_{\nabla X(t)}(0)J(t)$ and $\rho(t,\cdot)=\ell(t)h_t$ give the others.
		Integrating the bound for $\rho$ over $D$ in \eqref{eq:KR} gives the
		count statement.
	\end{proof}
	Applied to the regression parameters $\theta(t_n)$ and $\theta(t)$ in one
	local frame, with $t_n\to t$, Theorem~\ref{thm:singularcontinuity} shows
	that $t\mapsto h_t$
	is continuous in $L^1(\R)$ and uniformly, including where $\Omega$ loses
	rank, and, wherever the jet is locally Lipschitz, Lipschitz in
	$L^1(\R)$ and H\"older of order $3/4$ uniformly.
	
	\paragraph{Nonsingular parameters.} If $\Theta_0$ contains no singular
	parameter, then $\varepsilon_*\ge\varepsilon_0>0$ on $\Theta_0$. In Step~5
	of the proof of Theorem~\ref{thm:singularcontinuity},
	\eqref{eq:Psiholder} and $\E[\Psi'(0)\mathcal N]=0$ show that
	$\varkappa^2\mapsto\E\Psi(\varkappa\mathcal N)$ has derivative of order
	$(1+|x|)^{5/2}\varkappa^{-1/2}$, and $\varkappa^2\ge c_0\varepsilon_0$ on
	$\Theta_0$. The uniform bound in \eqref{eq:lipschitzTV} therefore holds
	with exponent $1$ and a constant of order $\varepsilon_0^{-1/4}$, while the
	total variation bound is unchanged; the exponent $3/4$ is the price of
	uniformity as $\varepsilon_0\to0$. Likewise, if the jet law of $X$ is
	nonsingular at every point of $\bar D$ in
	Corollary~\ref{cor:fieldparametercontinuity}, then
	$\norm{h_{n,t}-h_t}_\infty\le cd_n$.

	\begin{remark}[Saddle and minimum height laws]
		Lemma~\ref{lem:kernelcoupling}, Proposition~\ref{prop:tiltedrep}
		except for its support statement, the bounds \eqref{eq:lipschitzTV}
		of Theorem~\ref{thm:singularcontinuity}, and
		Corollary~\ref{cor:fieldparametercontinuity} hold for saddles, with
		$\mathfrak D$ replaced by the saddle weight
		$\mathfrak D_{\mathrm{sdl}}=-p\,\ind_{\{p<0\}}$ of
		Section~\ref{subsec:saddles}. This weight has the same growth and
		difference bounds, and $p$ vanishes on the boundary $\{p=0\}$ of the
		saddle region, which is the union of two round cones, so Steps~1
		and~5 of the proof of Theorem~\ref{thm:singularcontinuity} apply to
		each cone; the saddle normalizer
		$\Kcal_{\mathrm{sdl}}(\bar m,\widetilde V)$ is positive because
		$\widetilde V\succ0$.
		Applying the peak results to $-X$ gives the conclusions for minima,
		including both optimality statements.
	\end{remark}

	\subsection{Endpoint behavior and optimality of the uniform rate}
	\label{subsec:endpoint}
	When $Q$ is semidefinite, the local height density vanishes at the
	endpoint $x_*$ of its support, at an order fixed by the type of the
	constraint. This order also determines the optimal uniform rate in
	Theorem~\ref{thm:singularcontinuity}.
	
	\begin{proposition}[Order of contact at the endpoint]
		\label{prop:contact}
		Fix $t$ with $\varepsilon_*=0$ and $Q$ semidefinite, use the
		coordinates of Theorem~\ref{thm:canonicalweights}, and write
		$\delta=\varsigma(x-x_*)$ with $\varsigma Q\succeq0$. There are
		$c_t,\delta_t>0$ such that, for $0<\delta<\delta_t$,
		\begin{equation}\label{eq:contact}
			h_t(x)=\begin{cases}
				c_t\,\delta^4\,\{1+O(\delta)\},&Q\text{ definite},\\
				c_t\,\delta^{3/2}\,\{1+O(\delta)\},&Q\text{ of rank one},
			\end{cases}
		\end{equation}
		while $h_t(x)=0$ for $\delta\le0$. In the definite case,
		$c_t=\pi c_*\varphi_{\Lambda_M}(\eta_M(x_*))/32$, where
		$c_*=\phin{x_*}{\bar\mu}{s^2}/\{\varpi_PJ(t)\}$ and $\varphi_{\Lambda_M}$ is the
		$\mathsf N_2(0,\Lambda_M)$ density; in the rank-one case,
		\[
		c_t=\frac43\,c_*\int_{-\infty}^0
		\phin{z}{\mu_{\tau_M}(x_*)}{\sigma_{\tau_M}^2}\,|z|^{3/2}\,
		\phin{0}{\beta_{B_0}(z;x_*)}{\lambda_{B_0}}\dd z.
		\]
	\end{proposition}
	\begin{proof}
		By \eqref{eq:peak-height} and Theorem~\ref{thm:canonicalweights},
		$h_t(x)=\phin{x}{\bar\mu}{s^2}\,\widehat J(x)/\{\varpi_PJ(t)\}$, and the
		Gaussian factor equals $\phin{x_*}{\bar\mu}{s^2}\{1+O(\delta)\}$. For
		$\delta\le0$, $h_t(x)=0$ by Theorem~\ref{thm:singularsupport}.
		
		In the definite case, $\xi\sim\mathsf N_2(\eta_M(x),\Lambda_M)$ has density
		$\varphi_{\Lambda_M}(\cdot-\eta_M(x))$, whose linear term at $0$
		integrates to zero over disks centered at $0$, so
		$\PP(\norm\xi^2\le y)=\pi y\varphi_{\Lambda_M}(\eta_M(x))+O(y^2)$, and \eqref{eq:Alayer} gives
		$\Acal(a;\eta_M(x),\Lambda_M)=\pi \varphi_{\Lambda_M}(\eta_M(x))a^2/2+O(a^3)$. Taking $a=\delta^2/4$
		in \eqref{eq:Q-definite}, and using that $\eta_M$ is affine in $x$, gives
		the first case.
		
		In the rank-one case, for $b>0$,
		\[
		\mathcal H(b;\beta,\lambda)=\int_{-\sqrt b}^{\sqrt b}(b-w^2)\phin{w}{\beta}{\lambda}\dd w
		=\frac43b^{3/2}\phin{0}{\beta}{\lambda}+O(b^{5/2})
		\]
		uniformly in $\beta$, because the linear term of
		$\phin{w}{\beta}{\lambda}-\phin{0}{\beta}{\lambda}$ integrates to zero and the
		second derivative of $\phin{\cdot}{\beta}{\lambda}$ is bounded by
		$c\lambda^{-3/2}$. In \eqref{eq:rankonesimplifiedweight},
		$b=-\delta(\delta+z)=\delta(|z|-\delta)\le\delta|z|$ for $z<-\delta$, so
		\[
		\widehat J(x)=\frac43\delta^{3/2}\int_{-\infty}^{-\delta}
		\phin{z}{\mu_{\tau_M}(x)}{\sigma_{\tau_M}^2}(|z|-\delta)^{3/2}
		\phin{0}{\beta_{B_0}(z;x)}{\lambda_{B_0}}\dd z+O(\delta^{5/2}).
		\]
		The means $\mu_{\tau_M}(x)$ and $\beta_{B_0}(z;x)$ are affine in $x$, with
		$x$-coefficients that do not depend on $z$, and
		$(|z|-\delta)^{3/2}=|z|^{3/2}+O(\delta|z|^{1/2})$. Hence the integral equals
		its value at $x=x_*$, taken over $(-\infty,0)$, up to $O(\delta)$; that
		value is positive, which gives the second case.
	\end{proof}
	
	For indefinite $Q$, $h_t$ is positive on $\R$ by
	Theorem~\ref{thm:singularsupport}, and for $\varepsilon_*>0$ it is
	real-analytic by Proposition~\ref{prop:tiltedrep}. In
	Section~\ref{subsubsec:singularfamily}, $J_{1/2}(x)\sim x^4/6$, and for
	$\gamma=1$ a numerical evaluation gives $h(x)/x^{3/2}\to0.1595$, the value
	of $c_t$ in \eqref{eq:contact}. The rank-one exponent $3/2$ determines the
	optimal uniform rate.
	
	\begin{proposition}[Optimality of the uniform exponent]
		\label{prop:sharprate}
		The exponent $3/4$ in \eqref{eq:lipschitzTV} cannot be increased. Let $Z$
		be the field of Section~\ref{subsubsec:singularfamily} with $\gamma=1$, and
		for $0<\epsilon\le1$ let $X_\epsilon=(Z+\epsilon\mathcal N_{\mathrm c})/\sqrt{1+\epsilon^2}$,
		where $\mathcal N_{\mathrm c}\sim\mathsf N(0,1)$ does not depend on $t$ and
		is independent of $Z$. Then $X_\epsilon$ satisfies
		Assumption~\ref{ass:regularity}, its regression parameters $\theta_\epsilon$
		at any $t$ ($\theta_0$ being those of $Z$) satisfy $\norm{\theta_\epsilon-\theta_0}\le c\epsilon^2$, and, for
		small $\epsilon$,
		\[
		\norm{h_{\theta_\epsilon}-h_{\theta_0}}_\infty\ge h_{\theta_\epsilon}(0)
		\ge c'\epsilon^{3/2}\ge c''\norm{\theta_\epsilon-\theta_0}^{3/4},
		\]
		with constants $c,c',c''>0$. 
	\end{proposition}
	
	\begin{proof}
		Since $\mathcal N_{\mathrm c}$ does not depend on $t$, $X_\epsilon$ has
		unit variance, its gradient and Hessian are those of $Z$ divided by
		$(1+\epsilon^2)^{1/2}$, and $\Cov(\nabla X_\epsilon(t),X_\epsilon(t))=0$,
		so Assumption~\ref{ass:regularity} holds. Its regression parameters are
		$\bar\mu=0$, $s^2=1$, $\bar m=0$, $\widetilde\iota=\iota/(1+\epsilon^2)$
		and $\widetilde V=V/(1+\epsilon^2)$, with $\iota\ne0$ and $V$ from
		Section~\ref{subsubsec:singularfamily}. Hence
		$0<\norm{\theta_\epsilon-\theta_0}\le c\epsilon^2$ with
		$c=\norm\iota+\norm V_{\mathrm F}$, $b_*=V^{-1}\iota$ does not depend on
		$\epsilon$, and $\varepsilon_*=\sigma_\epsilon^2$ with
		$\sigma_\epsilon=\epsilon/\sqrt{1+\epsilon^2}$. As $\mathfrak D$ is
		homogeneous of degree two, the tilted law $\nu_{\theta_\epsilon}$ is the
		image of $\nu_{\theta_0}$ under $w\mapsto w/\sqrt{1+\epsilon^2}$, and
		Proposition~\ref{prop:tiltedrep} shows that $h_{\theta_\epsilon}$ is the
		density of $(Y+\epsilon\mathcal N)/\sqrt{1+\epsilon^2}$, with
		$Y\sim h_{\theta_0}$ and $\mathcal N\sim\mathsf N(0,1)$ independent. Here
		$Q=\diag(1,0)$ has rank one and $x_*=0$, so
		Proposition~\ref{prop:contact} gives $h_{\theta_0}(0)=0$, which yields the
		first inequality, and $h_{\theta_0}(y)\ge c_ty^{3/2}/2$ for
		$0<y<\delta_t$, after decreasing $\delta_t$ if necessary. For
		$0<y<\sigma_\epsilon$, the normal density with mean
		$y/\sqrt{1+\epsilon^2}$ and variance $\sigma_\epsilon^2$ is at least
		$e^{-1/2}/(\sigma_\epsilon\sqrt{2\pi})$ at $0$, so restricting the
		convolution integral to $(0,\sigma_\epsilon)$ gives
		$h_{\theta_\epsilon}(0)\ge c_te^{-1/2}\sigma_\epsilon^{3/2}/(5\sqrt{2\pi})
		\ge c'\epsilon^{3/2}$ for small $\epsilon$, since
		$\sigma_\epsilon\ge\epsilon/\sqrt2$; the bound
		$\norm{\theta_\epsilon-\theta_0}\le c\epsilon^2$ then gives the last
		inequality. Finally, all $\theta_\epsilon$, $0\le\epsilon\le1$, lie in a
		compact subset of $\Theta$, on which a bound of order
		$\norm{\theta-\theta'}^{\beta}$ with $\beta>3/4$ in
		\eqref{eq:lipschitzTV} would give
		$c''\le L_0\norm{\theta_\epsilon-\theta_0}^{\beta-3/4}\to0$, so the exponent $3/4$ in
		\eqref{eq:lipschitzTV} cannot be increased.
	\end{proof}
	
	The lower bound comes from the endpoint $x_*$: at a fixed distance from
	it, the smoothing in the
	proof changes $h_{\theta_0}$ by
	$O(\epsilon^2)=O(\norm{\theta_\epsilon-\theta_0})$, whereas at $x_*$ the
	change is of order
	$\epsilon^{3/2}=\norm{\theta_\epsilon-\theta_0}^{3/4}$, because
	$h_{\theta_0}(y)\asymp y^{3/2}$ as $y\downarrow x_*$. For $\gamma=1/2$ in the
	Fourier family of Section~\ref{subsubsec:singularfamily}, $h_t$ vanishes
	to order~$4$ at $x_*$, and the same construction changes $h_{\theta_0}$
	by $O(\epsilon^2)$ uniformly, a linear uniform rate.

	\begin{remark}[Optimality of the total variation rate]
		\label{rem:optimalityscope}
		The linear rate in total variation in \eqref{eq:lipschitzTV} is optimal as well. Adding a constant $c$ to $X$ changes $\theta$ only through $\bar\mu\mapsto\bar\mu+c$ and shifts the height density, $h_{\theta'}(x)=h_\theta(x-c)$. For nonsingular $\theta$, $h_\theta$ is smooth with integrable derivative, so $\norm{h_\theta(\cdot-c)-h_\theta}_{L^1(\R)}/|c|\to\norm{h_\theta'}_{L^1(\R)}>0$ as $c\to0$.
		\end{remark}

	\section{The Euclidean plane}
	\label{sec:stationary}
	Here $\mathcal M=T\subset\R^2$ carries the Euclidean metric, the
	covariant Hessian is the ordinary Hessian, and $\dd A_{\mathsf g}(t)=\dd t$.
	Section~\ref{subsec:stationarytrend} treats stationary fields, including
	isotropic ones, and Section~\ref{subsec:quadratic} treats a
	nonstationary field.
	
	\subsection{Stationary fields}
	\label{subsec:stationarytrend}
	Let $X=\mu+Z$, where $Z$ is centered and stationary with covariance
	$C_0(t-t')$. Evenness of $C_0$ and covariance differentiation give $r=A=0$, $\iota=-\vech(G)$ and $V_{12}=V_{33}$.
	Thus the gradient is independent of the height and Hessian, and
	$v$, $G$, and $V$ are constant. At a fixed $t$, put
	$m_W=\vech(\nabla^2\mu(t))$. Gaussian regression and
	Theorem~\ref{thm:main} reduce to $J(t,x)=\Kcal(m(x),\Omega)$ and $J(t)=\Kcal(m_W,V)$ with $m(x)=m_W+[x-\mu(t)]\iota /v$ and 
	$\Omega=V-\iota\iota^\top/v$.
	The remaining Gaussian factors are
	\[
	p_{X(t)\mid\nabla X(t)=0}(x)=\phin{x}{\mu(t)}v,\qquad
	p_{\nabla X(t)}(0)=
	\frac{\exp\{-\tfrac12(\nabla\mu(t))^\top G^{-1}\nabla\mu(t)\}}
	{2\pi\sqrt{\det G}}.
	\]
	Substitution in \eqref{eq:KR}--\eqref{eq:peak-height} gives the peak
	intensities, height density, and expected counts.
	
	\paragraph{Singular stationary constraints.}
	Here the singularity criterion is
	$\varepsilon_*=v-\iota^\top V^{-1}\iota=0$.
	In this case, Theorem~\ref{thm:singularsupport} gives the constraint
	operator $Q$ and endpoint $x_*(t)=\mu(t)+\tr\bigl(Q\nabla^2\mu(t)\bigr)$.
	The identity $\tr(Q\nabla^2Z)=-Z$, together with
	$\Cov(Z,\nabla^2Z)=-G$, yields $\tr(QG)=v>0$.
	Since $G\succ0$, the constraint cannot be negative semidefinite.
	Equivalently, by Corollary~\ref{cor:spectralconic}, the spectral
	measure lives on the conic $\{k^\top Qk=1\}$, which is empty when
	$Q\preceq0$.
	The stationary possibilities are therefore
	\[
	\{x:h_t(x)>0\}=
	\begin{cases}
		(x_*(t),\infty),&Q\succeq0,\\
		\R,&Q\text{ is indefinite}.
	\end{cases}
	\]
	The original trace is deterministic given height precisely when
	$Q=(v/\tr G)I_2$; its value is then $\tau_x=\tr(\nabla^2\mu(t))
	-(\tr G)\bigl(x-\mu(t)\bigr)/v$.
	Its conditional peak weight follows from the deterministic trace branch of
	\eqref{eq:substitution}. Theorem~\ref{thm:canonicalweights} evaluates
	the other constraints, and saddle heights have support $\R$ by
	Theorem~\ref{thm:singularsupport}.
	
	\paragraph{Total counts for centered fields.}
	\label{subsec:stationarycounts}
	For centered stationary planar fields, the total density of maxima was
	obtained by Longuet-Higgins \citep[Section~2.4]{LonguetHiggins1957} in
	terms of complete elliptic integrals. It is the stationary case of
	Corollary~\ref{cor:centerednormalizer}.
	\begin{proposition}[Longuet-Higgins]
		\label{prop:stationarycounts}
		Let $X$ be a centered stationary Gaussian field satisfying
		Assumption~\ref{ass:regularity}, and let $\vartheta>0,-\lambda_1,-\lambda_2$
		be the eigenvalues of $V^{1/2}\Upsilon  V^{1/2}$ in
		\eqref{eq:determinantspectrum}, with $\lambda_1\ge\lambda_2>0$. Then
		$\vartheta=\lambda_1+\lambda_2$ and, with the complete elliptic
		integrals $\mathrm E(k)=\int_0^{\pi/2}(1-k^2\sin^2\psi)^{1/2}\dd\psi$ and $\mathrm K(k)=\int_0^{\pi/2}(1-k^2\sin^2\psi)^{-1/2}\dd\psi$,
		\begin{equation}\label{eq:stationaryellipticnormalizer}
			J=\frac{\sqrt{\lambda_1(\lambda_1+2\lambda_2)}}{\pi}
			\left[\mathrm E(k_\lambda)
			-\frac{\lambda_2}{\lambda_1+2\lambda_2}\mathrm K(k_\lambda)\right],
			\qquad
			k_\lambda^2=\frac{\lambda_1^2-\lambda_2^2}{\lambda_1(\lambda_1+2\lambda_2)}.
		\end{equation}
		Moreover, $\ell=J/(2\pi\sqrt{\det G})$, $\ell_{\mathrm{min}}=\ell$ and
		$\ell_{\mathrm{sdl}}=2\ell$.
	\end{proposition}
	\begin{proof}
		Since $V_{12}=V_{33}$, $\vartheta-\lambda_1-\lambda_2=\tr(\Upsilon  V)
		=V_{12}-V_{33}=0$. Polar coordinates give
		$\Acal(y;0,\diag(\lambda_1,\lambda_2))=\frac2\pi\int_0^{\pi/2}
		\Acal(y;0,\lambda_\psi I_2)\dd\psi$ with
		$\lambda_\psi=\lambda_1\cos^2\psi+\lambda_2\sin^2\psi$. Inserting this
		into \eqref{eq:centerednormalizer} and using
		\eqref{eq:centerednormalizercircular} with $\vartheta=\lambda_1+\lambda_2$
		gives 
		\eqref{eq:stationaryellipticnormalizer}. Finally, the symmetry
		$H\overset{\mathrm d}=-H$ with $\E[\det H]=V_{12}-V_{33}=0$ gives
		$\ell_{\mathrm{min}}=\ell$ and $\ell_{\mathrm{sdl}}=2\ell$.
	\end{proof}
	Since $\vartheta,-\lambda_2,-\lambda_1$ are the roots of Longuet-Higgins'
	cubic \citep[(2.4.27)]{LonguetHiggins1957}, 
	$\vartheta=\lambda_1+\lambda_2$ is his relation $l_1+l_2+l_3=0$, and
	\eqref{eq:stationaryellipticnormalizer} reproduces his peak density.
	His height distributions are restricted to narrow spectra, whereas the
	conditional peak weights $J(t,x)=\Kcal(m(x),\Omega)$ above apply to any stationary field,
	with or without a trend; Proposition~\ref{prop:stationarycounts} itself
	concerns centered fields.
	
	\paragraph{Isotropic fields.}
	\label{subsec:isotropic}
	Let $X$ be centered, stationary and isotropic with unit variance,
	$C(t,t')=\varrho(\norm{t-t'}^2)$, $a=-2\varrho'(0)$, $b=4\varrho''(0)$ and
	$\kappa=a/\sqrt b\in(0,\sqrt2]$. Then $G=aI_2$ and, given height $x$ and
	zero gradient, $\tau$ and $\xi$ are independent with $\tau\sim\mathsf N\bigl(-2ax,4b(2-\kappa^2)\bigr)$ and $\xi\sim\mathsf N_2(0,bI_2)$.
	Theorem~\ref{thm:main} with the circular kernel
	\eqref{eq:centeredcircularkernel} gives $J(x)$ as a single Gaussian
	integral, whose evaluation recovers the known isotropic peak height
	density \citep[Examples~3.8 and~3.15]{ChengSchwartzman2018}, a closed
	form in $\phi$ and $\Phi$ for $0<\kappa<\sqrt2$.
	Proposition~\ref{prop:stationarycounts} with $\lambda_1=\lambda_2=b$ gives
	$J=b/\sqrt3$ and $\ell=b/(2\pi\sqrt3\,a)$. The endpoint $\kappa=\sqrt2$
	(random plane waves) is the definite singular case of
	Theorem~\ref{thm:singularsupport} with $Q=I_2/(2a)$: the trace equals
	$-2ax$, and
	\begin{equation}\label{eq:isotropicboundaryheight}
		h(x)=2\sqrt3\,\phi(x)\bigl[x^2-1+e^{-x^2}\bigr]\ind_{\{x>0\}},
	\end{equation}
	the uniform and total variation limit of the nonsingular isotropic laws as
	$\kappa\uparrow\sqrt2$ by Theorem~\ref{thm:singularcontinuity}. The saddle
	height density is $\phin{x}{0}{1-\kappa^2/3}$ for $0<\kappa\le\sqrt2$
	\citep{ChengSchwartzman2018}.
	
	\subsection{A nonstationary quadratic field}
	\label{subsec:varyingvariance}
	\label{subsec:quadratic}
	Let $\zeta_0,\ldots,\zeta_4$ be independent standard normal variables,
	and consider
	\[
	X(t)=\zeta_0\left(1+\frac{\norm t^2}{2}\right)
	+\zeta_1t_1+\zeta_2t_2
	+\frac{\zeta_3}{2}(t_1^2-t_2^2)+\zeta_4t_1t_2,
	\qquad t\in\R^2.
	\]
	Write $q=\norm t^2$. The field is centered and rotationally invariant,
	but its variance varies with location. In the original Hessian entry
	order, direct differentiation gives
	\[
	\begin{gathered}
		v=1+2q+\frac{q^2}{2},\qquad r=(q+2)t,\qquad
		G=(1+q)I_2+tt^\top,\\
		\iota=(1+t_1^2,1+t_2^2,t_1t_2)^\top,\qquad
		A=\begin{pmatrix}2t_1&0\\0&2t_2\\t_2&t_1\end{pmatrix},
		\qquad V=\diag(2,2,1).
	\end{gathered}
	\]
	Thus $r\ne0$ and $A\ne0$ away from the origin. Conditioning on zero
	gradient gives
	\begin{equation}\label{eq:variablevarianceheight}
		X(t)\mid\nabla X(t)=0\sim\mathsf N(0,s_q^2),\qquad
		s_q^2=\frac{q^2+2}{2(1+2q)}.
	\end{equation}
	
	The assumptions hold on all of $\R^2$: the paths are quadratic, the
	Hessian increments vanish, $\det\Cov\bigl((\nabla X(t)^\top,W^\top)^\top\bigr)=4$ and $\det\Cov\bigl((X(t),\nabla X(t)^\top)^\top\bigr)=(1+q)(q^2+2)/2>0$.
	
	For $t\ne0$, use the radial frame
	$e_1=t/\sqrt q$, $e_2=(-t_2,t_1)^\top/\sqrt q$, and $O=(e_1,e_2)$;
	at the origin take $O=I_2$. In this frame, Gaussian regression gives
	\[
	O^\top\nabla^2X(t)O\mid\{X(t)=x,\nabla X(t)=0\}
	\overset{\mathrm d}=
	\begin{pmatrix}H_{\mathrm{rr}}&H_{\mathrm{rt}}\\H_{\mathrm{rt}}&(q-1)H_{\mathrm{rr}}+2x\end{pmatrix},
	\]
	where
	\[
	H_{\mathrm{rr}}\sim\mathsf N\!\left(\frac{2(1-q)x}{q^2+2},\frac2{q^2+2}\right),
	\qquad H_{\mathrm{rt}}\sim\mathsf N\!\left(0,\frac1{1+q}\right)
	\]
	are independent. Thus the height-conditioned Hessian covariance has
	rank two at every location.
	
	The exact identity $X(t)=t^\top\nabla X(t)
	+\tr\!\big[(I_2-tt^\top)\nabla^2X(t)\big]/2$
	identifies the constraint matrix as $Q(t)=(tt^\top-I_2)/2$, with
	$x_*=0$. Its radial and tangential eigenvalues are $(q-1)/2$ and
	$-1/2$. By Theorem~\ref{thm:singularsupport}, the positivity set of
	$h_t$ is $(-\infty,0)$ for $q\le1$ and $\R$ for $q>1$. Thus this one
	field realizes negative definite, rank-one negative semidefinite, and
	indefinite constraints as the location passes through the unit circle.
	In the terms of Theorem~\ref{thm:singularpde}, $Z$ solves on $\R^2$
	the equation $Z-t^\top\nabla Z-\tfrac12\tr\{(I_2-tt^\top)\nabla^2Z\}=0$,
	with $Q(t)=(tt^\top-I_2)/2$ and $\beta(t)=-t$. It is elliptic in the open
	unit disk, parabolic on the unit circle and hyperbolic outside it, a
	mixed-type equation like the Tricomi equation, and the peak height
	support changes accordingly.
	
	Sylvester's criterion gives negative definiteness exactly when
	$H_{\mathrm{rr}}<0$ and $(q-1)H_{\mathrm{rr}}^2+2xH_{\mathrm{rr}}-H_{\mathrm{rt}}^2>0$. Conditioning on $H_{\mathrm{rr}}=z$ therefore
	gives one scalar-kernel formula for every $q\ge0$:
	\begin{equation*}
		J(t,x)=\int_{-\infty}^0
		\phin{z}{2(1-q)x/(q^2+2)}{2/(q^2+2)}
		\mathcal H\!\left((q-1)z^2+2xz;0,\frac1{1+q}\right)\dd z.
	\end{equation*}
	The convention $\mathcal H(b;0,\lambda)=0$ for $b\le0$ enforces the
	determinant constraint, including at the origin.
	Given zero gradient, the radial Hessian entries are independent and
	centered with variances $2/(1+2q)$, $2$ and $1/(1+q)$, so
	Corollary~\ref{cor:centerednormalizer} applies with
	$\vartheta=\lambda_1=(1+2q)^{-1/2}$ and $\lambda_2=(1+q)^{-1}$. With
	\eqref{eq:variablevarianceheight} and
	$p_{\nabla X(t)}(0)=[2\pi\sqrt{(1+q)(1+2q)}]^{-1}$ this gives all peak
	quantities, which depend on $t$ only through $q$.
	Figure~\ref{fig:surfaces}(a) shows the local height densities on both
	sides of the unit circle.

	\section{Gaussian fields on curved surfaces}
	\label{sec:sphere}
	In Sections~\ref{subsec:sphericalisotropic}--\ref{subsec:sphericalchecks},
	$\mathcal M=\mathbb S^2$, the unit sphere with the round metric, and
	$\Delta_{\mathbb S^2}$ denotes the Laplace--Beltrami operator. The jet can
	be computed intrinsically, for instance in spherical coordinates. For
	fields defined through ambient coordinates, such as the quadrupole below,
	it is convenient to use a $C^2$ extension $\widetilde X$ of $X$ to a
	neighborhood of $\mathbb S^2$ and the matrix $\mathsf E(t)\in\R^{3\times2}$
	of ambient components of an orthonormal tangent frame, so that
	$\mathsf E^\top\mathsf E=I_2$ and $\mathsf E\mathsf E^\top=I_3-tt^\top$. Then
	\begin{equation}\label{eq:sphericalambienthessian}
		\nabla X(t)=\mathsf E(t)^\top\nabla_{\R^3}\widetilde X(t),\qquad
		H(t)=\mathsf E(t)^\top\nabla_{\R^3}^2\widetilde X(t)\mathsf E(t)
		-\bigl(t^\top\nabla_{\R^3}\widetilde X(t)\bigr)I_2,
	\end{equation}
	which depends only on $X$ and the round metric, not on the extension.
	Section~\ref{subsec:ellipsoid} treats an ellipsoid, whose curvature is not constant.
	
	\subsection{Isotropic fields on $\mathbb S^2$}
	\label{subsec:sphericalisotropic}
	Let $X$ be centered, isotropic and of unit variance,
	$C(t,t')=\varrho_{\mathbb S^2}(t^\top t')$, $a=\varrho_{\mathbb S^2}'(1)$ and
	$b=\varrho_{\mathbb S^2}''(1)$. Then $r=A=0$, $G=aI_2$, $\iota=-a(1,1,0)^\top$,
	and $V$ has rows $(3b+a,b+a,0)$, $(b+a,3b+a,0)$ and $(0,0,b)$
	\citep[Lemma~4.1]{ChengSchwartzman2018}, so that, given height $x$ and
	zero gradient, $\tau\sim\mathsf N(-2ax,4(a+2b-a^2))$ and
	$\xi\sim\mathsf N_2(0,bI_2)$ independently. Theorem~\ref{thm:main} with
	\eqref{eq:centeredcircularkernel} and \eqref{eq:centerednormalizercircular}
	recovers the height density and
	$\E[N(\mathbb S^2)]=1+2b^{3/2}/(a\sqrt{a+3b})$ of
	\citep[Example~4.6]{ChengSchwartzman2018}. The case $a+2b=a^2$ is that of
	a single harmonic degree $j\ge2$, with $a=j(j+1)/2$ and
	$\Delta_{\mathbb S^2}X=-2aX$ \citep[Proposition~4.8]{ChengSchwartzman2018}: the
	definite case of \eqref{eq:generaleigenfunctionconstraint}, with peak
	height support $[0,\infty)$. Adding a random constant,
	$Y_{j,\epsilon}=(Y_j+\epsilon\mathcal N_{\mathrm c})/\sqrt{1+\epsilon^2}$, where
	$\mathcal N_{\mathrm c}\sim\mathsf N(0,1)$ does not depend on $t$ and is
	independent of $Y_j$,
	leaves the critical points of a degree-$j$ field $Y_j$ unchanged and turns
	its peak height density $h_j$ into the Gaussian mixture
	$\int_0^\infty h_j(y)\phin{x}{y/\sqrt{1+\epsilon^2}}{\epsilon^2/(1+\epsilon^2)}\dd y$,
	a nonsingular law that converges to $h_j$ as in
	Theorem~\ref{thm:singularcontinuity}; Proposition~\ref{prop:sharprate}
	uses the same construction on the plane. For $j=2$, $a=b=3$ and
	\begin{equation}\label{eq:sphericalisotropicquadrupole}
		h(x)=\phi(x)\bigl(3x^2-2+2e^{-3x^2/2}\bigr)\ind_{\{x>0\}},
		\qquad \ell=(2\pi)^{-1},
	\end{equation}
	so that the expected number of maxima is two; the saddle height law is
	$\mathsf N(0,1/4)$.
	
	\subsection{A nonisotropic quadrupole on $\mathbb S^2$}
	\label{subsec:sphericalquadrupole}
	Let $\zeta=(\zeta_1,\ldots,\zeta_5)^\top\sim\mathsf N_5(0,\Sigma_\zeta)$ with
	$\Sigma_\zeta\succ0$, and let $\Xi$ be the traceless symmetric matrix with
	rows $(\zeta_1+\zeta_2,\zeta_3,\zeta_4)$, $(\zeta_3,\zeta_1-\zeta_2,\zeta_5)$
	and $(\zeta_4,\zeta_5,-2\zeta_1)$. The field and its covariance are
	\begin{equation}\label{eq:sphericalquadrupolemodel}
		X(t)=t^\top\Xi t,\qquad C(t,t')=f(t)^\top\Sigma_\zeta f(t'),
	\end{equation}
	with $f(t)=(t_1^2+t_2^2-2t_3^2,\,t_1^2-t_2^2,\,2t_1t_2,\,2t_1t_3,\,2t_2t_3)^\top$;
	for general $\Sigma_\zeta$ the variance is nonconstant and the covariance
	nonisotropic. By \eqref{eq:sphericalambienthessian} with
	$\widetilde X(t)=t^\top\Xi t$,
	\begin{equation}\label{eq:sphericalquadrupolejet}
		\nabla X(t)=2\mathsf E(t)^\top\Xi t,\qquad
		H(t)=2\mathsf E(t)^\top\Xi\mathsf E(t)-2X(t)I_2,\qquad
		\tau(t)=-6X(t),
	\end{equation}
	the last identity because $\tr\Xi=0$; it is the degree-two eigenfunction
	equation and holds for every realization.
	
	\begin{proposition}[Nonisotropic quadrupole]
		\label{prop:sphericalquadrupole}
		The field \eqref{eq:sphericalquadrupolemodel} satisfies
		Assumption~\ref{ass:regularity}. Given $X(t)=x$ and $\nabla X(t)=0$, the
		trace equals $-6x$ and
		$\xi(t)\sim\mathsf N_2(m_\xi(t,x),S_\xi(t))$, with $m_\xi$ linear in $x$
		and $S_\xi(t)\succ0$; the conditional Hessian covariance has rank two.
		Hence
		\begin{equation}\label{eq:sphericalquadrupoleweights}
			J(t,x)=\ind_{\{x>0\}}\Acal\bigl(9x^2;m_\xi(t,x),S_\xi(t)\bigr),\qquad
			\rho(t,x)=\frac{\phin{x}{0}{s^2(t)}}{2\pi\sqrt{\det G(t)}}\,J(t,x),
		\end{equation}
		and the local peak height density is positive exactly on $(0,\infty)$.
	\end{proposition}
	\begin{proof}
		At $e_3=(0,0,1)^\top$, with the first two coordinate vectors as frame,
		$X(e_3)=-2\zeta_1$, $\nabla X(e_3)=2(\zeta_4,\zeta_5)^\top$ and
		$\vech H(e_3)=(6\zeta_1+2\zeta_2,6\zeta_1-2\zeta_2,2\zeta_3)^\top$, so that
		$\xi(e_3)=2(\zeta_2,\zeta_3)^\top$. Hence $\zeta\mapsto(\nabla X(e_3),\vech H(e_3))$
		and $\zeta\mapsto(X(e_3),\nabla X(e_3),\xi(e_3))$ are invertible, which gives
		\eqref{eq:non-degenerate}, $s^2(e_3)>0$ and $S_\xi(e_3)\succ0$. A rotation
		moving $t$ to $e_3$ acts on $\Xi$ by conjugation, an invertible linear map
		on traceless symmetric matrices, so the same holds at every $t$. The
		paths are restrictions of quadratic polynomials, which gives the
		regularity condition. Since the trace equals $-6x$, the deterministic
		trace branch of Theorem~\ref{thm:main} gives
		\eqref{eq:sphericalquadrupoleweights}; for $x\le0$ no negative definite
		Hessian is possible.
	\end{proof}
	All parameters follow from the jet covariance $\Sigma(t)$, obtained by
	applying the linear jet maps to $\Sigma_\zeta$; correlations between the
	traceless Hessian, height and gradient generally make $m_\xi(t,x)\ne0$.
	For independent coefficients with $\Var(\zeta_j)=\sigma_j^2$, the quantities
	at $e_3$ are independent, with $s^2(e_3)=4\sigma_1^2$,
	$G(e_3)=4\diag(\sigma_4^2,\sigma_5^2)$, $m_\xi(e_3,x)=0$ and
	$S_\xi(e_3)=4\diag(\sigma_2^2,\sigma_3^2)$; if $\sigma_2^2\ne\sigma_3^2$, the
	closed form \eqref{eq:centeredordinarykernel} evaluates $J(e_3,x)$ through
	three ordinary noncentral chi square CDFs. Unequal traceless variances thus
	affect only the kernel, whereas the support comes from the eigenfunction
	equation. The choice $\sigma_1^2=1/4$, $\sigma_2^2=\cdots=\sigma_5^2=3/4$ is
	the rotation-invariant law of $\Xi$; then $C(t,t')=P_2(t^\top t')$, where
	$P_2(u)=(3u^2-1)/2$ is the Legendre polynomial of degree two, and the
	peak law is \eqref{eq:sphericalisotropicquadrupole}.
	
	\paragraph{Trends.}
	\label{subsec:sphericaltrend}
	For $X=\mu+Z$ with $Z$ the centered quadrupole and $\mu\in C^2(\mathbb S^2)$, the
	identity $\Delta_{\mathbb S^2}Z=-6Z$ gives $\tau(t)=-6\bigl(X(t)-x_*(t)\bigr)$ and $x_*(t)=\mu(t)+\Delta_{\mathbb S^2}\mu(t)/6$,
	so that $J(t,x)=\ind_{\{x>x_*(t)\}}\Acal(9[x-x_*(t)]^2;m_\xi(t,x),S_\xi(t))$,
	with $S_\xi$ unchanged and $m_\xi$ now affine in $x$. This is the definite
	case of Theorem~\ref{thm:singularsupport}, with height support
	$[x_*(t),\infty)$. The intensity uses $\phin{x}{\bar\mu(t)}{s^2(t)}$ and the
	general gradient density; $\bar\mu(t)$ need not equal $\mu(t)$ or
	$x_*(t)$. A constant trend shifts the endpoint by that constant, and a
	traceless quadratic trend, for which $(\Delta_{\mathbb S^2}+6)\mu=0$, leaves it at
	zero.
	
	\subsection{An exact eigenvalue check}
	\label{subsec:sphericalchecks}
\begin{figure}[htbp]
		\centering
		\includegraphics[width=\textwidth]{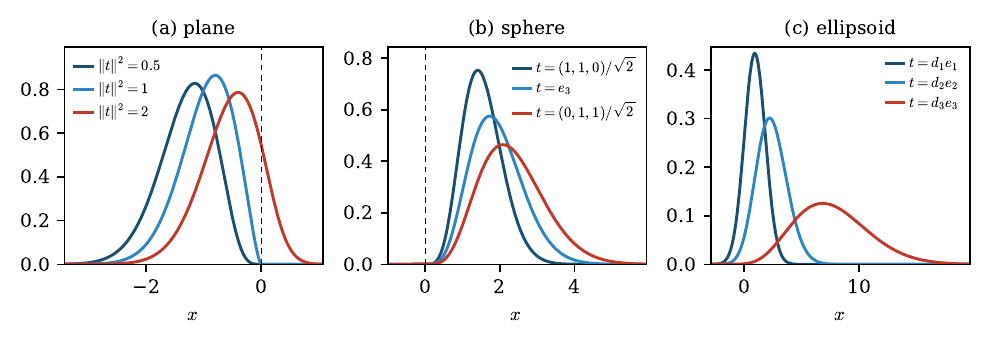}
		\caption{Local peak height densities $h_t$; dashed vertical lines mark
			the endpoint $0$ of the support where it is a half line. (a)
			Quadratic field of
			Section~\ref{subsec:quadratic}, inside, on and outside the unit
			circle; the support is $(-\infty,0]$ for $\norm t\le1$ and $\R$ for
			$\norm t>1$. (b) Nonisotropic quadrupole of
			Section~\ref{subsec:sphericalquadrupole} with independent
			coefficients of variances $(0.30,1.00,0.40,0.70,1.50)$, at
			$t=(1,1,0)^\top/\sqrt2$, $t=e_3$ and $t=(0,1,1)^\top/\sqrt2$; the first and last are close
			to the smallest and largest mean peak heights over $\mathbb S^2$,
			$1.540$ and $2.345$. By Proposition~\ref{prop:sphericalquadrupole},
			every local density vanishes on $(-\infty,0]$.
			(c) Quadratic form on the ellipsoid with semi-axes $(1,1.5,2)$,
			case (i) of Proposition~\ref{prop:ellipsoid} with independent
			entries of variances $(1.0,0.6,1.4,0.8,0.5,1.2)$ for
			$(\Xi_{11},\Xi_{22},\Xi_{33},\Xi_{12},\Xi_{13},\Xi_{23})$, at the
			three vertices $d_ie_i$.}
		\label{fig:surfaces}
	\end{figure}

	The pure quadrupole admits a check of the height-resolved local formulas
	that does not use the Kac--Rice formula. The eigenvalues
	$\alpha_{(1)}<\alpha_{(2)}<\alpha_{(3)}$ of $\Xi$ are distinct almost surely,
	and $t$ is critical exactly when $\Xi t$ is parallel to $t$. At a unit
	eigenvector with eigenvalue $\alpha_{(i)}$, the height is $\alpha_{(i)}$
	and the intrinsic Hessian eigenvalues are $2(\alpha_{(j)}-\alpha_{(i)})$,
	$j\ne i$. Hence the field is almost surely Morse, with two maxima at height
	$\alpha_{(3)}$, two saddles at height $\alpha_{(2)}$ and two minima at height
	$\alpha_{(1)}$, and for every $\Sigma_\zeta\succ0$,
	\begin{equation}\label{eq:sphericalquadrupolecountcheck}
		\begin{split}
			\frac12\int_{\mathbb S^2}\rho(t,x)\,\dd A_{\mathsf g}(t)&=p_{(3)}(x),\qquad
			\frac12\int_{\mathbb S^2}\rho_{\mathrm{sdl}}(t,x)\,\dd A_{\mathsf g}(t)=p_{(2)}(x),\\
			\frac12\int_{\mathbb S^2}\rho_{\mathrm{min}}(t,x)\,\dd A_{\mathsf g}(t)&=p_{(1)}(x),
		\end{split}
	\end{equation}
	where $p_{(i)}$ is the density of $\alpha_{(i)}$. The left sides use the full
	local machinery (jet covariance, two-stage regression and norm kernel); the
	right sides are densities of the ordered eigenvalues of a $3\times3$
	Gaussian matrix. Figure~\ref{fig:surfaces}(b) shows three local peak
	height densities of the quadrupole with independent coefficients of
	variances $(0.30,1.00,0.40,0.70,1.50)$: all are supported on
	$[0,\infty)$, but nonisotropy moves them markedly.
	
	In the isotropic case, \eqref{eq:sphericalquadrupolecountcheck} identifies
	\eqref{eq:sphericalisotropicquadrupole} as the density of the largest
	eigenvalue of a rotation-invariant traceless $3\times3$ Gaussian matrix, and
	$\mathsf N(0,1/4)$ as the law of its middle eigenvalue. The latter can be
	confirmed directly: the joint eigenvalue density is proportional to
	$\prod_{i<j}(\alpha_{(j)}-\alpha_{(i)})\,e^{-\tr\Xi^2/3}$ on
	$\{\alpha_{(1)}+\alpha_{(2)}+\alpha_{(3)}=0\}$; with $y=\alpha_{(2)}$ and
	$d=\alpha_{(3)}-\alpha_{(1)}>3|y|$, the product equals $d(d^2-9y^2)/4$ and
	$\tr\Xi^2=(3y^2+d^2)/2$, and integrating over $d$ by the substitution
	$d^2-9y^2=v$ leaves a density proportional to $e^{-2y^2}$. These six counts
	are preserved by adding a constant or a traceless quadratic trend.
	They satisfy the Morse relation
	$N(\mathbb S^2)-N_{\mathrm{sdl}}(\mathbb S^2)+N_{\mathrm{min}}(\mathbb S^2)
	=\chi(\mathbb S^2)=2$, which holds for every Morse function on
	$\mathbb S^2$; the symmetry $X\mapsto-X$ is needed only to replace
	$\E[N_{\mathrm{min}}(\mathbb S^2)]$ by $\E[N(\mathbb S^2)]$ in the expected
	form $\E[N_{\mathrm{sdl}}(\mathbb S^2)]=2\E[N(\mathbb S^2)]-\chi(\mathbb S^2)$
	\citep{AdlerTaylor2007}.

	\subsection{A quadratic form on an ellipsoid}\label{subsec:ellipsoid}
	
	For a surface $\mathcal M\subset\R^3$ with the induced metric, unit normal
	$n(t)$ and second fundamental form
	$\mathrm{I\!I}_t(u,v)=\langle\overline\nabla_un,v\rangle$, where
	$\overline\nabla$ is the Euclidean derivative, the Gauss formula
	$\nabla_uv=\overline\nabla_uv+\mathrm{I\!I}(u,v)\,n$ for tangent fields
	extends \eqref{eq:sphericalambienthessian} to
	\begin{equation}\label{eq:embeddedhessian}
		\nabla X(t)=\mathsf E(t)^\top\nabla_{\R^3}\widetilde X(t),\qquad
		H(t)=\mathsf E(t)^\top\nabla_{\R^3}^2\widetilde X(t)\mathsf E(t)
		-\bigl(n(t)^\top\nabla_{\R^3}\widetilde X(t)\bigr)\,\mathrm{I\!I}(t),
	\end{equation}
	where $\mathrm{I\!I}(t)$ is the matrix of $\mathrm{I\!I}_t$ in the frame and
	the Gaussian curvature is $\mathsf K(t)=\det\mathrm{I\!I}(t)$; for
	$\mathbb S^2$, $n(t)=t$ and $\mathrm{I\!I}(t)=I_2$, which recovers
	\eqref{eq:sphericalambienthessian}.
	
	Let $D_d =\diag(d_1,d_2,d_3)$ with $d_j>0$, and let
	$\mathcal M=\{t\in\R^3:t^\top D_d ^{-2}t=1\}$ carry the induced metric. Then
	\[
	n(t)=\frac{D_d ^{-2}t}{\norm{D_d ^{-2}t}},\qquad
	\mathrm{I\!I}(t)=\frac{\mathsf E^\top D_d ^{-2}\mathsf E}{\norm{D_d ^{-2}t}},\qquad
	\mathsf K(t)=\frac{1}{d_1^2d_2^2d_3^2\norm{D_d ^{-2}t}^4},
	\]
	so the curvature is not constant unless $d_1=d_2=d_3$; for
	$d=(1,1.5,2)$ it ranges over $[1/9,16/9]$. Let $\Xi$ be a symmetric
	$3\times3$ Gaussian matrix whose distinct entries
	$\zeta=(\Xi_{11},\Xi_{22},\Xi_{33},\Xi_{12},\Xi_{13},\Xi_{23})^\top$ are
	centered with covariance $\Sigma_\zeta\succ0$, and set $X(t)=t^\top\Xi t$ for
	$t\in\mathcal M$. By \eqref{eq:embeddedhessian} with
	$\widetilde X(t)=t^\top\Xi t$,
	\begin{equation}\label{eq:ellipsoidjet}
		\nabla X(t)=2\mathsf E^\top\Xi t,\qquad
		H(t)=2\mathsf E^\top\Xi\mathsf E
		-\frac{2\,t^\top D_d ^{-2}\Xi t}{\norm{D_d ^{-2}t}^2}\,\mathsf E^\top D_d ^{-2}\mathsf E,
	\end{equation}
	so the second fundamental form, and with it the curvature, enters the
	jet at every point. For $D_d =I_3$, the gradient and Hessian formulas reduce to the first two identities in \eqref{eq:sphericalquadrupolejet}. If $\Xi$ is also traceless, then $\tau=-6X$, recovering the third identity as well.
	
	\begin{proposition}[Quadratic form on an ellipsoid]
		\label{prop:ellipsoid}
		(i) Without constraints on $\Xi$, the linear jet map
		$\zeta\mapsto(X(t),\nabla X(t),\vech H(t))$ is invertible for every
		$t\in\mathcal M$. Hence the jet law is nonsingular,
		Assumption~\ref{ass:regularity} holds, $\Omega(t)\succ0$, and the local
		peak height density is positive on $\R$.
		
		(ii) Suppose instead that $\tr(D_d \Xi D_d )=0$, with the remaining five
		coefficients having a positive definite covariance. Then
		Assumption~\ref{ass:regularity} holds and the jet law is singular at
		every $t\in\mathcal M$. Moreover, $x_*=0$, and the constraint operator of
		Theorem~\ref{thm:singularsupport} is
		\begin{equation}\label{eq:ellipsoidQ}
			Q(t)=\tfrac16B_t^{-1}=\tfrac16\,\mathsf E^\top(D_d ^2-tt^\top)\mathsf E\succ0,
			\qquad B_t=\mathsf E^\top D_d ^{-2}\mathsf E=\norm{D_d ^{-2}t}\,\mathrm{I\!I}(t).
		\end{equation}
		Consequently, the local peak height density is positive exactly on
		$(0,\infty)$.
	\end{proposition}
	\begin{proof}
		The paths are restrictions of quadratic polynomials, which gives the
		regularity condition. Fix $t$. The jet at $t$ is linear in $\Xi$, so we
		evaluate it at a fixed symmetric matrix $\Xi$. Since
		$t^\top n(t)=\norm{D_d ^{-2}t}^{-1}>0$, we have $t\notin T_t\mathcal M$, and
		hence $\R^3=T_t\mathcal M\oplus\operatorname{span}\{t\}$. If $\Xi$ gives
		$\nabla X(t)=0$, then $\Xi t$ is normal, hence a multiple of $D_d ^{-2}t$;
		multiplying by $t^\top$ and using $t^\top D_d ^{-2}t=1$ gives
		$\Xi t=X(t)D_d ^{-2}t$. The coefficient of $\mathsf E^\top D_d ^{-2}\mathsf E$ in
		\eqref{eq:ellipsoidjet} then equals $2X(t)$, and
		$H(t)=2\mathsf E^\top\{\Xi-X(t)D_d ^{-2}\}\mathsf E$.
		
		(i) Let $\Xi$ lie in the kernel of the jet map, so that $X(t)=0$,
		$\nabla X(t)=0$ and $H(t)=0$. Then $\Xi t=0$ and $u^\top\Xi u'=0$ for
		$u,u'\in T_t\mathcal M$, so $\Xi=0$. Thus the jet map is injective, the
		jet covariance is positive definite, both conditions in
		\eqref{eq:non-degenerate} hold, and $\varepsilon_*>0$ gives $\Omega\succ0$
		by Theorem~\ref{thm:singularrank}.
		
		(ii) If $\Xi$ satisfies the constraint and gives $\nabla X(t)=0$ and
		$H(t)=0$, the same argument applied to $\Xi-X(t)D_d ^{-2}$ gives
		$\Xi=X(t)D_d ^{-2}$, and $\tr(D_d \Xi D_d )=3X(t)=0$ forces $\Xi=0$; hence
		$\Cov\bigl((\nabla X(t)^\top,W^\top)^\top\bigr)\succ0$. For $w\ne0$ with $w^\top t=0$, the matrix
		$\Xi=D_d ^{-2}-3ww^\top/\norm{D_d w}^2$ satisfies the constraint and gives
		$\nabla X(t)=0$ but $X(t)=1$, so $\Var(X(t)\mid\nabla X(t))>0$. Only five
		coefficients drive the six-dimensional jet, so the jet law is singular.
		Now let $\nabla X(t)=0$, and write $X=X(t)$ and $H=H(t)$. Since
		$\Xi t=XD_d^{-2}t$, in the basis $F=(e_1,e_2,t)$ the matrix $F^\top\Xi F$ is
		block diagonal with blocks $H/2+XB_t$ and $X$.
		Likewise, $F^\top D_d ^{-2}F=\diag(B_t,1)$, because $e_i^\top D_d ^{-2}t=0$ and
		$t^\top D_d ^{-2}t=1$, so that $F^{-1}D_d ^2F^{-\top}=\diag(B_t^{-1},1)$. Hence
		\[
		0=\tr(D_d \Xi D_d )=\tr\bigl(F^\top\Xi F\,F^{-1}D_d ^2F^{-\top}\bigr)
		=\tfrac12\tr(B_t^{-1}H)+X\{\tr(I_2)+1\}
		=\tfrac12\tr(B_t^{-1}H)+3X,
		\]
		which is $\tr(QH)=-X$ with $Q=B_t^{-1}/6\succ0$. By the uniqueness in
		Theorem~\ref{thm:singularrank}, $Q$ is the constraint operator and
		$x_*=0$. The rows of $F^{-1}$
		are $(e_i-(e_i^\top t)D_d ^{-2}t)^\top$, $i=1,2$, and $(D_d ^{-2}t)^\top$, so the
		upper-left block of $F^{-1}D_d ^2F^{-\top}$ is also
		$\mathsf E^\top(D_d ^2-tt^\top)\mathsf E$; this gives the second expression in
		\eqref{eq:ellipsoidQ}, and $B_t=\norm{D_d ^{-2}t}\,\mathrm{I\!I}(t)$ follows from
		the form of $\mathrm{I\!I}(t)$ above.
	\end{proof}
	
	By \eqref{eq:ellipsoidQ}, $Q(t)$ is proportional to $\mathrm{I\!I}(t)^{-1}$,
	so its eigenvalue ratio is the reciprocal of the ratio of the principal curvatures of
	$\mathcal M$ at $t$; at the vertex $(0,0,d_3)$, for instance,
	$Q=\diag(d_1^2,d_2^2)/6$. The constraint in (ii) is therefore definite,
	with equal eigenvalues exactly at the umbilic points of $\mathcal M$, of
	which there are four when the semi-axes are distinct. Elsewhere the trace
	is random in the original coordinates and deterministic in the canonical
	coordinates of Theorem~\ref{thm:canonicalweights}(i). For $D_d =I_3$,
	$B_t=I_2$ and \eqref{eq:ellipsoidQ} reduces to $Q=I_2/6$, as for the
	quadrupole of Section~\ref{subsec:sphericalquadrupole}.
	By Theorem~\ref{thm:singularpde}(ii), the field of case (ii) solves
	on all of $\mathcal M$ an elliptic equation with principal part
	$Q(t)=B_t^{-1}/6$, proportional to $\mathrm{I\!I}(t)^{-1}$.
	
	The exact check of Section~\ref{subsec:sphericalchecks} extends to this
	surface. Critical points satisfy $\Xi t=\alpha D_d ^{-2}t$, that is,
	$(D_d \Xi D_d )q=\alpha q$ for $q=D_d ^{-1}t\in \mathbb S^2$, and the height is $\alpha$.
	The map $q\mapsto D_d q$ is a diffeomorphism of $\mathbb S^2$ onto $\mathcal M$ that
	carries $X$ to the quadratic form $q^\top(D_d \Xi D_d )q$, so critical points and
	their types correspond. Hence, with $\alpha_{(1)}<\alpha_{(2)}<\alpha_{(3)}$
	the eigenvalues of $D_d \Xi D_d $, which are distinct almost surely, the field
	is almost surely Morse with two maxima at height $\alpha_{(3)}$, two
	saddles at height $\alpha_{(2)}$ and two minima at height $\alpha_{(1)}$,
	and \eqref{eq:sphericalquadrupolecountcheck} holds with $\mathbb S^2$ replaced by
	$\mathcal M$ and $p_{(i)}$ the density of $\alpha_{(i)}$. The local
	intensities, however, are computed in the metric of $\mathcal M$ through
	\eqref{eq:ellipsoidjet}, so this identity tests the curvature terms.
	Figure~\ref{fig:surfaces}(c) shows, for $d=(1,1.5,2)$ in case (i),
	the local peak height densities at the three vertices $d_ie_i$, where
	$X(d_ie_i)=d_i^2\Xi_{ii}$; the jet at each vertex combines the
	coefficient variances, scaled by the semi-axes, with the second
	fundamental form in \eqref{eq:ellipsoidjet}.

	\section{Discussion}
	\label{sec:discussion}
	For a smooth Gaussian field on a Riemannian surface, the peak
	intensity, the expected number of peaks above a level and the local peak
	height law reduce to the expected peak weight $\Kcal$ at the
	regression parameters of Section~\ref{sec:regression}; Sections
	\ref{sec:trace}--\ref{sec:continuity} evaluate $\Kcal$, classify its
	singular laws, and show that the resulting height laws are stable across
	degeneracy.
	Riemannian geometry enters through the covariant jet, through identities
	such as \eqref{eq:generaleigenfunctionconstraint}, and through global
	balances such as the Morse relation of Section~\ref{subsec:sphericalchecks}; the pointwise reduction and support
	classification are the same on every surface.
	
	On the computational side, the norm kernel is a finite combination of
	ordinary noncentral chi square CDFs for equal norm eigenvalues or a
	centered norm (Propositions~\ref{prop:ordinaryequal}
	and~\ref{prop:ordinarycentered}), and of generalized ones in general, so
	that conditional peak weights require at most two scalar integrations. For a nondegenerate Hessian law, Theorem~\ref{thm:spectral}
	reduces them to a single integral.
	Corollary~\ref{cor:fieldparametercontinuity} also has a statistical
	reading: estimates of the mean and covariance jets with uniform
	error $d_n$ yield plug-in estimates of local peak height laws and
	expected counts with error $O(d_n)$ in total variation and
	$O(d_n^{3/4})$ uniformly, where both exponents are optimal in general. Since the
	covariance jet involves derivatives of $C$ of total order four, the
	attainable rate of $d_n$ then limits the accuracy of such plug-in
	estimates in peak-based inference \citep{ChengSchwartzman2017}.
	
	The affine support and cone arguments extend to higher-dimensional
	tangent spaces. So does Theorem~\ref{thm:singularpde},
	with the conics of Corollary~\ref{cor:spectralconic} replaced by the
	quadrics $\{k^\top Qk=1\}$ in $\R^n$; a half-line support then
	corresponds to a semidefinite principal part. The trace reductions,
		however, use a feature specific to dimension two: the
	traceless Hessian has eigenvalues $\pm\norm\xi$, so the trace and
	squared traceless norm determine both the peak constraint and the
	determinant. In dimension $n\ge3$, negative definiteness instead
	requires
	\[
	\lambda_{\max}\!\left(
	\nabla^2X(t)-\frac{\tau}{n}I_n\right)<-\frac{\tau}{n},
	\qquad \tau=\tr\nabla^2X(t).
	\]
	The squared traceless norm does not determine this largest
	eigenvalue. For example, $\diag(1,-5,-5)$ and $\diag(-7,-1,-1)$
	both have trace $-9$ and squared traceless Frobenius norm $24$,
	but only the second is negative definite.
	Higher spectral invariants therefore enter the calculation, and
	the determinant is a polynomial of degree $n$, rather than the
	quadratic expression used here. Isotropic models admit random matrix reductions
	\citep{AzaisDelmas2022,KurikiEtAl2026}; general fields in higher
	dimensions may require conditioning on further spectral invariants or
	additional covariance symmetries. Degenerate Hessian laws given the
	gradient (Remark~\ref{rem:degeneratehessian}) likewise require a Morse
	criterion and a stability theory that allows the normalizer to vanish.

	\section*{Acknowledgments} This work was supported by the NSF project ``ATD: Collaborative Research: A Geostatistical
	Framework for Spatiotemporal Extremes''
	(DMS-2220529 and DMS-2220523).

	\appendix
	
	\section{Gaussian norm moments}
	\label{app:normmoments}
	This appendix records the moment calculations used to evaluate the
	kernels in Section~\ref{subsec:finiteCDF}.
	
	\subsection{Spectral representation and chi square moments}
	
	Write $\chi^2_k(\omega)$ for the noncentral chi square distribution with
	$k>0$ degrees of freedom and noncentrality $\omega\ge0$, and
	$F^{\chi^2}_{k,\omega}$ for its CDF. For
	$\xi\sim\mathsf N_2(\eta,\Lambda)$ with $\Lambda\succ0$, choose an
	orthogonal matrix $O$ such that
	$\Lambda=O\diag(\lambda_1,\lambda_2)O^\top$, $\lambda_1,\lambda_2>0$,
	and set $\widehat\eta=O^\top\eta$ and
	$\omega_j=\widehat\eta_j^2/\lambda_j$, $j=1,2$.
	Orthogonal transformations preserve the norm, so for independent
	standard normals $\mathcal N_1,\mathcal N_2$,
	\begin{equation}\label{eq:normspectral}
		\norm\xi^2\overset{\mathrm d}=
		\sum_{j=1}^2\lambda_j
		\left(\mathcal N_j+\frac{\widehat\eta_j}{\sqrt{\lambda_j}}\right)^2.
	\end{equation}
	Each squared term has law $\chi_1^2(\omega_j)$.
	The following density identity will be used to obtain
	the finite CDF formulas.
	
	For the scalar moments used in these formulas, let $f_{k,\omega}$
	be the density of $\chi_k^2(\omega)$. The Poisson mixture and
	the central identity $qf_{k,0}(q)=kf_{k+2,0}(q)$ give
	\begin{equation}\label{eq:chisqdensityshift}
		qf_{k,\omega}(q)
		=kf_{k+2,\omega}(q)+\omega f_{k+4,\omega}(q),\qquad q>0.
	\end{equation}
	Indeed, writing $p_\ell=e^{-\omega/2}(\omega/2)^\ell/\ell!$,
	one has $f_{k,\omega}=\sum_{\ell\ge0}p_\ell f_{k+2\ell,0}$
	and $2\ell p_\ell=\omega p_{\ell-1}$ for $\ell\ge1$.
	Expanding the left side of \eqref{eq:chisqdensityshift} and
	reindexing its $2\ell$ term proves the identity, including
	$\omega=0$.
	
	The following propositions use \eqref{eq:normspectral}
	with the notation above and $a,y\ge0$.
	
	\subsection{Circular norm laws}
	
	\begin{proposition}[Equal eigenvalues]
		\label{prop:ordinaryequal}
		If $\lambda_1=\lambda_2=\lambda$, let
		$\omega=\omega_1+\omega_2=\norm\eta^2/\lambda$. Then
		\begin{equation}\label{eq:circularordinary}
			F_{\eta,\lambda I_2}(y)
			=F^{\chi^2}_{2,\omega}(y/\lambda),
			\qquad
			\Acal(a;\eta,\lambda I_2)=aF^{\chi^2}_{2,\omega}(a/\lambda)
			-\lambda\left\{2F^{\chi^2}_{4,\omega}(a/\lambda)
			+\omega F^{\chi^2}_{6,\omega}(a/\lambda)\right\}.
		\end{equation}
		If also $\eta=0$, these formulas simplify to
		\eqref{eq:centeredcircularkernel}.
		
	\end{proposition}
	\begin{proof}
		Additivity of independent noncentral chi square variables gives
		$Y=\|\xi\|^2/\lambda\sim\chi_2^2(\omega)$ and the CDF in
		\eqref{eq:circularordinary}. Integrating
		\eqref{eq:chisqdensityshift} with $k=2$ over $0<q\le a/\lambda$
		gives its norm moment, since $\Acal(a;\eta,\lambda I_2)
		=a\PP(Y\le a/\lambda)
		-\lambda\E[Y\ind_{\{Y\le a/\lambda\}}]$.
		The exponential law of $\lambda\chi_2^2$ gives
		\eqref{eq:centeredcircularkernel}.
	\end{proof}
	
	\subsection{Centered norms with unequal eigenvalues}
	
	\begin{proposition}[Unequal eigenvalues with centered norm]
		\label{prop:ordinarycentered}
		Suppose $\lambda_1>\lambda_2>0$ and $\eta=0$. Let $\gamma_\pm=
		\bigl(\lambda_2^{-1/2}\pm\lambda_1^{-1/2}\bigr)/2$.
		Then, for $a,y\ge0$,
		\begin{equation}\label{eq:centeredordinarycdf}
			F_{0,\Lambda}(y)
			=F^{\chi^2}_{2,\gamma_-^2y}(\gamma_+^2y)
			-F^{\chi^2}_{2,\gamma_+^2y}(\gamma_-^2y),
		\end{equation}
		and
		\begin{equation}\label{eq:centeredordinarykernel}
			\begin{split}
				\Acal(a;0,\Lambda)
				={}&(a-\lambda_1-\lambda_2)
				\left[
				F^{\chi^2}_{2,\gamma_-^2a}(\gamma_+^2a)
				-F^{\chi^2}_{2,\gamma_+^2a}(\gamma_-^2a)
				\right]\\
				&+\frac{a(\lambda_1+\lambda_2)}
				{2\sqrt{\lambda_1\lambda_2}}
				\left[
				1-F^{\chi^2}_{2,\gamma_-^2a}(\gamma_+^2a)
				-F^{\chi^2}_{2,\gamma_+^2a}(\gamma_-^2a)
				\right]\\
				&+\frac{a(\lambda_1+\lambda_2
					-2\sqrt{\lambda_1\lambda_2})}
				{2\sqrt{\lambda_1\lambda_2}}
				\left[
				F^{\chi^2}_{2,\gamma_-^2a}(\gamma_+^2a)
				-F^{\chi^2}_{4,\gamma_-^2a}(\gamma_+^2a)
				\right].
			\end{split}
		\end{equation}
	\end{proposition}
	
	The CDF \eqref{eq:centeredordinarycdf} is the squared Hoyt CDF of
	\citep{Paris2009}: with $F^{\chi^2}_{2,\omega}(x)=1-\mathrm Q_1(\sqrt\omega,\sqrt x)$
	for the Marcum function $\mathrm Q_1$, it equals
	$\mathrm Q_1(\gamma_+\sqrt y,\gamma_-\sqrt y)-\mathrm Q_1(\gamma_-\sqrt y,\gamma_+\sqrt y)$.
	
	\begin{proof}
		Write $F(y)=F_{0,\Lambda}(y)$, $p_\lambda=\sqrt{\lambda_1\lambda_2}$,
		and $s_\lambda=\lambda_1+\lambda_2$, and set
		\[
		\begin{split}
			E_0(y)=1-F^{\chi^2}_{2,\gamma_-^2y}(\gamma_+^2y)
			-F^{\chi^2}_{2,\gamma_+^2y}(\gamma_-^2y),\qquad
			E_1(y)=\frac{\gamma_-}{\gamma_+}
			\left[F^{\chi^2}_{2,\gamma_-^2y}(\gamma_+^2y)
			-F^{\chi^2}_{4,\gamma_-^2y}(\gamma_+^2y)\right].
		\end{split}
		\]
		Differentiation of the ordinary noncentral chi square CDFs gives
		\[
		F'(y)=\frac{E_0(y)}{2p_\lambda},\qquad
		\frac{\dd}{\dd y}
		\left\{y[s_\lambda E_0(y)+(\lambda_1-\lambda_2)E_1(y)]\right\}
		=(s_\lambda-y)E_0(y).
		\]
		Let $R(y)$ denote the right-hand side of
		\eqref{eq:centeredordinarykernel}, with $a=y$.
		Since $(\lambda_1-\lambda_2)\gamma_-/\gamma_+=s_\lambda-2p_\lambda$,
		\[
		R(y)=(y-s_\lambda)F(y)
		+\frac{y}{2p_\lambda}
		[s_\lambda E_0(y)+(\lambda_1-\lambda_2)E_1(y)].
		\]
		Thus $R(0)=0$ and $R'(y)=F(y)$. Then, \eqref{eq:Alayer}
		yields $\Acal(a;0,\Lambda)=R(a)$.
	\end{proof}
	
	\subsection{Scalar moments and noncentered norms}
	Recall the scalar kernel $\mathcal H$ defined in \eqref{eq:kernelH}.
	For $\lambda>0$, let $\omega=\mu_0^2/\lambda$. Then
	$(\mu_0+\sqrt\lambda\,\mathcal N)^2/\lambda\sim\chi_1^2(\omega)$.
	Integrating \eqref{eq:chisqdensityshift} with $k=1$ gives, for $b>0$,
	\begin{equation}\label{eq:Hnoncentralcdf}
		\mathcal H(b;\mu_0,\lambda)
		=bF^{\chi^2}_{1,\omega}(b/\lambda)
		-\lambda\left\{F^{\chi^2}_{3,\omega}(b/\lambda)
		+\omega F^{\chi^2}_{5,\omega}(b/\lambda)\right\}.
	\end{equation}
	For a centered scalar Gaussian variable, only two central chi square
	CDFs are needed:
	\begin{equation}\label{eq:Hcentralcdf}
		\mathcal H(b;0,\lambda)
		=bF^{\chi^2}_{1,0}(b/\lambda)
		-\lambda F^{\chi^2}_{3,0}(b/\lambda),\qquad b>0.
	\end{equation}
	\begin{remark}[Noncentered norms with unequal variances]
		\label{remark:ordinarynoncentered}
		If $\eta\ne0$ and $\Lambda\succ0$ has unequal eigenvalues,
		\eqref{eq:normconvolution} with \eqref{eq:Hnoncentralcdf} evaluates
		$\Acal$ by one scalar integral of ordinary noncentral chi square CDFs.
		Alternatively, let $\omega_j=\widehat\eta_j^2/\lambda_j$ and
		$F^{\mathrm{gen}}_{k_1,k_2}(y)=\PP\{\lambda_1Y_1+\lambda_2Y_2\le y\}$
		for independent $Y_j\sim\chi^2_{k_j}(\omega_j)$, so that
		$F^{\mathrm{gen}}_{1,1}=F_{\eta,\Lambda}$. Applying
		\eqref{eq:chisqdensityshift} to each coordinate gives
		\begin{equation}\label{eq:normgeneralizedcdf}
			\Acal(a;\eta,\Lambda)
			=aF^{\mathrm{gen}}_{1,1}(a)
			-\lambda_1\{F^{\mathrm{gen}}_{3,1}(a)+\omega_1F^{\mathrm{gen}}_{5,1}(a)\}
			-\lambda_2\{F^{\mathrm{gen}}_{1,3}(a)+\omega_2F^{\mathrm{gen}}_{1,5}(a)\},
		\end{equation}
		which requires generalized chi square CDFs; Ruben's expansion
		\citep{Ruben1962} is a classical tool for these.
	\end{remark}
	
	\subsection{Singular norm covariances}
	
	\begin{proposition}[Norm kernel with a singular covariance]
		\label{prop:singularcdfnorm}
		Let $\xi\sim\mathsf N_2(\eta,\Lambda)$ with $\Lambda\succeq0$ singular.
		
		(i) If $\operatorname{rank}(\Lambda)=1$, choose an orthogonal
		matrix $O$ such that $\Lambda=O\diag(\lambda,0)O^\top$, $\lambda>0$,
		and set $\widehat\eta=O^\top\eta$,
		$d_0=\widehat\eta_2^2$, and $\omega=\widehat\eta_1^2/\lambda$.
		Then
		\begin{equation}\label{eq:singularnormcdf}
			F_{\eta,\Lambda}(y)=
			\begin{cases}
				F^{\chi^2}_{1,\omega}((y-d_0)/\lambda),&y>d_0,\\
				0,&y\le d_0,
			\end{cases}
		\end{equation}
		and the norm kernel is simply
		\begin{equation}\label{eq:Hrankonereduction}
			\Acal(a;\eta,\Lambda)
			=\mathcal H(a-d_0;\widehat\eta_1,\lambda),\qquad a\ge0.
		\end{equation}
		The convention in \eqref{eq:kernelH} covers $a\le d_0$.
		
		(ii) If $\Lambda=0$, then $\|\xi\|^2=\|\eta\|^2$ almost surely,
		$F_{\eta,0}(y)=\ind_{\{y\ge\|\eta\|^2\}}$, and
		$\Acal(a;\eta,0)=(a-\|\eta\|^2)_+$.
	\end{proposition}
	\begin{proof}
		In (i), orthogonal transformation gives
		\[
		\|\xi\|^2\overset{\mathrm d}=
		d_0+(\widehat\eta_1+\sqrt\lambda\,\mathcal N)^2,
		\qquad \mathcal N\sim\mathsf N(0,1).
		\]
		The ordinary noncentral chi square law of the random square gives
		\eqref{eq:singularnormcdf}, and \eqref{eq:kernelH} gives
		\eqref{eq:Hrankonereduction} directly. In (ii), $\xi=\eta$
		almost surely.
	\end{proof}

\end{document}